\documentclass[11pt,a4paper]{article}

\usepackage{amsmath,amssymb,latexsym,pstricks,mathrsfs,comment,amsthm,graphicx,tikz,tikz-cd,enumerate,accents,pgffor,cite,wrapfig,multicol,float}
\usepackage{bibspacing}
\newcommand{\nc}{\newcommand}
\nc{\rnc}{\renewcommand}

\usepackage{geometry}  \rnc{\ss}{\smallskip} \nc{\ms}{\medskip} \nc{\bs}{\bigskip} \nc{\nss}{\vspace{-3mm}}

\makeatletter

\DeclareMathSymbol{\widehatsym}{\mathord}{largesymbols}{"62}
\newcommand\lowerwidehatsym{%
  \text{\smash{\raisebox{-1.3ex}{%
    $\widehatsym$}}}}
\newcommand\fixwidehat[1]{%
  \mathchoice
    {\accentset{\displaystyle\lowerwidehatsym}{#1}}
    {\accentset{\textstyle\lowerwidehatsym}{#1}}
    {\accentset{\scriptstyle\lowerwidehatsym}{#1}}
    {\accentset{\scriptscriptstyle\lowerwidehatsym}{#1}}
}
\rnc{\widehat}{\fixwidehat}

\begin{document}

\nc{\M}{\mathcal M}
\nc{\G}{\mathcal G}
\nc{\F}{\mathbb F}
\nc{\MnJ}{\mathcal M_n^J}
\nc{\EnJ}{\mathcal E_n^J}
\nc{\Mat}{\operatorname{Mat}}
\nc{\RegMnJ}{\Reg(\MnJ)}
\nc{\row}{\mathfrak r}
\nc{\col}{\mathfrak c}
\nc{\Row}{\operatorname{Row}}
\nc{\Col}{\operatorname{Col}}
\nc{\Span}{\operatorname{span}}
\nc{\mat}[4]{\left[\begin{matrix}#1&#2\\#3&#4\end{matrix}\right]}
\nc{\tmat}[4]{\left[\begin{smallmatrix}#1&#2\\#3&#4\end{smallmatrix}\right]}
\nc{\ttmat}[4]{{\tiny \left[\begin{smallmatrix}#1&#2\\#3&#4\end{smallmatrix}\right]}}
\nc{\tmatt}[9]{\left[\begin{smallmatrix}#1&#2&#3\\#4&#5&#6\\#7&#8&#9\end{smallmatrix}\right]}
\nc{\ttmatt}[9]{{\tiny \left[\begin{smallmatrix}#1&#2&#3\\#4&#5&#6\\#7&#8&#9\end{smallmatrix}\right]}}
\nc{\MnGn}{\M_n\sm\G_n}
\nc{\MrGr}{\M_r\sm\G_r}
\nc{\qbin}[2]{\left[\begin{matrix}#1\\#2\end{matrix}\right]_q}
\nc{\tqbin}[2]{\left[\begin{smallmatrix}#1\\#2\end{smallmatrix}\right]_q}
\nc{\qbinx}[3]{\left[\begin{matrix}#1\\#2\end{matrix}\right]_{#3}}
\nc{\tqbinx}[3]{\left[\begin{smallmatrix}#1\\#2\end{smallmatrix}\right]_{#3}}
\nc{\MNJ}{\M_nJ}
\nc{\JMN}{J\M_n}
\nc{\RegMNJ}{\Reg(\MNJ)}
\nc{\RegJMN}{\Reg(\JMN)}
\nc{\RegMMNJ}{\Reg(\MMNJ)}
\nc{\RegJMMN}{\Reg(\JMMN)}
\nc{\Wb}{\overline{W}}
\nc{\Xb}{\overline{X}}
\nc{\Yb}{\overline{Y}}
\nc{\Zb}{\overline{Z}}
\nc{\Sib}{\overline{\Si}}
\nc{\Om}{\Omega}
\nc{\Omb}{\overline{\Om}}
\nc{\Gab}{\overline{\Ga}}
\nc{\qfact}[1]{[#1]_q!}
\nc{\smat}[2]{\left[\begin{matrix}#1&#2\end{matrix}\right]}
\nc{\tsmat}[2]{\left[\begin{smallmatrix}#1&#2\end{smallmatrix}\right]}
\nc{\hmat}[2]{\left[\begin{matrix}#1\\#2\end{matrix}\right]}
\nc{\thmat}[2]{\left[\begin{smallmatrix}#1\\#2\end{smallmatrix}\right]}
\nc{\LVW}{\mathcal L(V,W)}
\nc{\KVW}{\mathcal K(V,W)}
\nc{\LV}{\mathcal L(V)}
\nc{\RegLVW}{\Reg(\LVW)}
\nc{\sM}{\mathscr M}
\nc{\sN}{\mathscr N}
\rnc{\iff}{\Leftrightarrow}
\nc{\Hom}{\operatorname{Hom}}
\nc{\End}{\operatorname{End}}
\nc{\Aut}{\operatorname{Aut}}
\nc{\Lin}{\mathcal L}
\nc{\Hommn}{\Hom(V_m,V_n)}
\nc{\Homnm}{\Hom(V_n,V_m)}
\nc{\Homnl}{\Hom(V_n,V_l)}
\nc{\Homkm}{\Hom(V_k,V_m)}
\nc{\Endm}{\End(V_m)}
\nc{\Endn}{\End(V_n)}
\nc{\Endr}{\End(V_r)}
\nc{\Autm}{\Aut(V_m)}
\nc{\Autn}{\Aut(V_n)}
\nc{\MmnJ}{\M_{mn}^J}
\nc{\MmnA}{\M_{mn}^A}
\nc{\MmnB}{\M_{mn}^B}
\nc{\Mmn}{\M_{mn}}
\nc{\Mkl}{\M_{kl}}
\nc{\Mnm}{\M_{nm}}
\nc{\EmnJ}{\mathcal E_{mn}^J}
\nc{\MmGm}{\M_m\sm\G_m}
\nc{\RegMmnJ}{\Reg(\MmnJ)}
\rnc{\implies}{\Rightarrow}
\nc{\DMmn}[1]{D_{#1}(\Mmn)}
\nc{\DMmnJ}[1]{D_{#1}(\MmnJ)}
\nc{\MMNJ}{\Mmn J}
\nc{\JMMN}{J\Mmn}
\nc{\JMMNJ}{J\Mmn J}
\nc{\Inr}{\mathcal I(V_n,W_r)}
\nc{\Lnr}{\mathcal L(V_n,W_r)}
\nc{\Knr}{\mathcal K(V_n,W_r)}
\nc{\Imr}{\mathcal I(V_m,W_r)}
\nc{\Kmr}{\mathcal K(V_m,W_r)}
\nc{\Lmr}{\mathcal L(V_m,W_r)}
\nc{\Kmmr}{\mathcal K(V_m,W_{m-r})}
\nc{\tr}{{\operatorname{T}}}
\nc{\MMN}{\MmnA(\F_1)}
\nc{\MKL}{\Mkl^B(\F_2)}
\nc{\RegMMN}{\Reg(\MmnA(\F_1))}
\nc{\RegMKL}{\Reg(\Mkl^B(\F_2))}
\nc{\gRhA}{\widehat{\mathscr R}^A}
\nc{\gRhB}{\widehat{\mathscr R}^B}
\nc{\gLhA}{\widehat{\mathscr L}^A}
\nc{\gLhB}{\widehat{\mathscr L}^B}
\nc{\timplies}{\Rightarrow}
\nc{\tiff}{\Leftrightarrow}
\nc{\Sija}{S_{ij}^a}
\nc{\dmat}[8]{\draw(#1*1.5,#2)node{$\left[\begin{smallmatrix}#3&#4&#5\\#6&#7&#8\end{smallmatrix}\right]$};}
\nc{\bdmat}[8]{\draw(#1*1.5,#2)node{${\mathbf{\left[\begin{smallmatrix}#3&#4&#5\\#6&#7&#8\end{smallmatrix}\right]}}$};}
\nc{\rdmat}[8]{\draw(#1*1.5,#2)node{\rotatebox{90}{$\left[\begin{smallmatrix}#3&#4&#5\\#6&#7&#8\end{smallmatrix}\right]$}};}
\nc{\rldmat}[8]{\draw(#1*1.5-0.375,#2)node{\rotatebox{90}{$\left[\begin{smallmatrix}#3&#4&#5\\#6&#7&#8\end{smallmatrix}\right]$}};}
\nc{\rrdmat}[8]{\draw(#1*1.5+.375,#2)node{\rotatebox{90}{$\left[\begin{smallmatrix}#3&#4&#5\\#6&#7&#8\end{smallmatrix}\right]$}};}
\nc{\rfldmat}[8]{\draw(#1*1.5-0.375+.15,#2)node{\rotatebox{90}{$\left[\begin{smallmatrix}#3&#4&#5\\#6&#7&#8\end{smallmatrix}\right]$}};}
\nc{\rfrdmat}[8]{\draw(#1*1.5+.375-.15,#2)node{\rotatebox{90}{$\left[\begin{smallmatrix}#3&#4&#5\\#6&#7&#8\end{smallmatrix}\right]$}};}
\nc{\xL}{[x]_{\! _\gL}}\nc{\yL}{[y]_{\! _\gL}}\nc{\xR}{[x]_{\! _\gR}}\nc{\yR}{[y]_{\! _\gR}}\nc{\xH}{[x]_{\! _\gH}}\nc{\yH}{[y]_{\! _\gH}}\nc{\XK}{[X]_{\! _\gK}}\nc{\xK}{[x]_{\! _\gK}}
\nc{\RegSija}{\Reg(\Sija)}
\nc{\MnmK}{\M_{nm}^K}
\nc{\cC}{\mathcal C}
\nc{\cR}{\mathcal R}
\nc{\Ckl}{\cC_k(l)}
\nc{\Rkl}{\cR_k(l)}
\nc{\Cmr}{\cC_m(r)}
\nc{\Rmr}{\cR_m(r)}
\nc{\Cnr}{\cC_n(r)}
\nc{\Rnr}{\cR_n(r)}
\nc{\Z}{\mathbb Z}

\nc{\Reg}{\operatorname{Reg}}
\nc{\RP}{\operatorname{RP}}
\nc{\TXa}{\T_X^a}
\nc{\TXA}{\T(X,A)}
\nc{\TXal}{\T(X,\al)}
\nc{\RegTXa}{\Reg(\TXa)}
\nc{\RegTXA}{\Reg(\TXA)}
\nc{\RegTXal}{\Reg(\TXal)}
\nc{\PalX}{\P_\al(X)}
\nc{\EAX}{\E_A(X)}
\nc{\Bb}{\overline{B}}
\nc{\bb}{\overline{\be}}
\nc{\bw}{{\bf w}}
\nc{\bz}{{\bf z}}
\nc{\TASA}{\T_A\sm\S_A}
\nc{\Ub}{\overline{U}}
\nc{\Vb}{\overline{V}}
\nc{\eb}{\overline{e}}
\nc{\EXa}{\E_X^a}
\nc{\oijr}{1\leq i<j\leq r}
\nc{\veb}{\overline{\ve}}
\nc{\bbT}{\mathbb T}
\nc{\Surj}{\operatorname{Surj}}
\nc{\Sone}{S^{(1)}}
\nc{\fillbox}[2]{\draw[fill=gray!30](#1,#2)--(#1+1,#2)--(#1+1,#2+1)--(#1,#2+1)--(#1,#2);}
\nc{\raa}{\rangle_J}
\nc{\raJ}{\rangle_J}
\nc{\Ea}{E_J}
\nc{\EJ}{E_J}
\nc{\ep}{\epsilon} \nc{\ve}{\varepsilon}
\nc{\IXa}{\I_X^a}
\nc{\RegIXa}{\Reg(\IXa)}
\nc{\JXa}{\J_X^a}
\nc{\RegJXa}{\Reg(\JXa)}
\nc{\IXA}{\I(X,A)}
\nc{\IAX}{\I(A,X)}
\nc{\RegIXA}{\Reg(\IXA)}
\nc{\RegIAX}{\Reg(\IAX)}
\nc{\trans}[2]{\left(\begin{smallmatrix} #1 \\ #2 \end{smallmatrix}\right)}
\nc{\bigtrans}[2]{\left(\begin{matrix} #1 \\ #2 \end{matrix}\right)}
\nc{\lmap}[1]{\mapstochar \xrightarrow {\ #1\ }}
\nc{\EaTXa}{E}

\nc{\gL}{\mathscr L}
\nc{\gR}{\mathscr R}
\nc{\gH}{\mathscr H}
\nc{\gJ}{\mathscr J}
\nc{\gD}{\mathscr D}
\nc{\gK}{\mathscr K}
\nc{\gLa}{\mathscr L^a}
\nc{\gRa}{\mathscr R^a}
\nc{\gHa}{\mathscr H^a}
\nc{\gJa}{\mathscr J^a}
\nc{\gDa}{\mathscr D^a}
\nc{\gKa}{\mathscr K^a}
\nc{\gLJ}{\mathscr L^J}
\nc{\gRJ}{\mathscr R^J}
\nc{\gHJ}{\mathscr H^J}
\nc{\gJJ}{\mathscr J^J}
\nc{\gDJ}{\mathscr D^J}
\nc{\gKJ}{\mathscr K^J}
\nc{\gLh}{\widehat{\mathscr L}^J}
\nc{\gRh}{\widehat{\mathscr R}^J}
\nc{\gHh}{\widehat{\mathscr H}^J}
\nc{\gJh}{\widehat{\mathscr J}^J}
\nc{\gDh}{\widehat{\mathscr D}^J}
\nc{\gKh}{\widehat{\mathscr K}^J}
\nc{\Lh}{\widehat{L}^J}
\nc{\Rh}{\widehat{R}^J}
\nc{\Hh}{\widehat{H}^J}
\nc{\Jh}{\widehat{J}^J}
\nc{\Dh}{\widehat{D}^J}
\nc{\Kh}{\widehat{K}^J}
\nc{\gLb}{\widehat{\mathscr L}}
\nc{\gRb}{\widehat{\mathscr R}}
\nc{\gHb}{\widehat{\mathscr H}}
\nc{\gJb}{\widehat{\mathscr J}}
\nc{\gDb}{\widehat{\mathscr D}}
\nc{\gKb}{\widehat{\mathscr K}}
\nc{\Lb}{\widehat{L}^J}
\nc{\Rb}{\widehat{R}^J}
\nc{\Hb}{\widehat{H}^J}
\nc{\Jb}{\widehat{J}^J}
\nc{\Db}{\overline{D}}
\nc{\Kb}{\widehat{K}}

\hyphenation{mon-oid mon-oids}

\nc{\itemit}[1]{\item[\emph{(#1)}]}
\nc{\E}{\mathcal E}
\nc{\TX}{\T(X)}
\nc{\TXP}{\T(X,\P)}
\nc{\EX}{\E(X)}
\nc{\EXP}{\E(X,\P)}
\nc{\SX}{\S(X)}
\nc{\SXP}{\S(X,\P)}
\nc{\Sing}{\operatorname{Sing}}
\nc{\idrank}{\operatorname{idrank}}
\nc{\SingXP}{\Sing(X,\P)}
\nc{\De}{\Delta}
\nc{\sgp}{\operatorname{sgp}}
\nc{\mon}{\operatorname{mon}}
\nc{\Dn}{\mathcal D_n}
\nc{\Dm}{\mathcal D_m}

\nc{\lline}[1]{\draw(3*#1,0)--(3*#1+2,0);}
\nc{\uline}[1]{\draw(3*#1,5)--(3*#1+2,5);}
\nc{\thickline}[2]{\draw(3*#1,5)--(3*#2,0); \draw(3*#1+2,5)--(3*#2+2,0) ;}
\nc{\thicklabel}[3]{\draw(3*#1+1+3*#2*0.15-3*#1*0.15,4.25)node{{\tiny $#3$}};}

\nc{\slline}[3]{\draw(3*#1+#3,0+#2)--(3*#1+2+#3,0+#2);}
\nc{\suline}[3]{\draw(3*#1+#3,5+#2)--(3*#1+2+#3,5+#2);}
\nc{\sthickline}[4]{\draw(3*#1+#4,5+#3)--(3*#2+#4,0+#3); \draw(3*#1+2+#4,5+#3)--(3*#2+2+#4,0+#3) ;}
\nc{\sthicklabel}[5]{\draw(3*#1+1+3*#2*0.15-3*#1*0.15+#5,4.25+#4)node{{\tiny $#3$}};}

\nc{\stll}[5]{\sthickline{#1}{#2}{#4}{#5} \sthicklabel{#1}{#2}{#3}{#4}{#5}}
\nc{\tll}[3]{\stll{#1}{#2}{#3}00}

\nc{\mfourpic}[9]{
\slline1{#9}0
\slline3{#9}0
\slline4{#9}0
\slline5{#9}0
\suline1{#9}0
\suline3{#9}0
\suline4{#9}0
\suline5{#9}0
\stll1{#1}{#5}{#9}{0}
\stll3{#2}{#6}{#9}{0}
\stll4{#3}{#7}{#9}{0}
\stll5{#4}{#8}{#9}{0}
\draw[dotted](6,0+#9)--(8,0+#9);
\draw[dotted](6,5+#9)--(8,5+#9);
}
\nc{\vdotted}[1]{
\draw[dotted](3*#1,10)--(3*#1,15);
\draw[dotted](3*#1+2,10)--(3*#1+2,15);
}

\nc{\Clab}[2]{
\sthicklabel{#1}{#1}{{}_{\phantom{#1}}C_{#1}}{1.25+5*#2}0
}
\nc{\sClab}[3]{
\sthicklabel{#1}{#1}{{}_{\phantom{#1}}C_{#1}}{1.25+5*#2}{#3}
}
\nc{\Clabl}[3]{
\sthicklabel{#1}{#1}{{}_{\phantom{#3}}C_{#3}}{1.25+5*#2}0
}
\nc{\sClabl}[4]{
\sthicklabel{#1}{#1}{{}_{\phantom{#4}}C_{#4}}{1.25+5*#2}{#3}
}
\nc{\Clabll}[3]{
\sthicklabel{#1}{#1}{C_{#3}}{1.25+5*#2}0
}
\nc{\sClabll}[4]{
\sthicklabel{#1}{#1}{C_{#3}}{1.25+5*#2}{#3}
}

\nc{\mtwopic}[6]{
\slline1{#6*5}{#5}
\slline2{#6*5}{#5}
\suline1{#6*5}{#5}
\suline2{#6*5}{#5}
\stll1{#1}{#3}{#6*5}{#5}
\stll2{#2}{#4}{#6*5}{#5}
}
\nc{\mtwopicl}[6]{
\slline1{#6*5}{#5}
\slline2{#6*5}{#5}
\suline1{#6*5}{#5}
\suline2{#6*5}{#5}
\stll1{#1}{#3}{#6*5}{#5}
\stll2{#2}{#4}{#6*5}{#5}
\sClabl1{#6}{#5}{i}
\sClabl2{#6}{#5}{j}
}

\nc{\keru}{\operatorname{ker}^\wedge} \nc{\kerl}{\operatorname{ker}_\vee}

\nc{\coker}{\operatorname{coker}}
\nc{\KER}{\ker}
\nc{\N}{\mathbb N}
\nc{\LaBn}{L_\al(\B_n)}
\nc{\RaBn}{R_\al(\B_n)}
\nc{\LaPBn}{L_\al(\PB_n)}
\nc{\RaPBn}{R_\al(\PB_n)}
\nc{\rhorBn}{\rho_r(\B_n)}
\nc{\DrBn}{D_r(\B_n)}
\nc{\DrPn}{D_r(\P_n)}
\nc{\DrPBn}{D_r(\PB_n)}
\nc{\DrKn}{D_r(\K_n)}
\nc{\alb}{\al_{\vee}}
\nc{\beb}{\be^{\wedge}}
\nc{\bnf}{\bn^\flat}
\nc{\Bal}{\operatorname{Bal}}
\nc{\Red}{\operatorname{Red}}
\nc{\Pnxi}{\P_n^\xi}
\nc{\Bnxi}{\B_n^\xi}
\nc{\PBnxi}{\PB_n^\xi}
\nc{\Knxi}{\K_n^\xi}
\nc{\C}{\mathscr C}
\nc{\exi}{e^\xi}
\nc{\Exi}{E^\xi}
\nc{\eximu}{e^\xi_\mu}
\nc{\Eximu}{E^\xi_\mu}
\nc{\REF}{ {\red [Ref?]} }
\nc{\GL}{\operatorname{GL}}
\rnc{\O}{\operatorname{O}}

\nc{\vtx}[2]{\fill (#1,#2)circle(.2);}
\nc{\lvtx}[2]{\fill (#1,0)circle(.2);}
\nc{\uvtx}[2]{\fill (#1,1.5)circle(.2);}

\nc{\Eq}{\mathfrak{Eq}}
\nc{\Gau}{\Ga^\wedge} \nc{\Gal}{\Ga_\vee}
\nc{\Lamu}{\Lam^\wedge} \nc{\Laml}{\Lam_\vee}
\nc{\bX}{{\bf X}}
\nc{\bY}{{\bf Y}}
\nc{\ds}{\displaystyle}

\nc{\uvert}[1]{\fill (#1,1.5)circle(.2);}
\nc{\uuvert}[1]{\fill (#1,3)circle(.2);}
\nc{\uuuvert}[1]{\fill (#1,4.5)circle(.2);}
\rnc{\lvert}[1]{\fill (#1,0)circle(.2);}
\nc{\overt}[1]{\fill (#1,0)circle(.1);}
\nc{\overtl}[3]{\node[vertex] (#3) at (#1,0) {  {\tiny $#2$} };}
\nc{\cv}[2]{\draw(#1,1.5) to [out=270,in=90] (#2,0);}
\nc{\cvs}[2]{\draw(#1,1.5) to [out=270+30,in=90+30] (#2,0);}
\nc{\ucv}[2]{\draw(#1,3) to [out=270,in=90] (#2,1.5);}
\nc{\uucv}[2]{\draw(#1,4.5) to [out=270,in=90] (#2,3);}
\nc{\textpartn}[1]{{\lower0.45 ex\hbox{\begin{tikzpicture}[xscale=.2,yscale=0.2] #1 \end{tikzpicture}}}}
\nc{\textpartnx}[2]{{\lower1.0 ex\hbox{\begin{tikzpicture}[xscale=.3,yscale=0.3] 
\foreach \x in {1,...,#1}
{ \uvert{\x} \lvert{\x} }
#2 \end{tikzpicture}}}}
\nc{\disppartnx}[2]{{\lower1.0 ex\hbox{\begin{tikzpicture}[scale=0.3] 
\foreach \x in {1,...,#1}
{ \uvert{\x} \lvert{\x} }
#2 \end{tikzpicture}}}}
\nc{\disppartnxd}[2]{{\lower2.1 ex\hbox{\begin{tikzpicture}[scale=0.3] 
\foreach \x in {1,...,#1}
{ \uuvert{\x} \uvert{\x} \lvert{\x} }
#2 \end{tikzpicture}}}}
\nc{\disppartnxdn}[2]{{\lower2.1 ex\hbox{\begin{tikzpicture}[scale=0.3] 
\foreach \x in {1,...,#1}
{ \uuvert{\x} \lvert{\x} }
#2 \end{tikzpicture}}}}
\nc{\disppartnxdd}[2]{{\lower3.6 ex\hbox{\begin{tikzpicture}[scale=0.3] 
\foreach \x in {1,...,#1}
{ \uuuvert{\x} \uuvert{\x} \uvert{\x} \lvert{\x} }
#2 \end{tikzpicture}}}}

\nc{\dispgax}[2]{{\lower0.0 ex\hbox{\begin{tikzpicture}[scale=0.3] 
#2
\foreach \x in {1,...,#1}
{\lvert{\x} }
 \end{tikzpicture}}}}
\nc{\textgax}[2]{{\lower0.4 ex\hbox{\begin{tikzpicture}[scale=0.3] 
#2
\foreach \x in {1,...,#1}
{\lvert{\x} }
 \end{tikzpicture}}}}
\nc{\textlinegraph}[2]{{\raise#1 ex\hbox{\begin{tikzpicture}[scale=0.8] 
#2
 \end{tikzpicture}}}}
\nc{\textlinegraphl}[2]{{\raise#1 ex\hbox{\begin{tikzpicture}[scale=0.8] 
\tikzstyle{vertex}=[circle,draw=black, fill=white, inner sep = 0.07cm]
#2
 \end{tikzpicture}}}}
\nc{\displinegraph}[1]{{\lower0.0 ex\hbox{\begin{tikzpicture}[scale=0.6] 
#1
 \end{tikzpicture}}}}
 
\nc{\disppartnthreeone}[1]{{\lower1.0 ex\hbox{\begin{tikzpicture}[scale=0.3] 
\foreach \x in {1,2,3,5,6}
{ \uvert{\x} }
\foreach \x in {1,2,4,5,6}
{ \lvert{\x} }
\draw[dotted] (3.5,1.5)--(4.5,1.5);
\draw[dotted] (2.5,0)--(3.5,0);
#1 \end{tikzpicture}}}}

\nc{\partn}[4]{\left( \begin{array}{c|c} 
#1 \ & \ #3 \ \ \\ \cline{2-2}
#2 \ & \ #4 \ \
\end{array} \!\!\! \right)}
\nc{\partnlong}[6]{\partn{#1}{#2}{#3,\ #4}{#5,\ #6}} 
\nc{\partnsh}[2]{\left( \begin{array}{c} 
#1 \\
#2 
\end{array} \right)}
\nc{\partncodefz}[3]{\partn{#1}{#2}{#3}{\emptyset}}
\nc{\partndefz}[3]{{\partn{#1}{#2}{\emptyset}{#3}}}
\nc{\partnlast}[2]{\left( \begin{array}{c|c}
#1 \ &  \ #2 \\
#1 \ &  \ #2
\end{array} \right)}

\nc{\uuarcx}[3]{\draw(#1,3)arc(180:270:#3) (#1+#3,3-#3)--(#2-#3,3-#3) (#2-#3,3-#3) arc(270:360:#3);}
\nc{\uuarc}[2]{\uuarcx{#1}{#2}{.4}}
\nc{\uuuarcx}[3]{\draw(#1,4.5)arc(180:270:#3) (#1+#3,4.5-#3)--(#2-#3,4.5-#3) (#2-#3,4.5-#3) arc(270:360:#3);}
\nc{\uuuarc}[2]{\uuuarcx{#1}{#2}{.4}}
\nc{\darcx}[3]{\draw(#1,0)arc(180:90:#3) (#1+#3,#3)--(#2-#3,#3) (#2-#3,#3) arc(90:0:#3);}
\nc{\darc}[2]{\darcx{#1}{#2}{.4}}
\nc{\udarcx}[3]{\draw(#1,1.5)arc(180:90:#3) (#1+#3,1.5+#3)--(#2-#3,1.5+#3) (#2-#3,1.5+#3) arc(90:0:#3);}
\nc{\udarc}[2]{\udarcx{#1}{#2}{.4}}
\nc{\uudarcx}[3]{\draw(#1,3)arc(180:90:#3) (#1+#3,3+#3)--(#2-#3,3+#3) (#2-#3,3+#3) arc(90:0:#3);}
\nc{\uudarc}[2]{\uudarcx{#1}{#2}{.4}}
\nc{\uarcx}[3]{\draw(#1,1.5)arc(180:270:#3) (#1+#3,1.5-#3)--(#2-#3,1.5-#3) (#2-#3,1.5-#3) arc(270:360:#3);}
\nc{\uarc}[2]{\uarcx{#1}{#2}{.4}}
\nc{\darcxhalf}[3]{\draw(#1,0)arc(180:90:#3) (#1+#3,#3)--(#2,#3) ;}
\nc{\darchalf}[2]{\darcxhalf{#1}{#2}{.4}}
\nc{\uarcxhalf}[3]{\draw(#1,1.5)arc(180:270:#3) (#1+#3,1.5-#3)--(#2,1.5-#3) ;}
\nc{\uarchalf}[2]{\uarcxhalf{#1}{#2}{.4}}
\nc{\uarcxhalfr}[3]{\draw (#1+#3,1.5-#3)--(#2-#3,1.5-#3) (#2-#3,1.5-#3) arc(270:360:#3);}
\nc{\uarchalfr}[2]{\uarcxhalfr{#1}{#2}{.4}}

\nc{\bdarcx}[3]{\draw[blue](#1,0)arc(180:90:#3) (#1+#3,#3)--(#2-#3,#3) (#2-#3,#3) arc(90:0:#3);}
\nc{\bdarc}[2]{\darcx{#1}{#2}{.4}}
\nc{\rduarcx}[3]{\draw[red](#1,0)arc(180:270:#3) (#1+#3,0-#3)--(#2-#3,0-#3) (#2-#3,0-#3) arc(270:360:#3);}
\nc{\rduarc}[2]{\uarcx{#1}{#2}{.4}}
\nc{\duarcx}[3]{\draw(#1,0)arc(180:270:#3) (#1+#3,0-#3)--(#2-#3,0-#3) (#2-#3,0-#3) arc(270:360:#3);}
\nc{\duarc}[2]{\uarcx{#1}{#2}{.4}}

\nc{\uv}[1]{\fill (#1,2)circle(.1);}
\nc{\lv}[1]{\fill (#1,0)circle(.1);}
\nc{\stline}[2]{\draw(#1,2)--(#2,0);}
\nc{\tlab}[2]{\draw(#1,2)node[above]{\tiny $#2$};}
\nc{\tudots}[1]{\draw(#1,2)node{$\cdots$};}
\nc{\tldots}[1]{\draw(#1,0)node{$\cdots$};}

\nc{\uvw}[1]{\fill[white] (#1,2)circle(.1);}
\nc{\huv}[1]{\fill (#1,1)circle(.1);}
\nc{\llv}[1]{\fill (#1,-2)circle(.1);}
\nc{\arcup}[2]{
\draw(#1,2)arc(180:270:.4) (#1+.4,1.6)--(#2-.4,1.6) (#2-.4,1.6) arc(270:360:.4);
}
\nc{\harcup}[2]{
\draw(#1,1)arc(180:270:.4) (#1+.4,.6)--(#2-.4,.6) (#2-.4,.6) arc(270:360:.4);
}
\nc{\arcdn}[2]{
\draw(#1,0)arc(180:90:.4) (#1+.4,.4)--(#2-.4,.4) (#2-.4,.4) arc(90:0:.4);
}
\nc{\cve}[2]{
\draw(#1,2) to [out=270,in=90] (#2,0);
}
\nc{\hcve}[2]{
\draw(#1,1) to [out=270,in=90] (#2,0);
}
\nc{\catarc}[3]{
\draw(#1,2)arc(180:270:#3) (#1+#3,2-#3)--(#2-#3,2-#3) (#2-#3,2-#3) arc(270:360:#3);
}

\nc{\arcr}[2]{
\draw[red](#1,0)arc(180:90:.4) (#1+.4,.4)--(#2-.4,.4) (#2-.4,.4) arc(90:0:.4);
}
\nc{\arcb}[2]{
\draw[blue](#1,2-2)arc(180:270:.4) (#1+.4,1.6-2)--(#2-.4,1.6-2) (#2-.4,1.6-2) arc(270:360:.4);
}
\nc{\loopr}[1]{
\draw[blue](#1,-2) edge [out=130,in=50,loop] ();
}
\nc{\loopb}[1]{
\draw[red](#1,-2) edge [out=180+130,in=180+50,loop] ();
}
\nc{\redto}[2]{\draw[red](#1,0)--(#2,0);}
\nc{\bluto}[2]{\draw[blue](#1,0)--(#2,0);}
\nc{\dotto}[2]{\draw[dotted](#1,0)--(#2,0);}
\nc{\lloopr}[2]{\draw[red](#1,0)arc(0:360:#2);}
\nc{\lloopb}[2]{\draw[blue](#1,0)arc(0:360:#2);}
\nc{\rloopr}[2]{\draw[red](#1,0)arc(-180:180:#2);}
\nc{\rloopb}[2]{\draw[blue](#1,0)arc(-180:180:#2);}
\nc{\uloopr}[2]{\draw[red](#1,0)arc(-270:270:#2);}
\nc{\uloopb}[2]{\draw[blue](#1,0)arc(-270:270:#2);}
\nc{\dloopr}[2]{\draw[red](#1,0)arc(-90:270:#2);}
\nc{\dloopb}[2]{\draw[blue](#1,0)arc(-90:270:#2);}
\nc{\llloopr}[2]{\draw[red](#1,0-2)arc(0:360:#2);}
\nc{\llloopb}[2]{\draw[blue](#1,0-2)arc(0:360:#2);}
\nc{\lrloopr}[2]{\draw[red](#1,0-2)arc(-180:180:#2);}
\nc{\lrloopb}[2]{\draw[blue](#1,0-2)arc(-180:180:#2);}
\nc{\ldloopr}[2]{\draw[red](#1,0-2)arc(-270:270:#2);}
\nc{\ldloopb}[2]{\draw[blue](#1,0-2)arc(-270:270:#2);}
\nc{\luloopr}[2]{\draw[red](#1,0-2)arc(-90:270:#2);}
\nc{\luloopb}[2]{\draw[blue](#1,0-2)arc(-90:270:#2);}

\nc{\larcb}[2]{
\draw[blue](#1,0-2)arc(180:90:.4) (#1+.4,.4-2)--(#2-.4,.4-2) (#2-.4,.4-2) arc(90:0:.4);
}
\nc{\larcr}[2]{
\draw[red](#1,2-2-2)arc(180:270:.4) (#1+.4,1.6-2-2)--(#2-.4,1.6-2-2) (#2-.4,1.6-2-2) arc(270:360:.4);
}

\rnc{\H}{\mathscr H}
\rnc{\L}{\mathscr L}
\nc{\R}{\mathscr R}
\nc{\D}{\mathscr D}
\nc{\J}{\mathscr J}

\nc{\ssim}{\mathrel{\raise0.25 ex\hbox{\oalign{$\approx$\crcr\noalign{\kern-0.84 ex}$\approx$}}}}
\nc{\POI}{\mathcal{POI}}
\nc{\wb}{\overline{w}}
\nc{\ub}{\overline{u}}
\nc{\vb}{\overline{v}}
\nc{\fb}{\overline{f}}
\nc{\gb}{\overline{g}}
\nc{\hb}{\overline{h}}
\nc{\pb}{\overline{p}}
\rnc{\sb}{\overline{s}}
\nc{\XR}{\pres{X}{R\,}}
\nc{\YQ}{\pres{Y}{Q}}
\nc{\ZP}{\pres{Z}{P\,}}
\nc{\XRone}{\pres{X_1}{R_1}}
\nc{\XRtwo}{\pres{X_2}{R_2}}
\nc{\XRthree}{\pres{X_1\cup X_2}{R_1\cup R_2\cup R_3}}
\nc{\er}{\eqref}
\nc{\larr}{\mathrel{\hspace{-0.35 ex}>\hspace{-1.1ex}-}\hspace{-0.35 ex}}
\nc{\rarr}{\mathrel{\hspace{-0.35 ex}-\hspace{-0.5ex}-\hspace{-2.3ex}>\hspace{-0.35 ex}}}
\nc{\lrarr}{\mathrel{\hspace{-0.35 ex}>\hspace{-1.1ex}-\hspace{-0.5ex}-\hspace{-2.3ex}>\hspace{-0.35 ex}}}
\nc{\nn}{\tag*{}}
\nc{\epfal}{\tag*{$\Box$}}
\nc{\tagd}[1]{\tag*{(#1)$'$}}
\nc{\ldb}{[\![}
\nc{\rdb}{]\!]}
\nc{\sm}{\setminus}
\nc{\I}{\mathcal I}
\nc{\InSn}{\I_n\setminus\S_n}
\nc{\dom}{\operatorname{dom}} \nc{\codom}{\operatorname{dom}}
\nc{\ojin}{1\leq j<i\leq n}
\nc{\eh}{\widehat{e}}
\nc{\wh}{\widehat{w}}
\nc{\uh}{\widehat{u}}
\nc{\vh}{\widehat{v}}
\nc{\sh}{\widehat{s}}
\nc{\fh}{\widehat{f}}
\nc{\textres}[1]{\text{\emph{#1}}}
\nc{\aand}{\emph{\ and \ }}
\nc{\iif}{\emph{\ if \ }}
\nc{\textlarr}{\ \larr\ }
\nc{\textrarr}{\ \rarr\ }
\nc{\textlrarr}{\ \lrarr\ }

\nc{\comma}{,\ }

\nc{\COMMA}{,\quad}
\nc{\TnSn}{\T_n\setminus\S_n} 
\nc{\TmSm}{\T_m\setminus\S_m} 
\nc{\TXSX}{\T_X\setminus\S_X} 
\rnc{\S}{\mathcal S}

\nc{\T}{\mathcal T} 
\nc{\A}{\mathscr A} 
\nc{\B}{\mathscr B} 
\rnc{\P}{\mathcal P} 
\nc{\K}{\mathcal K}
\nc{\PB}{\mathcal{PB}} 
\nc{\rank}{\operatorname{rank}}

\nc{\mtt}{\!\!\!\mt\!\!\!}

\nc{\sub}{\subseteq}
\nc{\la}{\langle}
\nc{\ra}{\rangle}
\nc{\mt}{\mapsto}
\nc{\im}{\mathrm{im}}
\nc{\id}{\mathrm{id}}
\nc{\bn}{\mathbf{n}}
\nc{\ba}{\mathbf{a}}
\nc{\bl}{\mathbf{l}}
\nc{\bm}{\mathbf{m}}
\nc{\bk}{\mathbf{k}}
\nc{\br}{\mathbf{r}}
\nc{\al}{\alpha}
\nc{\be}{\beta}
\nc{\ga}{\gamma}
\nc{\Ga}{\Gamma}
\nc{\de}{\delta}
\nc{\ka}{\kappa}
\nc{\lam}{\lambda}
\nc{\Lam}{\Lambda}
\nc{\si}{\sigma}
\nc{\Si}{\Sigma}
\nc{\oijn}{1\leq i<j\leq n}
\nc{\oijm}{1\leq i<j\leq m}

\nc{\comm}{\rightleftharpoons}
\nc{\AND}{\qquad\text{and}\qquad}

\nc{\bit}{\vspace{-3 truemm}\begin{itemize}}
\nc{\bitbmc}{\begin{itemize}\begin{multicols}}
\nc{\bmc}{\vspace{-3 truemm}\begin{itemize}\begin{multicols}}
\nc{\emc}{\end{multicols}\end{itemize}\vspace{-3 truemm}}
\nc{\eit}{\end{itemize}\vspace{-3 truemm}}
\nc{\ben}{\vspace{-3 truemm}\begin{enumerate}}
\nc{\een}{\end{enumerate}\vspace{-3 truemm}}
\nc{\eitres}{\end{itemize}}

\nc{\set}[2]{\{ {#1} : {#2} \}} 
\nc{\bigset}[2]{\big\{ {#1}: {#2} \big\}} 
\nc{\Bigset}[2]{\left\{ \,{#1} :{#2}\, \right\}}

\nc{\pres}[2]{\la {#1} \,|\, {#2} \ra}
\nc{\bigpres}[2]{\big\la {#1} \,\big|\, {#2} \big\ra}
\nc{\Bigpres}[2]{\Big\la \,{#1}\, \,\Big|\, \,{#2}\, \Big\ra}
\nc{\Biggpres}[2]{\Bigg\la {#1} \,\Bigg|\, {#2} \Bigg\ra}

\nc{\pf}{\noindent{\bf Proof.}  }
\nc{\epf}{\hfill$\Box$\bigskip}
\nc{\epfres}{\hfill$\Box$}
\nc{\pfnb}{\pf}
\nc{\epfnb}{\bigskip}
\nc{\pfthm}[1]{\bigskip \noindent{\bf Proof of Theorem \ref{#1}}\,\,  } 
\nc{\pfprop}[1]{\bigskip \noindent{\bf Proof of Proposition \ref{#1}}\,\,  } 
\nc{\epfreseq}{\tag*{$\Box$}}

\makeatletter
\newcommand\footnoteref[1]{\protected@xdef\@thefnmark{\ref{#1}}\@footnotemark}
\makeatother

\numberwithin{equation}{section}

\newtheorem{thm}[equation]{Theorem}
\newtheorem{lemma}[equation]{Lemma}
\newtheorem{cor}[equation]{Corollary}
\newtheorem{prop}[equation]{Proposition}

\theoremstyle{definition}

\newtheorem{rem}[equation]{Remark}
\newtheorem{defn}[equation]{Definition}
\newtheorem{eg}[equation]{Example}
\newtheorem{ass}[equation]{Assumption}

\title{Semigroups of rectangular matrices under a sandwich operation}
\author{
Igor Dolinka%
\\
{\footnotesize \emph{Department of Mathematics and Informatics}}\\
{\footnotesize \emph{University of Novi Sad, Trg Dositeja Obradovi\'ca 4, 21101 Novi Sad, Serbia}}\\
{\footnotesize {\tt dockie\,@\,dmi.uns.ac.rs}}\\~\\
James East\\
{\footnotesize \emph{Centre for Research in Mathematics; School of Computing, Engineering and Mathematics}}\\
{\footnotesize \emph{Western Sydney University, Locked Bag 1797, Penrith NSW 2751, Australia}}\\
{\footnotesize {\tt J.East\,@\,WesternSydney.edu.au}}
}


\maketitle

\vspace{-0.5cm}

\begin{abstract}
%
Let $\mathcal M_{mn}=\mathcal M_{mn}(\mathbb F)$ denote the set of all $m\times n$ matrices over a field $\mathbb F$, and fix some $n\times m$ matrix $A\in\mathcal M_{nm}$.  An associative operation $\star$ may be defined on $\mathcal M_{mn}$ by $X\star Y=XAY$ for all $X,Y\in\mathcal M_{mn}$, and the resulting \emph{sandwich semigroup} is denoted $\mathcal M_{mn}^A=\mathcal M_{mn}^A(\mathbb F)$.  These semigroups are closely related to Munn rings, which are fundamental tools in the representation theory of finite semigroups.  In this article, we study $\mathcal M_{mn}^A$ as well as its subsemigroups $\operatorname{Reg}(\mathcal M_{mn}^A)$ and $\mathcal E_{mn}^A$ (consisting of all regular elements and products of idempotents, respectively), as well as the ideals of $\operatorname{Reg}(\mathcal M_{mn}^A)$.  Among other results, we: characterise the regular elements; determine Green's relations and preorders; calculate the minimal number of matrices (or idempotent matrices, if applicable) required to generate each semigroup we consider; and classify the isomorphisms between finite sandwich semigroups $\mathcal M_{mn}^A(\mathbb F_1)$ and $\mathcal M_{kl}^B(\mathbb F_2)$.  Along the way, we develop a general theory of sandwich semigroups in a suitably defined class of \emph{partial semigroups} related to Ehresmann-style ``arrows only'' categories; we hope this framework will be useful in studies of sandwich semigroups in other categories.  We note that all our results have applications to the \emph{variants} $\mathcal M_n^A$ of the full linear monoid $\mathcal M_n$ (in the case $m=n$), and to certain semigroups of linear transformations of restricted range or kernel (in the case that $\operatorname{rank}(A)$ is equal to one of $m,n$).


\end{abstract}

\section{Introduction}\label{sect:intro}

%
%
%

In the classical representation theory of finite semigroups, a key role is played by the so-called \emph{Munn rings}.  These are rings of $m\times n$ matrices (where $m$ and $n$ need not be equal) with the familiar addition operation but with a \emph{sandwich} multiplication defined by $X\star Y=XAY$, where $A$ is a fixed $n\times m$ matrix.  These rings are so named, because of Douglas Munn's 1955 paper \cite{Munn1955}, in which it was shown that: (1) the representation theory of a finite semigroup is determined by the representations of certain \emph{completely $0$-simple semigroups} arising from its ideal structure, and (2) the semigroup algebra of such a finite completely $0$-simple semigroup is isomorphic to an appropriate Munn ring over the group algebra of a naturally associated maximal subgroup; conditions were also given for such a Munn ring to be semisimple.  (Here, the \emph{sandwich matrix} $A$ arises from the celebrated Rees structure theorem \cite{Rees1940} for completely $0$-simple semigroups.)  Since their introduction in \cite{Munn1955}, Munn rings have been studied by numerous authors, and continue to heavily inflence the theory of semigroup representations: for classical studies, see  \cite{Munn1957,Munn1955,Munn1960, Ponizovskii1956, LP1969,Hall1970, McAlister1971, McAlister1971a,McAlister1971b, Clifford1942,Clifford1960,CPbook}; for modern accounts, see for example \cite{IRS2011,GMS2009,Putcha1998,Putcha1996,OP1991,Steinberg2006,Steinberg2008,AMSV2009}, and especially the monographs \cite{Okninski1991,Okninski1998,Renner2005,Putcha1988,SteinbergBook}.

In the same year as Munn's article \cite{Munn1955} was published, William Brown introduced the so-called \emph{generalised matrix algebras}~\cite{Brown1955}, motivated by a connection with classical groups \cite{Brown1956,Brauer1937,Weyl1939}.  These generalised matrix algebras are again rings of $m\times n$ matrices over a field, with multiplication determined by a fixed $n\times m$ \emph{sandwich matrix}.  Whereas the sandwich matrix in a Munn ring is taken to be the structure matrix of a completely $0$-simple semigroup (and so has a certain prescribed form), Brown considered arbitrary sandwich matrices.  As with Munn rings, these generalised matrix algebras have influenced representation theory to this day, and have been studied by numerous authors; see for example \cite{DL2004,Thrall1955,XW2010,LW2012,XW2014,GW2015,KX2001,Gavarini2008,KX1998}.

Shortly after the Munn and Brown articles \cite{Munn1955,Brown1955} appeared, Evgeny Lyapin's early monograph on semigroups~\cite{Lyapin} was published.  In \cite[Chapter VII]{Lyapin}, we find a number of interesting semigroup constructions, including the following.  Let $V$ and $W$ be arbitrary non-empty sets, and let $\theta:W\to V$ be an arbitrary (but fixed) function.  Then the set $\T(V,W)$ of all functions $V\to W$ forms a semigroup, denoted $\T^\theta(V,W)$, under the operation $\star_\theta$ defined by $f\star_\theta g=f\circ\theta\circ g$.
%
If it is assumed that $V$ and $W$ are vector spaces (over the same field) and~$\theta$ a linear transformation, then the subset $\mathcal L(V,W)\sub\T(V,W)$ of all linear transformations $V\to W$ is a subsemigroup of $\T^\theta(V,W)$.  This subsemigroup, denoted $\mathcal L^\theta(V,W)$ and referred to as a \emph{linear sandwich semigroup}, is clearly isomorphic to the underlying multiplicative semigroup of an associated generalised matrix algebra \cite{Brown1955}.  As noted above, the addition on a generalised matrix algebra is just the usual operation, so these linear sandwich semigroups capture and isolate (in a sense) the more complex of the operations on the algebras.  

The sandwich semigroups $\T^\theta(V,W)$ were first investigated in a series of articles by Magill and Subbiah \cite{Magill1967,MS1975,MS1978}, and more recent studies may be found in \cite{Sullivan_preprint,MGS2013,CC2008,WK2002}; most of these address structural concerns such as (von Neumann) regularity, Green's relations, ideals, classification up to isomorphism, and so on.
%
%
The linear sandwich semigroups $\mathcal L^\theta(V,W)$ have received less attention, though they have also been studied by a number of authors \cite{MGS2014,Kemprasit2002,Chinram2009,JCK2010}, with studies again focusing on basic structural properties.  This is regrettable, because these semigroups display a great deal of algebraic and combinatorial charm, as we hope to show in the current article.  
It is therefore our purpose to carry out a systematic investigation of the linear sandwich semigroups, bringing their study up to date, and focusing on modern themes, especially combinatorial invariant theory.  As does Brown \cite{Brown1955}, we focus on the case that $V$ and $W$ are finite dimensional; in fact, we study the equivalent sandwich semigroups $\MmnA=\MmnA(\F)$ consisting of all $m\times n$ matrices over the field $\F$ under the operation~$\star_A$ defined by $X\star_AY=XAY$, where $A$ is a fixed $n\times m$ matrix.

We speculate that the difficulty (until now) of systematically investigating the linear sandwich semigroups may be due to the lack of a consistent theoretical framework for studying sandwich semigroups in more generality.  In the case that $V=W$, the sets $\T(V,W)$ and $\mathcal L(V,W)$ are themselves semigroups (under composition); these are the \emph{full transformation semigroup} $\T_V$ \cite{GMbook,MGS2013,Gomes1987,GR2012,EMP2015, Howie1990,Howie1966,Howie1978} and the \emph{general linear monoid} $\mathcal L_V$ \cite{DG2014,Djokovic1968,Erdos1967,Dawlings81/82,AM2005, Dawlings1982,Gray2008,Waterhouse,Putcha2006, Laffey1983,Okninski1998}, respectively.  In turn, the semigroups $\mathcal T^\theta(V,V)$ and $\mathcal L^\theta(V,V)$ are special cases of the \emph{semigroup variant} construction.  The \emph{variant} of a semigroup $S$ with respect to an element $a\in S$ is the semigroup $S^a=(S,\star_a)$, with operation defined by $x\star_ay=xay$.  Variants were first explicitly studied by Hickey in the 1980s \cite{Hickey1983,Hickey1986}, though (as noted above) the idea goes back to Lyapin's monograph \cite{Lyapin}; a more recent study may be found in \cite{KL2001}.  The current authors developed the general theory of variants further in \cite{DE2}, and then used this as a starting point to explore the variants of the finite full transformation semigroups, obtaining a great deal of algebraic and combinatorial information about these semigroups.  Unfortunately, the theory of semigroup variants does not help with studying the more general sandwich semigroups $\mathcal T^\theta(V,W)$ and $\mathcal L^\theta(V,W)$, since the underlying sets $\T(V,W)$ and $\mathcal L(V,W)$ are not even semigroups if $V\not=W$.  One of the main goals of the current article, therefore, is to develop an appropriate general framework for working with arbitrary sandwich semigroups.  Namely, if $V$ and $W$ are objects in a (locally) small category $\C$, and if $\theta\in\Hom(W,V)$ is some fixed morphism, then the set $\Hom(V,W)$ becomes a semigroup under the sandwich operation defined by $f\star_\theta g=f\circ\theta\circ g$, for $f,g\in\Hom(V,W)$.  (In the case that $V=W$ and $\theta$ is the identity morphism, this construction reduces to the usual endomorphism monoid $\End(V)$.)  The semigroups $\mathcal T^\theta(V,W)$ and $\mathcal L^\theta(V,W)$ arise when $\C$ is the category of sets (and mappings) or vector spaces (and linear transformations), respectively.  In order to develop a general theory of sandwich semigroups in such categories, we first explain how many important semigroup theoretical techniques 
extend to the more general categorical setting; we note that there is only a little overlap with the theory of Green's relations in categories developed in \cite{LS2012}, which  focuses on issues more relevant to representation theory.  In order to avoid any confusion arising from terminology conflicts between semigroup and category theory, rather than speak of (locally small) categories, we focus on the equivalently defined class of \emph{partial semigroups}, which are related to Ehresmann-style ``arrows only'' categories \cite{Ehresmann1965}.  We hope that the general theory we develop will prove to be a useful starting point for future studies of sandwich semigroups in other categories.

The article is organised as follows.  
In Section \ref{sect:partial}, we develop a general theory of sandwich semigroups in partial semigroups (i.e., locally finite categories), extending certain important semigroup theoretic notions (such as Green's relations, regularity and stability, the definitions of which are given in Section~\ref{sect:partial}) to the more general context.  
In Section \ref{sect:preliminaries}, we gather results on the partial semigroup~$\M=\M(\F)$ of all (finite dimensional) matrices over the field $\F$, mainly focusing on regularity, stability and Green's relations, and we state some well-known results on (idempotent) generation and ideals of the general linear monoids~$\M_n$.  
We begin our investigation of the linear sandwich semigroups $\MmnA$ in Section \ref{sect:MmnJ}, the main results of this section being: a characterisation of the regular elements (Proposition~\ref{prop:P1P2}); a description of Green's relations (Theorem~\ref{green_thm}) and the ordering on $\D$-classes (Propositions \ref{prop:DorderMmnJ}, \ref{prop:DorderP} and \ref{prop_maximalD}); a classification of the isomorphism classes of sandwich semigroups over $\Mmn$ (Corollary \ref{cor:MmnAcongMmnB}); and the calculation of $\rank(\MmnA)$ (Theorems~\ref{thm:rankMmnJ} and~\ref{thm:rankMmnJ_r=m}).  
(Recall that the \emph{rank} of a semigroup $S$, denoted $\rank(S)$, is the minimum size of a generating set for $S$.)
Section \ref{sect:non-sandwich} explores the relationship between a sandwich semigroup $\MmnA$ and various (non-sandwich) matrix semigroups, the main structural results being Theorem \ref{thm:diamondsMmnJ} and Propositions~\ref{mono_prop} and~\ref{prop:MAN}.  
%
We then focus on the regular subsemigroup $P=\Reg(\MmnA)$ in Section \ref{sect:RegMmnJ}, where we: 
calculate the size of $P$ and various Green's classes (Proposition \ref{prop:DXJ_combinatorics} and Theorem \ref{inflation_thm}); classify the isomorphism classes of finite linear sandwich semigroups (Theorem \ref{thm:classification}); and calculate $\rank(P)$ (Theorem \ref{thm:rankP}).  
In Section \ref{sect:EMmnJ}, we investigate the idempotent generated subsemigroup $\E_{mn}^A$ of $\MmnA$, where we: enumerate the idempotents of $\MmnA$ (Proposition \ref{prop:enumeration_E}); show that $\E_{mn}^A$ consists of $P\sm D$ and the idempotents from $D$, where $D$ is the maximal $\D$-class (Theorem \ref{thm:EmnJ}); and calculate $\rank(\E_{mn}^A)$ and $\idrank(\E_{mn}^A)$, showing in particular that these are equal (Theorem \ref{thm:rankEmnJ}).  
(The \emph{idempotent rank} of an idempotent generated semigroup $S$, denoted $\idrank(S)$, is defined similarly to the rank, but with respect to idempotent generating sets for $S$.)
Finally, in Section \ref{sect:ideals}, we classify the proper ideals of $P$, showing that these are idempotent generated, and calculating their ranks and idempotent ranks, which are again equal (Theorem~\ref{thm:ideals}).
We note that all our results have applications to the \emph{variants} $\M_n^A$ of the full linear monoid $\M_n$ (in the case $m=n$), and to certain semigroups of linear transformations of restricted range or kernel (in the case that $\rank(A)$ is equal to one of $m,n$; see Remarks~\ref{rem:r=m} and \ref{rem:r=m2}).%

\section{Sandwich semigroups from partial semigroups}\label{sect:partial}

Recall that our main interest is in the linear sandwich semigroups $\MmnA=\MmnA(\F)$.  The underlying set of $\MmnA$ is $\Mmn$, the set of all $m\times n$ matrices over the field $\F$, which is not itself a semigroup (unless $m=n$).  However, $\Mmn$ is contained in $\M$, the set of \emph{all} (finite dimensional) matrices over $\F$.  While $\M$ is still not a semigroup, it does have the structure of a 
(small) category.  As we will see, in order to understand the linear sandwich semigroups $\MmnA$, we need to move beyond just $m\times n$ (and $n\times m$) matrices, and gain a fuller understanding of the whole category $\M$.  Some (but not all) of what we need
to know about $\M$ is true in a larger class of categories, and more general structures we call \emph{partial semigroups}, so we devote this section to the development of the general theory of these structures.  We begin with the definitions.



\ms
\begin{defn}\label{defn:S}
A \emph{partial semigroup} is a $5$-tuple $(S,\cdot,I,\lam,\rho)$ consisting of a set $S$, a partial binary operation $(x,y)\mt x\cdot y$ (defined on some subset of $S\times S$), a set $I$, and functions $\lam,\rho:S\to I$, such that, for all $x,y,z\in S$,
\bit
\item[(i)] $x\cdot y$ is defined if and only if $\rho(x)=\lam(y)$,
\item[(ii)] if $x\cdot y$ is defined, then $\lam(x\cdot y)=\lam(x)$ and $\rho(x\cdot y)=\rho(y)$,
\item[(iii)] if $x\cdot y$ and $y\cdot z$ are defined, then $(x\cdot y)\cdot z=x\cdot (y\cdot z)$.
\eit
We say that a partial semigroup $(S,\cdot,I,\lam,\rho)$ is \emph{monoidal} if in addition to (i--iii), 
\bit
\item[(iv)] there exists a function $I\to S:i\mt e_i$ such that, for all $x\in S$, $x\cdot e_{\rho(x)}=x=e_{\lam(x)}\cdot x$.
\eit
We say that a partial semigroup $(S,\cdot,I,\lam,\rho)$ is \emph{regular} if in addition to (i--iii),
\bit
\item[(v)] for all $x\in S$, there exists $y\in S$ such that $x=x\cdot y\cdot x$ and $y=y\cdot x\cdot y$.
\eit
\end{defn}


\ms
{\begin{rem}
We note that conditions (i--iv) amount to one of several equivalent ways to define (small) categories in an ``arrows only'' fashion.  See for example Ehresmann's monograph \cite{Ehresmann1965}, and also \cite{Hollings2012} for a historical discussion of the connections between category theory and (inverse) semigroup theory.  
\end{rem}

For a partial semigroup $(S,\cdot,I,\lam,\rho)$, and for $i,j\in I$, we write
\[
S_{ij} = \set{x\in S}{\lam(x)=i,\ \rho(x)=j} \AND S_i = S_{ii}.
\]
So $S=\bigcup_{i,j\in I}S_{ij}$.  Note that if $x\in S$, then $x\cdot x$ is defined if and only if $\lam(x)=\rho(x)$.  It follows that $S_i$ is a semigroup with respect to the induced binary operation (the restriction of $\cdot$ to $S_i\times S_i$) for each $i\in I$, but that $S_{ij}$ is not if $i\not=j$.  We will often slightly abuse notation and refer to ``the partial semigroup $S$'' if the rest of the data $(S,\cdot,I,\lam,\rho)$ is clear from context.  We also note that in what follows, we could allow~$S$ and~$I$ to be classes (rather than insist on them being sets); but we would still require $S_{ij}$ to be a set for each $i,j\in I$.  

Note that, as is the case with semigroups, condition (v) is equivalent to the (ostensibly) weaker condition:
\bit
\item[(v)$'$] for all $x\in S$, there exists $z\in S$ such that $x=x\cdot z\cdot x$.
\eit
Indeed, with $z$ as in (v)$'$, one easily checks that $y=z\cdot x\cdot z$ satisfies the condition of (v).

%
%
If $S$ is monoidal, then $S_i$ is a monoid with identity $e_i\in S_i$ for each $i$.  
If $S$ is not monoidal, then $S$ may be embedded in a monoidal partial semigroup $\Sone$ as follows: for each $i\in I$ we adjoin an element $e_i$ to $S_i$ and declare that $x\cdot e_i=x$ and $e_i\cdot y=y$ for all $x,y\in S$ with $\rho(x)=i$ and $\lam(y)=i$, if such an element $e_i\in S_i$ does not already exist.  In particular, if $S$ is monoidal, then $S=\Sone$.

Obviously any semigroup is a partial semigroup (with $|I|=1$); in particular, all results we prove in this section concerning partial semigroups hold for semigroups.
A great number of non-semigroup examples exist, but we will limit ourselves to describing just a few.

\ms
\begin{eg}
As a trivial example, let $\set{S_i}{i\in I}$ be any set of pairwise disjoint semigroups.  Then $S=\bigcup_{i\in I}S_i$ is a partial semigroup where we define $\lam,\rho:S\to I$ by $\lam(x)=\rho(x)=i$ for each $i\in I$ and  $x\in S_i$, and $x\cdot y$ is defined if and only if $x,y\in S_i$ for some $i$, in which case $x\cdot y$ is just the product of $x,y$ in $S_i$.  Note that this $S$ is regular (resp., monoidal) if and only if each $S_i$ is regular (resp., a monoid).
\end{eg}

\nc{\X}{\mathscr X}

\ms
\begin{eg}
Let $\X$ be some set, and $\P(\X)=\set{A}{A\sub\X}$ the power set of $\X$.  The set $T_\X=\set{(B,f,A)}{A,B\sub\X,\ \text{$f$ is a function $A\to B$}}$ is a regular monoidal partial semigroup.  We define $I=\P(\X)$, and $\lam(B,f,A)=B$ and $\rho(B,f,A)=A$, with $(D,g,C)\cdot(B,f,A)$ defined if and only if $B=C$, in which case $(D,g,C)\cdot(B,f,A)=(D,g\circ f,A)$.  
\end{eg}

The previous example may be extended in a number of ways, by replacing functions $f:A\to B$ by other objects such as binary relations \cite{Thornton1982, Chase1979}, partial functions \cite{Sullivan1975,Chinram2008}, partial bijections \cite{Chinram2008b}, block bijections \cite{FL1998}, partial braids \cite{EL2004}, partitions \cite{Martin1994}, Brauer diagrams \cite{Brauer1937}, etc., or by assuming the functions $f:A\to B$ preserve some kind of algebraic or geometric structure on the sets $A,B$.  The main example we will concentrate on in this article is as follows.

\ms
\begin{eg}\label{eg:M}
Let $\F$ be a field, and write $\M=\M(\F)$ for the set of all (finite dimensional, non-empty) matrices over~$\F$.  Then $\M$ has the structure of a regular monoidal partial semigroup.  We take $I=\mathbb N=\{1,2,3,\ldots\}$ to be the set of all natural numbers and, for $X\in\M$, we define $\lam(X)$ (resp., $\rho(X)$) to be the number of rows (resp., columns) of $X$.  For $m,n\in\mathbb N$, $\Mmn=\Mmn(\F)$ denotes the set of all $m\times n$ matrices over $\F$, and forms a semigroup if and only if $m=n$.  (Of course, $\M$ is isomorphic to a certain partial semigroup of linear transformations; we will have more to say about this later.)
\end{eg}

For the remainder of this section, we fix a partial semigroup $(S,\cdot,I,\lam,\rho)$, and we write $xy$ for the product $x\cdot y$ (whenever it is defined).  
Note that we may define a second partial binary operation $\bullet$ on $S$ by
\[
x\bullet y=y\cdot x \qquad\text{for each $x,y\in S$ with $\rho(y)=\lam(x)$.}
\]
We see then that $(S,\bullet,I,\rho,\lam)$ is a partial semigroup (note the swapping of $\lam$ and $\rho$), and we call this the \emph{dual partial semigroup} to $(S,\cdot,I,\lam,\rho)$.  As is frequently the case in semigroup theory, this duality will allow us to shorten several proofs.

Green's relations and preorders are crucial tools in semigroup theory (for general background on semigroups, see \cite{Hig,Howie}), and we will need to extend these to the partial semigroup setting.  
%
If $x,y\in S$, then we say
\bit
\item $x\leq_\R y$ if $x=ya$ for some $a\in\Sone$,
\item $x\leq_\L y$ if $x=ay$ for some $a\in \Sone$,
\item $x\leq_\J y$ if $x=ayb$ for some $a,b\in\Sone$.
\eit
Note that if $x\leq_\R y$ (resp., $x\leq_\L y$), then $\lam(x)=\lam(y)$ (resp., $\rho(x)=\rho(y)$).  
%
Note also that if $x\leq_\R y$, then $ux\leq_\R uy$ for any $u\in S$ with $\rho(u)=\lam(x)$; a dual statement holds for the $\leq_\L$ relation.
Finally, note that the use of $\Sone$ is merely for convenience since, for example, $x\leq_\R y$ means that $x=y$ or $x=ya$ for some $a\in S$.  
All three of the above relations are preorders (i.e., they are reflexive and transitive).
If $\gK$ is one of $\R$, $\L$, $\J$, we write $\gK={\leq_\gK}\cap{\geq_\gK}$ for the equivalence relation on $S$ induced by $\gK$.  So, for example, $x\R y$ if and only if $x=ya$ and $y=xb$ for some $a,b\in \Sone$.  
%
%
We also define equivalence relations
\[
\H=\R\cap\L \AND \D=\R\vee\L.
\]
(The join $\ve\vee\eta$ of two equivalences $\ve$ and $\eta$ is the transitive closure of $\ve\cup\eta$, and is itself an equivalence.)  
It is easy to see that $\D\sub\J$.  
The duality mentioned above means that $x\leq_\R y$ in $(S,\cdot,I,\lam,\rho)$ if and only if $x\leq_\L y$ in $(S,\bullet,I,\rho,\lam)$, and so on.

Analogously to the definition for semigroups \cite[Definition A.2.1]{RSbook}, we say that the partial semigroup $S$ is \emph{stable} if for all $x,y\in S$,
\[
x\J xy \ \iff\ x\R xy \AND x\J yx\ \iff\ x\L yx.
\]
The following simple but crucial observation is proved in analogous fashion to the corresponding results for semigroups; see for example \cite[Proposition 2.1.3]{Howie} and \cite[Corollary A.2.5]{RSbook}.  

\ms
\begin{lemma}\label{lem:Rol=LoR_S}
We have $\D=\R\circ\L=\L\circ\R$. If $S$ is stable, then $\D=\J$.  \epfres
\end{lemma}

If $x\in S_{ij}$ and $\gK$ is one of $\R$, $\L$, $\J$, $\D$, $\H$, we write 
\[
\xK = \set{y\in S}{x\gK y} \AND K_x=\xK\cap S_{ij} = \set{y\in S_{ij}}{x\gK y}.
\]
We call $\xK$ (resp., $K_x$) the \emph{$\gK$-class of $x$ in $S$} (resp., \emph{in $S_{ij}$}).
The next result is reminiscent of Green's Lemma, and may be proved in virtually identical fashion to \cite[Lemma 2.2.1]{Howie}.

\ms
\begin{lemma}\label{lem:GreensLemma_S}
Let $x,y\in S$.
\bit
\itemit{i} Suppose $x\R y$, and that $x=ya$ and $y=xb$ where $a,b\in\Sone$.  Then the maps $\xL\to \yL:w\mt wb$ and $ \yL\to  \xL:w\mt wa$ are mutually inverse bijections.  These maps restrict to mutually inverse bijections $\xH\to \yH$ and $\yH\to \xH$.
\itemit{ii} Suppose $x\L y$, and that $x=ay$ and $y=bx$ where $a,b\in\Sone$.  Then the maps $\xR\to \yR:w\mt bw$ and $\yR\to \xR:w\mt aw$ are mutually inverse bijections.  These maps restrict to mutually inverse bijections $\xH\to \yH$ and $\yH\to \xH$.  
\itemit{iii} If $x\D y$, then $\big|\xR\big|=\big|\yR\big|$, $\big|\xL\big|=\big|\yL\big|$ and $\big|\xH\big|=\big|\yH\big|$. \epfres
\eitres
\end{lemma}



Note that if $x,y\in S$ are such that $x\H y$, then $\lam(x)=\lam(y)$ and $\rho(x)=\rho(y)$.  It follows that $\xH=H_x$ for all $x\in S$.  

\ms
\begin{lemma}\label{lem:GreensLemma_Sij}
Let $x,y\in S_{ij}$.
\bit
\itemit{i} Suppose $x\R y$, and that $x=ya$ and $y=xb$ where $a,b\in\Sone$.  Then the maps $L_x\to L_y:w\mt wb$ and $L_y\to L_x:w\mt wa$ are mutually inverse bijections.  These maps restrict to mutually inverse bijections $H_x\to H_y$ and $H_y\to H_x$.
\itemit{ii} Suppose $x\L y$, and that $x=ay$ and $y=bx$ where $a,b\in\Sone$.  Then the maps $R_x\to R_y:w\mt bw$~and $R_y\to R_x:w\mt aw$ are mutually inverse bijections.  These maps restrict to mutually inverse bijections $H_x\to H_y$ and $H_y\to H_x$.
\itemit{iii} If $x\D y$, then $|R_x|=|R_y|$, $|L_x|=|L_y|$ and $|H_x|=|H_y|$.
\eitres
\end{lemma}

\pf Suppose $x\R y$, and that $x=ya$ and $y=xb$ where $a,b\in\Sone$.  We first show that the map $f:L_x\to S:w\mt wb$ does indeed map $L_x$ into $L_y$.  With this in mind, let $w\in L_x$.  We already know that $wb\in \yL$, by Lemma \ref{lem:GreensLemma_S}(i).  Also, $w=ux$ for some $u\in\Sone$, since $w\L x$.  Now, $\lam(wb)=\lam(w)=i$, and also $\rho(wb)=\rho(uxb)=\rho(uy)=\rho(y)=j$, showing that $wb\in \yL\cap S_{ij}=L_y$, as required.  By symmetry, it follows that $g:L_y\to S:w\mt wa$  maps $L_y$ into $L_x$.  By Lemma \ref{lem:GreensLemma_S}(i), we see that $f\circ g$ and $g\circ f$ are the identity maps on their respective domains.  This completes the proof of (i). 

Next, note that (ii) follows from (i) by duality.  Now suppose $x\D y$.  So $x\R z\L y$ for some $z\in S$.  Since $x\R z$, it follows that $\lam(z)=\lam(x)=i$; similarly, $\rho(z)=j$, so in fact,  $z\in S_{ij}$.  In particular, $R_x=R_z$ and $L_y=L_z$.  The statement about cardinalities then follows from parts (i) and (ii). \epf

As is the case for semigroups \cite{Hig,Howie}, Lemma \ref{lem:Rol=LoR_S} means that the elements of a $\D$-class $D$ of $S$ or $S_{ij}$ may be grouped together in a rectangular array of cells, which (for continuity with semigroup theory) we call an \emph{eggbox}.  We place all elements from $D$ in a box in such a way that $\R$-related (resp., $\L$-related) elements are in the same row (resp., column), and $\H$-related elements in the same cell.  An example is given in Figure~\ref{fig:eggbox} below for a $\D$-class of the linear partial semigroup $\M(\mathbb Z_3)$.

We now come to the definition of the main objects of our study, the \emph{sandwich semigroups}.

\ms
\begin{defn}
Let $(S,\cdot,I,\lam,\rho)$ be a partial semigroup.  Fix some $a\in S_{ji}$, where $i,j\in I$.  Define a binary operation $\star_a$ on $S_{ij}$ by $x\star_a y=xay$ for each $x,y\in S_{ij}$.  It is easily checked that $\star_a$ is associative.  We denote by $\Sija=(S_{ij},\star_a)$ the semigroup obtained in this way, and call $\Sija$ the \emph{sandwich semigroup} of $S_{ij}$ with respect to $a$.  (Note that when $i=j$, $\Sija=S_i^a$ is the well-known \emph{variant} \cite{Hickey1983,Hickey1986,KL2001} of $S_i$ with respect to $a\in S_i$.)  
\end{defn}

Recall that an element $x$ of a semigroup $T$ is \emph{regular} if $x=xyx$ and $y=yxy$ for some $y\in T$ (or, equivalently, if $x=xzx$ for some $z\in T$).  The set of all regular elements of $T$ is denoted by $\Reg(T)$, and we say $T$ is \emph{regular} if $T=\Reg(T)$.  (In general, $\Reg(T)$ need not even be a subsemigroup of $T$.)  Of crucial importance is that if any element of a $\D$-class $D$ of a semigroup $T$ is regular, then \emph{every} element of $D$ is regular, in which case every element of $D$ is $\L$-related to at least one idempotent (and also $\R$-related to a possibly different idempotent); the $\H$-class $H_e$ of an idempotent $e\in E(T)=\set{x\in T}{x=x^2}$ is a group, and $H_e\cong H_f$ for any two $\D$-related idempotents~$e,f\in E(T)$.  When drawing eggbox diagrams, group $\gH$-classes are usually shaded grey (see for example Figure \ref{fig:M0...M3}).  See \cite{Hig,Howie} for more details.  

If $S$ is a regular partial semigroup, then the sandwich semigroups $\Sija$ need not be regular themselves (although all of the semigroups $S_i$ are), but the set $\RegSija$ of all regular elements of $\Sija$ forms a subsemigroup, as we now show.

\begin{prop}\label{prop:regularSija}
Let $(S,\cdot,I,\lam,\rho)$ be a regular partial semigroup.  Then $\RegSija$ is a subsemigroup of $\Sija$ for all $i,j\in I$ and $a\in S_{ji}$.
\end{prop}

\pf Let $x,y\in\RegSija$, so $x=xauax$ and $y=yavay$ for some $u,v\in S_{ij}$.  Since $S$ is regular, there exists $w\in S$ such that $(auaxayava)w(auaxayava)=(auaxayava)$.  Then
\begin{align*}
(xay)a(vawau)a(xay) &= (xauaxay)a(vawau)a(xayavay) \\
&= x(auaxayava)w(auaxayava)y = x(auaxayava)y = xay,
\end{align*}
showing that $(x\star_ay)\star(v\star_aw\star_au)\star_a(x\star_ay)=x\star_ay$, and $x\star_ay\in\RegSija$. \epf

In order to say more about the regular elements and Green's relations of the sandwich semigroup $\Sija$, we define the sets
\[
P_1^a = \set{x\in S_{ij}}{xa\R x} \COMMA
P_2^a = \set{x\in S_{ij}}{ax\L x} \COMMA
P_3^a = \set{x\in S_{ij}}{axa\J x} \COMMA
P^a=P_1^a\cap P_2^a.
\]
The next result explains the relationships that hold between these sets; the various inclusions are pictured in Figure \ref{fig:P}.

\ms
\begin{figure}[ht]
\begin{center}
\scalebox{.8}{
\begin{tikzpicture}
  \tikzset{venn circle/.style={draw,circle,minimum width=6cm,fill=#1,opacity=0.4}}
  \tikzset{small venn circle/.style={draw,circle,minimum width=1cm,fill=#1,opacity=0.4}}
  \node [venn circle = black,opacity=0.2] (C) at (4,0) {$$};
  \node [venn circle = black,opacity=0.2] (A) at (0,0) {$$};
  \node [venn circle = black,opacity=0.2] (B) at (2,-0.76393) {$$};
  \node [small venn circle = black,opacity=0.2] (B) at (2,0) {$$};
  \node at (2,0) {$R$};   
  \node at (-0.76393*2,0.76393*2) {$P_1^a$};   
  \node at (0.76393*2+4,0.76393*2) {$P_2^a$};   
  \node at (2,-0.76393*4) {$P_3^a$};   
  \node at (2,-1.2) {$P^a$};   
\end{tikzpicture} 
\qquad\qquad
\begin{tikzpicture}
  \tikzset{venn circle/.style={draw,circle,minimum width=6cm,fill=#1,opacity=0.4}}
  \tikzset{small venn circle/.style={draw,circle,minimum width=1cm,fill=#1,opacity=0.4}}
  \node [venn circle = black,opacity=0.2] (C) at (4,0) {$$};
  \node [venn circle = black,opacity=0.2] (A) at (0,0) {$$};
  \node [venn circle = black,opacity=0.0] (B) at (2,-0.76393) {$$};
  \node [small venn circle = black,opacity=0.2] (B) at (2,0) {$$};
  \node at (2,0) {$R$};   
  \node at (-0.76393*2,0.76393*2) {$P_1^a$};     \node at (0.76393*2+4,0.76393*2) {$P_2^a$};   
  \node at (2,-1.2) {$P^a=P_3^a$};   
\end{tikzpicture} 
}
    \caption{Venn diagrams illustrating the various relationships between the sets $P_1^a$, $P_2^a$, $P_3^a$, $P^a=P_1^a\cap P_2^a$ and $\Reg(\Sija)$ in the general case (left) and the stable case (right); for clarity, we have written $R=\Reg(\Sija)$.}
    \label{fig:P}
   \end{center}
 \end{figure}

\ms
\begin{prop}\label{prop:Reg(Sija)}
Let $(S,\cdot,I,\lam,\rho)$ be a partial semigroup, and fix $i,j\in I$ and $a\in S_{ji}$.  Then
\begin{itemize}\begin{multicols}2
\itemit{i} $\RegSija \sub P^a \sub P_3^a$,
\itemit{ii} $P^a=P_3^a$ if $S$ is stable.
\emc
\end{prop}

\pf If $x\in\RegSija$, then $x=xayax$ for some $y\in S_{ij}$, giving $x\R xa$ and $x\L ax$, so that $x\in P_1^a\cap P_2^a=P^a$.  Next, suppose $x\in P^a=P_1^a\cap P_2^a$, so $x=xav=uax$ for some $u,v\in\Sone$.  It follows that $x=uaxav$, so $x\J axa$ and $x\in P_3^a$.  This completes the proof of (i).


Now suppose $S$ is stable, and let $x\in P_3^a$.  So $x=uaxav$ for some $u,v\in\Sone$.  It then follows that $x\J xa$ and $x\J ax$.  By stability, it follows that $x\R xa$ and $x\L ax$, so that $x\in P_1^a\cap P_2^a=P^a$, completing the proof of (ii). \epf

\begin{rem}
The assumption of regularity (resp., stability) could be greatly weakened in Proposition~\ref{prop:regularSija} (resp., Proposition \ref{prop:Reg(Sija)}(ii)).  
However, because the linear partial semigroup $\M$ is regular and stable
(see Lemmas \ref{lem:green<M} and \ref{lem:regularity}), we will not pursue this thought any further.  
\end{rem}




We now show how the sets $P_1^a$, $P_2^a$, $P_3^a$ and $P^a=P_1^a\cap P_2^a$ may be used to relate Green's relations on the sandwich semigroups $\Sija$ to the corresponding relations on $S$.
To avoid confusion, if $\gK$ is one of $\R$, $\L$, $\J$, $\D$, $\H$, we write $\gKa$ for the Green's $\gK$-relation on $\Sija$.  So, for example, if $x,y\in S_{ij}$, then 
\bit
\item $x\gRa y$ if and only if [$x=y$] or [$x=y\star_a u=yau$ and $y=x\star_a v=xav$ for some $u,v\in S_{ij}$].  
\eit
It is then clear that $\gRa\sub\R$, and the analogous statement is true for all of the other Green's relations.
If $x\in S_{ij}$, we write $K_x^a=\set{y\in S_{ij}}{x\gKa y}$ for the $\gKa$-class of $x$ in $\Sija$.  Since $\gKa\sub\gK$, it follows that $K_x^a\sub K_x$ for all $x\in S_{ij}$.  
%
%

\ms
\begin{thm}\label{thm:green_Sij}
Let $(S,\cdot,I,\lam,\rho)$ be a partial semigroup, and let $a\in S_{ji}$ where $i,j\in I$.  If $x\in S_{ij}$, then   
\ms
\begin{itemize}\begin{multicols}{2}
\itemit{i} $R_x^a = \begin{cases}
R_x\cap P_1^a &\text{if $x\in P_1^a$}\\
\{x\} &\text{if $x\in S_{ij}\sm P_1^a$,}
\end{cases}$
\itemit{ii} $L_x^a = \begin{cases}
L_x\cap P_2^a &\hspace{0.7mm}\text{if $x\in P_2^a$}\\
\{x\} &\hspace{0.7mm}\text{if $x\in S_{ij}\sm P_2^a$,}
\end{cases}
$
\itemit{iii} $H_x^a = \begin{cases}
H_x &\hspace{7.4mm}\text{if $x\in P^a$}\\
\{x\} &\hspace{7.4mm}\text{if $x\in S_{ij}\sm P^a$,}
\end{cases}$
\itemit{iv} $D_x^a = \begin{cases}
D_x\cap P^a &\text{if $x\in P^a$}\\
L_x^a &\text{if $x\in P_2^a\sm P_1^a$}\\
R_x^a &\text{if $x\in P_1^a\sm P_2^a$}\\
\{x\} &\text{if $x\in S_{ij}\sm (P_1^a\cup P_2^a)$,}
\end{cases}$
\itemit{v} $J_x^a = \begin{cases}
J_x\cap P_3^a &\hspace{2.2mm}\text{if $x\in P_3^a$}\\
D_x^a &\hspace{2.2mm}\text{if $x\in S_{ij}\sm P_3^a$.}
\end{cases}$
\end{multicols}\end{itemize}
Further, if $x\in S_{ij}\sm P^a$, then $H_x^a=\{x\}$ is a non-group $\gHa$-class of $\Sija$.  
\end{thm}

\pf The proof of \cite[Proposition 3.2]{DE2} may easily be adapted to prove (i--iv) and the final statement about $\gHa$-classes.  We now prove (v).  Let $x\in S_{ij}$.

Suppose $y\in J_x^a\sm\{x\}$.  So one of (a--c) and one of (d--f) holds:
\bitbmc2
\item[(a)] $y=sax$ for some $s\in S_{ij}$,
\item[(b)] $y=xat$ for some $t\in S_{ij}$,
\item[(c)] $y=saxat$ for some $s,t\in S_{ij}$,
\item[(d)] $x=uay$ for some $u\in S_{ij}$,
\item[(e)] $x=yav$ for some $v\in S_{ij}$,
\item[(f)] $x=uayav$ for some $u,v\in S_{ij}$.
\emc
Suppose first that (a) and (d) hold.  Then $x\gLa y$.  Since $x\not=y$, we deduce that $x\in P_2^a$ by (ii).  Since $L_x^a=L_y^a$, we also have $y\in P_2^a$.  Similarly, if (b) and (e) hold, then $x\gRa y$ and $x,y\in P_1^a$.  One may check that any other combination of (a--c) and (d--f) implies $x,y\in P_3^a$.  For example, if (a) and (e) hold, then
\[
y=sax=s(aya)v \AND x=yav=s(axa)v.
\]
In particular, we have shown that $|J_x^a|\geq2$ implies $x\in P_1^a\cup P_2^a\cup P_3^a$.  By the contrapositive of this last statement, if $z\in S_{ij}\sm(P_1^a\cup P_2^a\cup P_3^a)$, then $J_z^a=\{z\}=D_z^a$, with the last equality following from (iv).  

Next, suppose $x\in P_1^a\sm P_3^a$.  In particular, $x\not\in P_2^a$ since $P_1^a\cap P_2^a\sub P_3^a$ by Proposition \ref{prop:Reg(Sija)}(i).  Since $\gDa\sub\gJa$, we have $D_x^a\sub J_x^a$.  Conversely, suppose $y\in J_x^a$.  We must show that $y\in D_x^a$.  If $y=x$, then we are done, so suppose $y\not=x$.  As above, one of (a--c) and one of (d--f) holds.  If (b) and (e) hold, then $y\in R_x^a=D_x^a$, the second equality holding by (iv).  If any other combination of (a--c) and (d--f) holds then, as explained in the previous paragraph, $x$ (and $y$) would belong to $P_2^a$ or $P_3^a$, a contradiction.  This completes the proof that $J_x^a\sub D_x^a$.  A dual argument shows that $J_x^a=D_x^a$ if $x\in P_2^a\sm P_3^a$.

Finally, suppose $x\in P_3^a$.  Let $z\in J_x\cap P_3^a$.  So we have
\[
x=s'axat' \COMMA z=s''azat'' \COMMA z=u'xv' \COMMA x=u''zv'' \qquad\text{for some $s',s'',t',t'',u',u'',v',v''\in \Sone$.}
\]
We then calculate
$
z = u'xv' = u's'axat'v' = u's'a(s'axat')at'v' = (u's'as') \star_a x \star_a (t'at'v'),
$
and similarly $x = (u''s''as'') \star_a z \star_a (t''at''v'')$, showing that $z\gJa x$, and $J_x\cap P_3^a\sub J_x^a$.  To prove the reverse inclusion, since we have already observed that $J_x^a\sub J_x$, it suffices to show that $J_x^a\sub P_3^a$.  So suppose $y\in J_x^a$.  If $y=x$, then $y\in P_3^a$, so suppose $y\not=x$.  Then one of (a--c) and one of (d--f) above holds.  If (a) and (d) hold, then
\[
y = sax = sa s'axat' = sa s'auayat',
\]
showing that $y\in P_3^a$.  A similar argument covers the case in which (b) and (e) hold.  As we observed above, any other combination of (a--c) and (d--f) implies that $y\in P_3^a$.  This completes the proof. \epf

For a pictorial understanding of Theorem \ref{thm:green_Sij}, 
Figures \ref{fig:V3212_V3322} and \ref{fig:V2322_V2422} below give eggbox diagrams of various linear sandwich semigroups.%
Next, we show that stability of $S$ entails stability of all sandwich semigroups~$\Sija$.

\ms
\begin{prop}\label{prop:stabilityMija}
Let $(S,\cdot,I,\lam,\rho)$ be a stable partial semigroup.  Then $\Sija$ is stable for all $i,j\in I$ and $a\in S_{ji}$.
\end{prop}

\pf Let $x,y\in S_{ij}$.  We must show that 
\[
x\gJa x\star_ay \ \iff\ x\gRa x\star_ay \AND x\gJa y\star_ax\ \iff\ x\gLa y\star_ax.
\]
By duality, it suffices to prove the first of these.  Clearly, $x\gRa x\star_ay \ \implies\ x\gJa x\star_ay$.  Conversely, suppose $x\gJa x\star_ay$.  Then one of the following holds:
\vspace{-0.2cm}
\begin{multicols}2
\begin{enumerate}
\item[(i)] $x=xay$,
\item[(ii)] $x=xayav$ for some $v\in S_{ij}$,
\item[(iii)] $x=uaxay$ for some $u\in S_{ij}$,
\item[(iv)] $x=uaxayav$ for some $u,v\in S_{ij}$.
\end{enumerate}
\end{multicols}
Clearly, (i) or (ii) implies $x\gRa xay$.  Next, suppose (iv) holds.  Then $x\J xaya$, so that $x\R xaya$ by stability.  In particular, (a) $x=xaya$ or (b) $x=xayaw$ for some $w\in S_{ij}$.  If (a) holds, then $x=(xaya)aya$, so (b) holds with $w=aya$.  In particular, $x=(x\star_ay)\star_aw$, completing the proof that $x\gRa x\star_ay$.  Finally, if (iii) holds, then $x=ua(uaxay)ay$, so that case (iii) reduces to case (iv).  The proof is therefore complete.~\epf



We conclude this section with a result that shows how regularity of the sandwich element implies close relationships between certain sandwich semigroups $\Sija$ and $S_{ji}^b$ and certain (non-sandwich) subsemigroups of $S_i$ and $S_j$.




\ms
\begin{thm}\label{thm:diamonds}
Let $(S,\cdot,I,\lam,\rho)$ be a partial semigroup and let $i,j\in I$.  Let $a\in S_{ji}$ and $b\in S_{ij}$ be such that $a=aba$ and $b=bab$.  Then 
\bit
\itemit{i} $S_{ij}a$ and $aS_{ij}$ are subsemigroups of $S_i$  and $S_j$ (respectively),
\itemit{ii} $(aS_{ij}a,\star_b)$ and $(bS_{ji}b,\star_a)$ are monoids with identities $b$ and $a$ (respectively), and are subsemigroups of $S_{ji}^b$ and $\Sija$ (respectively),
\itemit{iii} the maps $aS_{ij}a\to bS_{ji}b:x\mt bxb$ and $bS_{ji}b\to aS_{ij}a:x\mt axa$ define mutually inverse isomorphisms between $(aS_{ij}a,\star_b)$ and $(bS_{ji}b,\star_a)$,
\itemit{iv} $a\RegSija a$ is contained in $\Reg(S_{ji}^b)$, 
\itemit{v} the following diagrams commute, with all maps being homomorphisms:
\[
\includegraphics{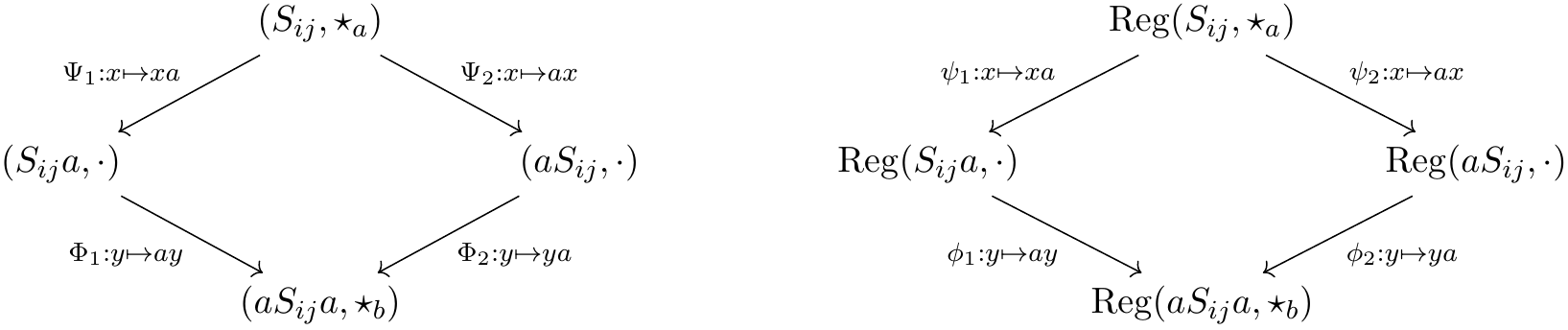}
\]
\eit
\end{thm}

\pf Part (i) is clear, and parts (ii) and (iii) are easily checked.  Next, suppose $x\in\RegSija$, so $x=xauax$ for some $u\in S_{ij}$.  Then
$
axa = axauaxa = axabauabaxa = (axa)\star_b(aua)\star_b(axa),
$
giving (iv).  Part (v) is all mostly easy to check.  That $\Phi_1$ is a homomorphism follows from
$\Phi_1((xa)(ya))=axaya = axabaya = \Phi_1(xa)\star_b\Phi_1(ya)$.
It is clear that $\psi_1$ maps $\RegSija$ into $\Reg(S_{ij}a)$.  
It follows from (iv) that $\phi_1$ maps $\Reg(S_{ij}a)$ into $\Reg(S_{ji}^b)$. \epf


\begin{rem}
Other relationships exist, such as $(baS_{ij}a,\cdot) = (bS_{ji}ba,\cdot)$, but these will not be explored any further.
\end{rem}

\section{The linear partial semigroup}\label{sect:preliminaries}

As noted earlier, to understand the linear sandwich semigroups $\MmnA$, it is crucial to first understand the partial semigroup $\M$.  So in this section, we gather the required material on $\M$, showing how the general framework of Section \ref{sect:partial} applies in this case.

We fix a field $\F$ for the remainder of the article.  For positive integers $m,n\in\mathbb N=\{1,2,3,\ldots\}$, we write $\M_{mn}=\Mmn(\F)$ for the set of all $m\times n$ matrices (i.e., all matrices with $m$ rows and $n$ columns) over $\F$.  We write $\M=\M(\F)=\bigcup_{m,n\in\mathbb N}\Mmn$ for the set of all (finite dimensional, non-empty) matrices over $\F$.  So $\M$ is a \emph{partial semigroup}, as noted in Example \ref{eg:M}.
By convention, we consider there to be a unique $m\times0$ and $0\times n$ matrix for any $m,n\geq0$, namely the empty matrix, which we denote by $\emptyset$.  So $\M_{mn}=\{\emptyset\}$ if $m=0$ or $n=0$.  But this is a matter of convenience, and we do not consider the empty matrix $\emptyset$ to be an element of $\M$.
%
We also write $\M_n=\M_n(\F)=\M_{nn}$ for any $n$%
, and denote by $\G_n=\G_n(\F)$ the group of $n\times n$ invertible matrices over $\F$.  So $\M_n$ and $\G_n$ are the \emph{full linear monoid} and \emph{general linear group} of degree $n$.  For background on the full linear monoids, the monograph \cite{Okninski1998} is highly recommended.

If $V$ and $W$ are vector spaces, we write $\Hom(V,W)$ for the set of all linear transformations from $V$ to~$W$.  As usual, if $\al\in\Hom(V,W)$, we write $\im(\al)=\set{\al(v)}{v\in V}$ and $\ker(\al)=\set{v\in V}{\al(v)=0}$ for the \emph{image} and \emph{kernel} of $\al$.  We write $\End(V)=\Hom(V,V)$ for the monoid of all endomorphisms of $V$ (i.e., all linear transformations $V\to V$), and $\Aut(V)$ for the group of all automorphisms of $V$ (i.e., all invertible endomorphisms of $V$).  For $n\geq0$, we write $V_n=\F^n$ for the vector space of all $n\times1$ column vectors over~$\F$.  We will identify $\M_{mn}$ with $\Homnm$ in the usual way.  Namely, if $X\in\Mmn$, we write $\lam_X\in\Homnm$ for the linear transformation $\lam_X:V_n\to V_m$ defined by $\lam_X(v)=Xv$ for all $v\in V_n$.  We will often prove statements about $\M_{mn}$ by proving the equivalent statement about $\Homnm$.  When $m=n$, the map $X\to\lam_X$ determines an isomorphism of monoids $\M_n\to\Endn$, and its restriction to $\G_n\sub\M_n$ determines an isomorphism of groups $\G_n\to\Autn$.  We write $\{e_{n1},\ldots,e_{nn}\}$ for the standard basis of $V_n$ ($e_{ni}$ has a $1$ in position $i$ and $0$'s elsewhere).  We also write $W_{ns}=\Span\{e_{n1},\ldots,e_{ns}\}$ for each $0\leq s\leq n$.  (We interpret $\Span\emptyset=\{0\}$, though the dimension of the ambient space must be understood from context.)

%
Our first aim is to characterise Green's relations ($\R$, $\L$, $\J$, $\D$, $\H$) and preorders ($\leq_\R$, $\leq_\L$, $\leq_\J$) on $\M$.  Because $\M$ is monoidal (see Definition~\ref{defn:S}), $\M=\M^{(1)}$.  So, for example, if $X,Y\in\M$ are two matrices (not necessarily of the same size), then $X\leq_\R Y$ if and only if $X=YA$ for some $A\in\M$.
Note that if $X\leq_\R Y$ (resp., $X\leq_\L Y$), then $X$ and $Y$ must have the same number of rows (resp., columns).  
%

Let $X\in\M_{mn}$.  For $1\leq i\leq m$ and $1\leq j\leq n$, we write $\row_i(X)$ and $\col_j(X)$ for the $i$th row and $j$th column of~$X$, respectively.  We write $\Row(X)=\Span\{\row_1(X),\ldots,\row_m(X)\}$ and $\Col(X)=\Span\{\col_1(X),\ldots,\col_n(X)\}$ for the \emph{row space} and \emph{column space} of $X$, respectively, and we write $\rank(X)=\dim(\Row(X))=\dim(\Col(X))$ for the \emph{rank} of $X$.  
Because of the transpose map $\M\to\M:A\mt A^\tr$, which is a bijection and satisfies $(AB)^\tr=B^\tr A^\tr$, the linear partial semigroup $\M$ is \emph{self-dual} (in the sense that it is anti-isomorphic to its own dual).  
Since $\Row(X^\tr)=\Col(X)$, any statement about row spaces implies a corresponding dual statement about column spaces (and vice versa).  (Without causing confusion, we will often blur the distinction between row vectors and column vectors, and think of $\Row(X)$ and $\Col(X)$ as subspaces of $\F^n$ and $\F^m$, respectively.)

The next result characterises Green's relations and preorders on $\M$ in terms of the parameters introduced above.  An equivalent formulation in the special case of square matrices may be found in \cite[Lemma 2.1]{Okninski1998}.



\ms
\begin{lemma}\label{lem:green<M}
Let $X,Y\in\M$.  Then
\begin{itemize}\begin{multicols}2
\itemit{i} $X\leq_\R Y \iff \Col(X)\sub\Col(Y)$, 
\itemit{ii} $X\leq_\L Y \iff \Row(X)\sub\Row(Y)$,
\itemit{iii} $X\leq_\J Y \iff \rank(X)\sub\rank(Y)$,
\itemit{iv} $X\R Y \tiff \Col(X)=\Col(Y)$,
\itemit{v} $X\L Y \tiff \Row(X)=\Row(Y)$,
\itemit{vi} $X\J Y \tiff \rank(X)=\rank(Y)$. 
\end{multicols} \end{itemize} \nss\ss
Further, $\M$ is stable, so $\J=\D$.
\end{lemma}

\pf Clearly, (iv--vi) follow from (i--iii).  Note that (ii) is the dual of (i), which is true because
\begin{align*}
X\leq_\R Y &\ \iff\ \text{$X=YA$ for some $A\in\M$}\\
&\ \iff\ \text{every column of $X$ is a linear combination of the columns of $Y$}\\
&\ \iff\ \Col(X)\sub\Col(Y).
\end{align*}
For (iii), 
if $X\leq_\J Y$, then $X=AYB$ for some $A,B\in\M$, giving $\rank(X)=\rank(AYB)\leq\rank(Y)$.  Conversely, suppose $\rank(X)\leq\rank(Y)$, and say $X\in\Mmn$ and $Y\in\Mkl$.  It is sufficient to show that $\lam_X=\al\circ\lam_Y\circ\be$ for some $\al\in\Homkm$ and $\be\in\Homnl$.  
Put $r=\rank(X)$ and $s=\rank(Y)$.  Choose bases $\B_1=\{u_1,\ldots,u_n\}$ and $\B_2=\{v_1,\ldots,v_l\}$ for $V_n$ and $V_l$ so that $\{u_{r+1},\ldots,u_n\}$ and $\{v_{s+1},\ldots,v_l\}$ are bases for $\ker(\lam_X)$ and $\ker(\lam_Y)$, respectively.  Extend (if necessary) the linearly independent sets $\{\lam_Y(v_1),\ldots,\lam_Y(v_r)\}$ and $\{\lam_X(u_1,\ldots,\lam_X(u_r)\}$ arbitrarily to bases
\[
\B_3=\{\lam_Y(v_1),\ldots,\lam_Y(v_r),w_{r+1},\ldots,w_k\} \AND \B_4=\{\lam_X(u_1),\ldots,\lam_X(u_r),x_{r+1},\ldots,x_m\}
\]
for $V_k$ and $V_m$.  Now let $\al\in\Homkm$ and $\be\in\Homnl$ be chosen arbitrarily so that
\[
\begin{array}{rclcrcll}
\al(\lam_Y(v_i)) \hspace{-.25cm}&=&\hspace{-.25cm} \lam_X(u_i), & &
\al(w_j) \hspace{-.25cm} &\in&\hspace{-.25cm} \Span\{x_{r+1},\ldots,x_m\}
&\qquad\text{for all $1\leq i\leq r$ and $r+1\leq j\leq k$,}
\\
\be(u_i) \hspace{-.25cm}&=&\hspace{-.25cm} v_i, & &
\be(u_j) \hspace{-.25cm} &\in&\hspace{-.25cm} \Span\{v_{s+1},\ldots,v_l\}
&\qquad\text{for all $1\leq i\leq r$ and $r+1\leq j\leq n$.}
\end{array}
\]
One easily checks that $\al\circ\lam_Y\circ\be=\lam_X$, by checking the respective actions on the basis $\B_1$ of $V_n$. 

To prove stability, we must show that for all $X,Y\in\M$,
\[
X\J XY \tiff X\R XY \AND X\J YX \tiff X\L YX.
\]
By duality, it suffices to prove the first equivalence.  Since $\R\sub\J$, it is enough to prove that $X\J XY \implies X\R XY$.  Now, $\Col(XY)\sub\Col(X)$.  But also $X\J XY$ gives $\dim(\Col(X))=\rank(X)=\rank(XY)=\dim(\Col(XY))$, so that $\Col(X)=\Col(XY)$, whence $X\R XY$. \epf


As we saw in Section \ref{sect:partial}, stability and regularity are very useful properties for a partial semigroup to have.  Now that we know $\M$ is stable, let us show that $\M$ is also regular.  




\ms
\begin{lemma}\label{lem:regularity}
The linear partial semigroup $\M$ is regular.  
\end{lemma}

\pf Let $X\in\M_{mn}$.  It suffices to show that there exists $\al\in\Hommn$ such that ${\lam_X=\lam_X\circ\al\circ\lam_X}$.  Let $\B=\{v_1,\ldots,v_n\}$ be a basis of $V_n$ such that $\{v_{r+1},\ldots,v_n\}$ is a basis of $\ker(\lam_X)$.  Extend (if necessary) the linearly independent set $\{\lam_X(v_1),\ldots,\lam_X(v_r)\}$ to a basis $\{\lam_X(v_1),\ldots,\lam_X(v_r),w_{r+1},\ldots,w_m\}$ of $V_m$.  Let $\al\in\Hommn$ be any linear transformation for which $\al(\lam_X(v_i))=v_i$ for each $1\leq i\leq r$.  Then one easily checks that $\lam_X=\lam_X\circ\al\circ\lam_X$ by calculating the action on the basis $\B$. \epf






As in Section \ref{sect:partial}, if $X\in\Mmn$ and $\gK$ is one of $\R$, $\L$, $\J$, $\D$, $\H$, we write $K_X=\set{Y\in\Mmn}{X\gK Y}$, and call $K_X$ the \emph{$\gK$-class} of $X$ in $\Mmn$.  Note that all matrices from $K_X$ have the same dimensions.  (We will have no need to consider the sets $\XK$ of \emph{all} matrices $\gK$-related to $X$.)  Recall that $\G_k$ denotes the group of all invertible $k\times k$ matrices over $\F$.  The next result gives an alternative description of various Green's classes in $\M$.


\ms
\begin{lemma}\label{lem:greenMmn}
Let $X\in\M_{mn}$.  Then
\bit
\itemit{i} $R_X  = \set{Y\in\M_{mn}}{\Col(X)=\Col(Y)} = X\G_n$,
\itemit{ii} $L_X  = \set{Y\in\M_{mn}}{\Row(X)=\Row(Y)} = \G_mX$,
\itemit{iii} $D_X = J_X = \set{Y\in\M_{mn}}{\rank(X)=\rank(Y)} = \G_mX\G_n$.
\eitres
\end{lemma}

\pf For (i), 
note that clearly $X\G_n\sub R_X$.  By Lemma \ref{lem:green<M}, it remains to show the reverse inclusion, so suppose $Y\in R_X$.  In particular, $X\J Y$, so $\rank(X)=\rank(Y)$.  Put $r=\rank(X)$.  We show that $\lam_Y=\lam_X\circ\al$ for some $\al\in\Autn$.  Since $X\R Y$, we already know that $\lam_Y=\lam_X\circ\be$ for some $\be\in\Endn$.  Let $\B_1=\{u_1,\ldots,u_n\}$ be a basis of $V_n$ such that $\{u_{r+1},\ldots,u_n\}$ is a basis of $\ker(\lam_Y)$.  So $\{\lam_X(\be(u_1)),\ldots,\lam_X(\be(u_r))\}=\{\lam_Y(u_1),\ldots,\lam_Y(u_r)\}$ is a basis of $\im(\lam_Y)$.  It follows that $\{\be(u_1),\ldots,\be(u_r)\}$ is linearly independent.  We may therefore extend this set to a basis $\B_2=\{\be(u_1),\ldots,\be(u_r),v_{r+1},\ldots,v_n\}$ of $V_n$, where $\{v_{r+1},\ldots,v_n\}$ is a basis of $\ker(\lam_X)$.  Now define $\al\in\Autn$ by
\[
\al(u_i) = \begin{cases}
\be(u_i) &\text{if $1\leq i\leq r$}\\
v_i &\text{if $r<i\leq n$.}
\end{cases}
\]
One easily checks that $\lam_Y=\lam_X\circ\al$.  This completes the proof of (i).  

Part (ii) is dual to (i).  
For (iii), 
clearly $\G_mX\G_n\sub J_X$, and the converse follows quickly from (i) and~(ii) and the fact that $\J=\D=\L\circ\R$. By Lemma \ref{lem:green<M}, this completes the proof. \epf

If $\gK$ is one of $\R$, $\L$, $\D=\J$, then the set $\Mmn/\gK$ of all $\gK$-classes of $\Mmn$ inherits a partial order:
\[
K_X \leq_\gK K_Y \ \iff \ X\leq_\gK Y.
\]
We typically write $\leq$ for the order $\leq_\J$ on the $\D=\J$-classes.  Of importance is the fact that these classes form a chain:
\[
\DMmn0<\DMmn1<\cdots<\DMmn l,
\]
where $\DMmn s=\set{X\in\Mmn}{\rank(X)=s}$ for all $0\leq s\leq l=\min(m,n)$.

Figure \ref{fig:eggbox} pictures an eggbox diagram (as explained in Section \ref{sect:partial}) of the $\D$-class $D_1(\M_{23}(\mathbb Z_3))$ of all $2\times3$ matrices of rank $1$ over the field $\F=\mathbb Z_3=\{0,1,2\}$
(see Lemma \ref{lem:combinatorics_Mmn} for an explanation of the number and sizes of the $\R$-, $\L$- and $\H$-classes).  The reader need not yet worry about the subdivisions within the eggbox; for now, it is enough to note that the matrices to the left (resp., top) of the vertical (resp., horizontal) divider satisfy the property that the first column (resp., row) spans the column space (resp., row space) of the matrix.

\begin{figure}[ht]
\begin{center}
\scalebox{.8}{
\begin{tikzpicture}[scale=0.9]
\dmat1{7.9}100000
\dmat1{7.1}200000
\dmat1{5.9}100100
\dmat1{5.1}200200
\dmat1{3.9}100200
\dmat1{3.1}200100
\dmat1{1.9}000100
\dmat1{1.1}000200
\dmat2{7.9}101000
\dmat2{7.1}202000
\dmat2{5.9}101101
\dmat2{5.1}202202
\dmat2{3.9}101202
\dmat2{3.1}202101
\dmat2{1.9}000101
\dmat2{1.1}000202
\dmat3{7.9}102000
\dmat3{7.1}201000
\dmat3{5.9}102102
\dmat3{5.1}201201
\dmat3{3.9}102201
\dmat3{3.1}201102
\dmat3{1.9}000102
\dmat3{1.1}000201
\dmat4{7.9}110000
\dmat4{7.1}220000
\dmat4{5.9}110110
\dmat4{5.1}220220
\dmat4{3.9}110220
\dmat4{3.1}220110
\dmat4{1.9}000110
\dmat4{1.1}000220
\dmat5{7.9}111000
\dmat5{7.1}222000
\dmat5{5.9}111111
\dmat5{5.1}222222
\dmat5{3.9}111222
\dmat5{3.1}222111
\dmat5{1.9}000111
\dmat5{1.1}000222
\dmat6{7.9}112000
\dmat6{7.1}221000
\dmat6{5.9}112112
\dmat6{5.1}221221
\dmat6{3.9}112221
\dmat6{3.1}221112
\dmat6{1.9}000112
\dmat6{1.1}000221
\dmat7{7.9}120000
\dmat7{7.1}210000
\dmat7{5.9}120120
\dmat7{5.1}210210
\dmat7{3.9}120210
\dmat7{3.1}210120
\dmat7{1.9}000120
\dmat7{1.1}000210
\dmat8{7.9}121000
\dmat8{7.1}212000
\dmat8{5.9}121121
\dmat8{5.1}212212
\dmat8{3.9}121212
\dmat8{3.1}212121
\dmat8{1.9}000121
\dmat8{1.1}000212
\dmat9{7.9}122000
\dmat9{7.1}211000
\dmat9{5.9}122122
\dmat9{5.1}211211
\dmat9{3.9}122211
\dmat9{3.1}211122
\dmat9{1.9}000122
\dmat9{1.1}000211
\dmat{10}{7.9}010000
\dmat{10}{7.1}020000
\dmat{10}{5.9}010010
\dmat{10}{5.1}020020
\dmat{10}{3.9}010020
\dmat{10}{3.1}020010
\dmat{10}{1.9}000010
\dmat{10}{1.1}000020
\dmat{11}{7.9}011000
\dmat{11}{7.1}022000
\dmat{11}{5.9}011011
\dmat{11}{5.1}022022
\dmat{11}{3.9}011022
\dmat{11}{3.1}022011
\dmat{11}{1.9}000011
\dmat{11}{1.1}000022
\dmat{12}{7.9}012000
\dmat{12}{7.1}021000
\dmat{12}{5.9}012012
\dmat{12}{5.1}021021
\dmat{12}{3.9}012021
\dmat{12}{3.1}021012
\dmat{12}{1.9}000012
\dmat{12}{1.1}000021
\dmat{13}{7.9}001000
\dmat{13}{7.1}002000
\dmat{13}{5.9}001001
\dmat{13}{5.1}002002
\dmat{13}{3.9}001002
\dmat{13}{3.1}002001
\dmat{13}{1.9}000001
\dmat{13}{1.1}000002
%
\foreach \x in {0,2,4,6,8}
\draw (0.75,\x+.5)--(20.25,\x+.5); 
\foreach \x in {0,...,13}
\draw (\x*1.5+0.75,.5)--(\x*1.5+0.75,8.5); 
\draw[line width=.5mm] (0.75,0.5)--(20.25,0.5)--(20.25,8.5)--(0.75,8.5)--(0.75,0.5)--(20.25,0.5);
\draw[line width=.5mm] (0.75,2.5)--(20.25,2.5);
\draw[line width=.5mm] (14.25,8.5)--(14.25,0.5);
\draw[|-|] (0.75,9.0)--(14.25,9.0);
\draw[|-|] (20.25,9.0)--(14.25,9.0);
\draw[|-|] (.25,.5)--(.25,2.5);
\draw[|-|] (.25,8.5)--(.25,2.5);
\draw(7.5,9.3)node{$\sub P_1\phantom{\sub}$};
\draw(17.25,9.3)node{$\not\sub P_1\phantom{\sub}$};
\draw(.25,5.5)node[left]{$\sub P_2\phantom{}$};
\draw(.25,1.5)node[left]{$\not\sub P_2\phantom{}$};
\end{tikzpicture}
}
    \caption{An eggbox diagram of the $\D$-class $D_1(\M_{23}(\mathbb Z_3))$.}
    \label{fig:eggbox}
   \end{center}
 \end{figure}

%
%
%
%

So $\Mmn$ has $\min(m,n)+1$ $\D$-classes.
It will also be convenient to have some more combinatorial information about the number and size of certain $\gK$-classes.  Recall that the $q$-factorials and $q$-binomial coefficients are defined by
\[
\qfact s = 1\cdot(1+q)\cdots(1+q+\cdots+q^{s-1}) = \frac{(q-1)(q^2-1)\cdots(q^s-1)}{(q-1)^s}
\]
and
\[
\qbin ms = \frac{\qfact m}{\qfact s \qfact{m-s}} = \frac{(q^m-1)(q^m-q)\cdots(q^m-q^{s-1})}{(q^s-1)(q^s-q)\cdots(q^s-q^{s-1})} = \frac{(q^m-1)(q^{m-1}-1)\cdots(q^{m-s+1}-1)}{(q^s-1)(q^{s-1}-1)\cdots(q-1)}.
\]
It is easy to check (and well-known) that when $|\F|=q<\infty$,
\[
|\G_s| = (q^s-1)(q^s-q)\cdots(q^s-q^{s-1}) = q^{{s\choose2}}(q-1)^s \qfact s.
\]
In what follows, a crucial role will be played by the matrices $J_{mns}\in\Mmn$ defined for $s\leq\min(m,n)$ by
\[
J_{mns} = \mat{I_s}{O_{s,n-s}}{O_{m-s,s}}{O_{m-s,n-s}} .
\]
Here and elsewhere, we write $I_s$ and $O_{kl}$ for the $s\times s$ identity matrix and $k\times l$ zero matrix (respectively).  If the dimensions are understood from context, we just write $O=O_{kl}$.  
So $J_{mns}$ is the $m\times n$ matrix with $1$'s in the first $s$ positions on the leading diagonal and $0$'s elsewhere.  Note that if $s=m\leq n$ (resp., $s=n\leq m$), then the matrices $O_{m-s,s}$ and $O_{m-s,n-s}$ (resp., $O_{s,n-s}$ and $O_{m-s,n-s}$) are empty, and $J_{mns}=[I_s\ O_{s,n-s}]$ (resp., $\thmat{I_s}{O_{m-s,s}}$).


\ms
\begin{lemma}\label{lem:combinatorics_Mmn}
Suppose $|\F|=q<\infty$, and let $0\leq s\leq \min(m,n)$.  Then
\bit
\itemit{i} $\DMmn s$ contains $\tqbin ms$ $\R$-classes,
\itemit{ii} $\DMmn s$ contains $\tqbin ns$ $\L$-classes,
\itemit{iii} $\DMmn s$ contains $\tqbin ms\tqbin ns$ $\H$-classes, each of which has size $|\G_s|=q^{{s\choose2}}(q-1)^s\qfact s$,
\itemit{iv} $|D_s(\Mmn)|=\tqbin ms\tqbin nsq^{{s\choose2}}(q-1)^s\qfact s$.
\eitres
\end{lemma}

\pf Parts (i) and (ii) follow immediately from parts (i) and (ii) of Lemma \ref{lem:greenMmn} and the well-known fact that $\tqbin ms$ is the number of $s$ dimensional subspaces of an $m$ dimensional vector space over $\F$.  The number of $\H$-classes follows immediately from (i) and (ii).  By Lemma \ref{lem:GreensLemma_Sij}, all the $\H$-classes in $\DMmn s$ have the same size, so it suffices to calculate the size of $H=H_{J_{mns}}$.
Let $X=\tmat ABCD\in H$, where $A\in\M_s$, $B\in\M_{s,n-s}$, and so on.  Since $\Row(X)=\Row(J_{mns})$, we see that $B$ and $D$ are zero matrices.  Considering column spaces, we see that~$C$ is also a zero matrix.  It follows that $X=\tmat AOOO$, and also $\rank(A)=\rank(X)=\rank(J_{mns})=s$.  Clearly every such matrix $X=\tmat AOOO$ with $\rank(A)=s$ belongs to $H$.  The condition that $\rank(A)=s$ is equivalent to $A\in\G_s$, so it follows that~$|H|=|\G_s|$.  Finally, (iv) follows from (iii).~\epf

Of course, by considering the size of $\Mmn$ when $|\F|=q<\infty$, we obtain the identity
\[
q^{mn} = \sum_{s=0}^l \qbin ms\qbin ns \qfact s(q-1)^sq^{{s\choose2}}.
\]
We conclude this section by stating some well-known results on the full linear monoids $\M_n$ and their ideals that we will require in what follows.
The set $E(\M_n)=\set{X\in\M_n}{X=X^2}$ of idempotents of~$\M_n$ is not a subsemigroup (unless $n\leq1$), but the 
subsemigroup $\E_n=\la E(\M_n)\ra$ of $\M_n$ generated by these idempotents has a neat description.  Namely, it was shown by Erdos \cite{Erdos1967} that any singular (i.e., non-invertible) matrix over~$\F$ is a product of idempotent matrices.  This result has been reproved by a number of authors \cite{FL1992,Djokovic1968,Dawlings81/82,AM2005,Laffey1983}.  The minimal number of (idempotent) matrices required to generate~$\E_n$ was determined by Dawlings \cite{Dawlings1982}.  Recall that the \emph{rank} (resp., \emph{idempotent rank}) of a semigroup (resp., idempotent generated semigroup)~$S$, denoted $\rank(S)$ (resp., $\idrank(S)$), is the minimal size of a generating set (resp., idempotent generating set) for $S$.  (The rank of a semigroup should not be confused with the rank of a matix.)  If $U$ is a subset of a semigroup $S$, we write $E(U)=E(S)\cap U$ for the set of all idempotents from $U$.



\ms
\begin{thm}[Erdos \cite{Erdos1967}, Dawlings \cite{Dawlings81/82,Dawlings1982}]\label{thm_MnGn}
We have
\[
\E_n=\la E(\M_n)\ra = (\MnGn)\cup\{I_n\} \AND \MnGn= \la E(D_{n-1}(\M_n)) \ra .
\]
Further, if $|\F|=q<\infty$, then
\[
\epfreseq
\rank(\MnGn)=\idrank(\MnGn)=(q^n-1)/(q-1).
\]
\end{thm}

The previous result has been extended by Gray \cite{Gray2008} to arbitrary ideals of $\M_n$.

\ms
\begin{thm}[Gray \cite{Gray2008}]\label{thm_ideals_Mn}
The ideals of $\M_n$ are precisely the sets
\[
I_s(\M_n)=D_0(\M_n)\cup\cdots\cup D_s(\M_n)=\set{X\in\M_n}{\rank(X)\leq s} \qquad\text{for $0\leq s\leq n$,}
\]
and they form a chain: $I_0(\M_n)\sub\cdots\sub I_n(\M_n)$.  If $0\leq s<n$, then $I_s(\M_n) = \la E(D_s(\M_n)) \ra$ is generated by the idempotents in its top $\D$-class.  Further, if $|\F|=q<\infty$, then
\[\epfreseq
\rank(I_s(\M_n))=\idrank(I_s(\M_n))=\qbin ns.
\]
\end{thm}

Note that $I_n(\M_n)=\M_n$, $D_n(\M_n)=\G_n$ and $I_{n-1}(\M_n)=\MnGn$, so Theorem \ref{thm_MnGn} is a special case of Theorem \ref{thm_ideals_Mn} since $\tqbin n{n-1}=(q^n-1)/(q-1)$.
%
%

On several occasions, we will need to make use of the fact that the general linear group $\G_n$ may be generated by two matrices, as was originally proved by Waterhouse \cite{Waterhouse}; see also \cite{Gill2015}, where minimal generating sets for $\G_n$ are explored in more detail.
%
Probabilistic generation of matrix groups is considered in \cite{GK2000,BGK2008}, for example, though the context is usually for classical groups.

\ms
\begin{thm}[Waterhouse \cite{Waterhouse}]\label{thm:waterhouse}
If $|\F|<\infty$, then
\bit
\itemit{i} $\rank(\G_1)=1$, and $\rank(\G_n)=2$ if $n\geq2$,
\itemit{ii} $\M_n=\la\G_n\cup\{X\}\ra$ for any $X\in D_{n-1}(\M_n)$,
\itemit{iii} $\rank(\M_1)=2$, and $\rank(\M_n)=3$ if $n\geq2$. \epfres
\eitres
\end{thm}


For convenience, eggbox diagrams are given for the full linear monoids $\M_n(\mathbb Z_2)$ for $0\leq n\leq3$ in Figure \ref{fig:M0...M3} below.  In the diagrams, group $\gH$-classes are shaded grey, and a label of {\tt k} indicates that the group $\gH$-class is isomorphic to $\G_k(\mathbb Z_2)$.

\begin{figure}[ht]
\begin{center}
\includegraphics[width=7.3mm]{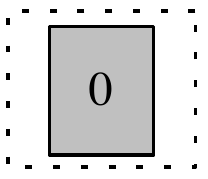}
\qquad
\includegraphics[width=7.5mm]{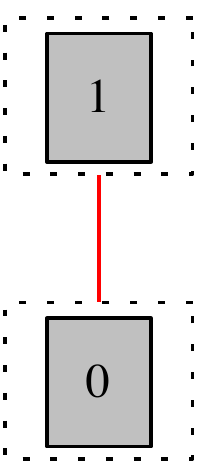}
\qquad
\includegraphics[width=13.5mm]{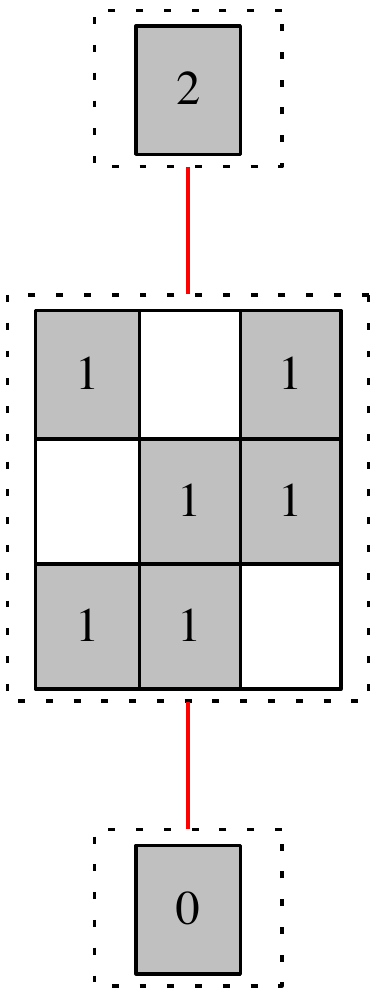}
\qquad
\includegraphics[width=21mm]{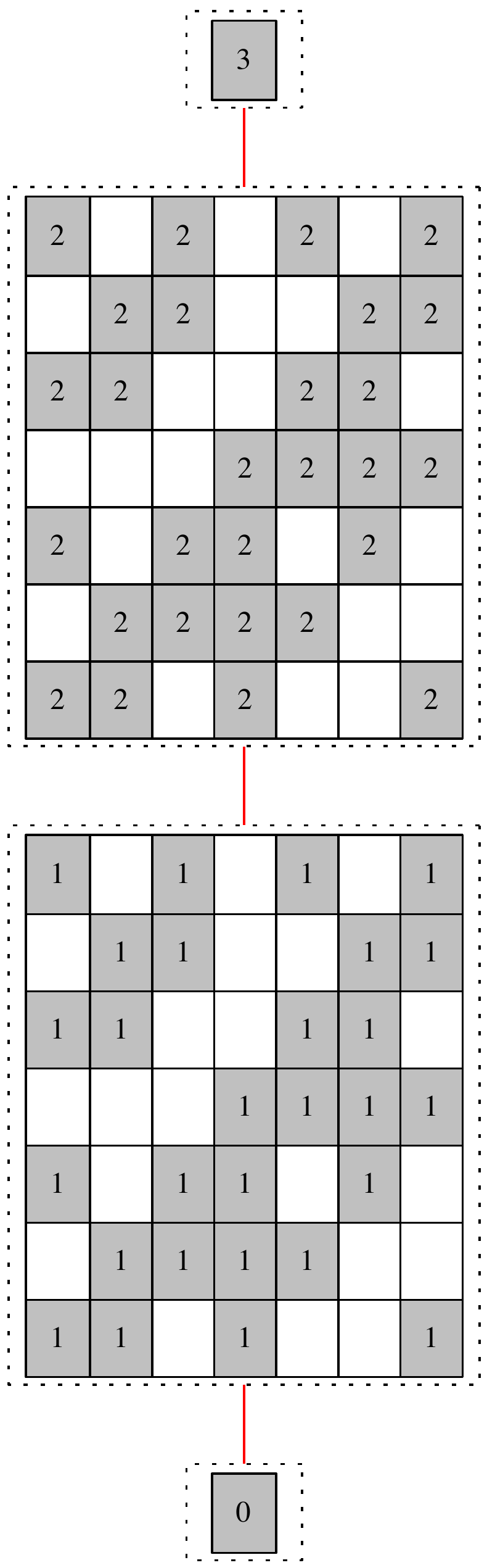}
    \caption{Egg box diagrams of the full linear semigroups $\M_0$, $\M_1$, $\M_2$, $\M_3$, all over $\mathbb Z_2$ (left to right).}
    \label{fig:M0...M3}
   \end{center}
 \end{figure}

\section{Linear sandwich semigroups}\label{sect:MmnJ}

Now that we have gathered the required material on $\M$, we may begin our study of
the linear sandwich semigroups.  From now on, we fix integers $m,n\geq1$ and an $n\times m$ matrix $A\in\M_{nm}$.  As in Section \ref{sect:partial}, we denote by $$\MmnA=\MmnA(\F)=(\Mmn,\star_A)$$ the sandwich semigroup of $\Mmn$  under the operation $\star_A$ defined by
\[
X\star_A Y = XAY \qquad\text{for $X,Y\in\M_{mn}$.}
\]
We note that if $m=n$, then $\MmnA=\M_n^A$ is a \emph{variant} \cite{Hickey1983} of the full linear monoid~$\M_n$, so everything we prove about linear sandwich semigroups holds for such linear variants also.  
%
%
We begin with a simple observation.

\ms\ms
\begin{lemma}\label{lem:MmnAMmnB}
\bit
\itemit{i} If $A\in\M_{nm}$, then $\M_{mn}^A\cong\M_{nm}^{A^\tr}$.
\itemit{ii} If $A,B\in\M_{nm}$ are such that $\rank(A)=\rank(B)$, then $\MmnA\cong\MmnB$.
\eitres
\end{lemma}

\pf It clear that $X\mt X^\tr$ defines an isomorphism $\M_{mn}^A\to\M_{nm}^{A^T}$, giving (i).  Next, if $\rank(A)=\rank(B)$, Lemma \ref{lem:greenMmn} gives $A=UBV$ for some $U\in\G_m$ and $V\in\G_n$.  But then one may check that $X\mt VXU$ defines an isomorphism $\MmnA\to\MmnB$, giving (ii). \epf

In particular, when studying the semigroup $\MmnA$ where $\rank(A)=r$, 
we may choose any $A\in\Mnm$ of rank~$r$.  For the rest of the article, we will therefore study the semigroup $\MmnJ$, where 
\[
J=J_{nmr} = \mat{I_r}{O_{r,m-r}}{O_{n-r,r}}{O_{n-r,m-r}}\in\M_{nm}.
\]
From now on, unless otherwise specified, whenever a $k\times l$ matrix $X$ (with $k,l\in\{m,n\}$) is written in $2\times2$ block form, $X=\tmat ABCD$, we will be tacitly assuming that $A\in\M_r$ (from which the dimensions of $B,C,D$ may be deduced).  So for example, we will usually just write $J=\tmat IOOO$.
For simplicity, we will write $\star$ for the operation $\star_J$ on $\MmnJ$, throughout.  One easily verifies the rule
\[
\mat ABCD \star \mat EFGH = \mat{AE}{AF}{CE}{CF}.
\]
Also note that if $X=\tmat ABCD$, then
\[
XJ=\mat AOCO \in\M_m \COMMA JX=\mat ABOO\in\M_n \COMMA JXJ = \mat AOOO\in\M_{nm}.
\]

\ms
\begin{rem}\label{rem:r=m}
In the special case that $r=m\leq n$, we have $J=[I\ O]$, and the product in $\MmnJ$ satisfies $[A\ B] \star[E\ F] = [AE\ AF]$.  But we just view this as a special case of the above rule, with the bottom rows --- i.e., $[C\ D]$, $[G\ H]$, $[CE\ CF]$ --- containing empty blocks.  A dual statement holds in the case $r=n\leq m$.  
In only one place will we need to consider the case in which $r=\min(m,n)$ separately (see Theorems~\ref{thm:rankMmnJ} and~\ref{thm:rankMmnJ_r=m}).  If $r=m=n$, then $\MmnJ$ is precisely the full linear monoid $\M_n$; since all the problems we investigate have already been solved for $\M_n$, we will typically assume that $r=m=n$ does not hold, though our results are true for the case $r<m=n$ (corresponding to \emph{variants} of the full linear monoids~$\M_n$).
See Remark~\ref{rem:r=m2}, where the above observations are used to show that the sandwich semigroups $\MmnJ$ are isomorphic to certain well-known (non-sandwich) matrix semigroups in the case that $r=\min(m,n)$.
\end{rem}

Green's relations and the regular elements of the sandwich semigroup $\MmnJ$ were calculated in \cite{Chinram2009,Kemprasit2002}.  We now show how these results may be recovered (and given a cleaner presentation) using the general theory developed in Section \ref{sect:partial}.
In particular, a crucial role is played by the sets
\[
P_1^J = \set{X\in\Mmn}{XJ\R X} ,\
P_2^J = \set{X\in\Mmn}{JX\L X} ,\
P_3^J = \set{X\in\Mmn}{JXJ\J X} ,\
P^J=P_1^J\cap P_2^J.
\]
For simplicity, we denote these sets simply by $P_1$, $P_2$, $P_3$, and $P=P_1\cap P_2$.  

Certain special matrices from~$\Mmn$ will be very important in what follows.  With this in mind, if $A\in\M_r$, $M\in\M_{m-r,r}$ and $N\in\M_{r,n-r}$, we write
\[
[M,A,N] = \mat{A}{AN}{MA}{MAN} \in\Mmn.
\]
One may check that when matrices of this form are multiplied in $\MmnJ$, they obey the rule
\[
[M,A,N]\star[K,B,L] = [M,AB,L].
\]

\ms
\begin{prop}\label{prop:P1P2}
~
\bit
\itemit{i} $P_1=\set{X\in\Mmn}{\rank(XJ)=\rank(X)}= \bigset{X\in\Mmn}{\Col(XJ)=\Col(X)}$,
\itemit{ii} $P_2=\set{X\in\Mmn}{\rank(JX)=\rank(X)}= \bigset{X\in\Mmn}{\Row(JX)=\Row(X)}$,
\itemit{iii} $P_3=P=\set{X\in\Mmn}{\rank(JXJ)=\rank(X)}$ 
\item[] $\phantom{P_3=P}=\bigset{[M,A,N]}{A\in\M_r,\ M\in\M_{m-r,r},\ N\in\M_{r,n-r}}$,
\itemit{iv} $P=\RegMmnJ$ is the set of all regular elements of $\MmnJ$, and is a subsemigroup of $\MmnJ$.
\eitres
\end{prop}

\pf Parts (i) and (ii) follow quickly from Lemma \ref{lem:green<M} (making crucial use of stability).
We now prove (iii).  Since $\M$ is stable, Proposition~\ref{prop:Reg(Sija)} and Lemma \ref{lem:greenMmn} give
$
P_3 = P = \set{X\in\Mmn}{\rank(JXJ)=\rank(X)}.
$
Now let $X=\tmat ABCD\in\Mmn$.  First, note that
\begin{align*}
X\in P_2 &\ \iff\  \Row(X)=\Row(JX)=\Row\tmat ABOO \\
&\ \iff\  \text{each row of $[C\ D]$ is a linear combination of the rows of $[A\ B]$} \\
&\ \iff\  \text{$[C\ D]=M[A\ B]=[MA\ MB]$ for some $M\in\M_{m-r,r}$.}
\intertext{Similarly,}
X\in P_1 &\ \iff\  \text{$\hmat BD= \hmat ACN=\hmat{AN}{CN}$ for some $N\in\M_{r,n-r}$.}
\end{align*}
Putting these together, we see that $X\in P=P_1\cap P_2$ if and only if $P=\tmat A{AN}{MA}{MAN}=[M,A,N]$, completing the proof of (iii).


For (iv), 
Proposition \ref{prop:Reg(Sija)} gives $\RegMmnJ\sub P$.
Conversely, suppose $X=[M,A,N]\in P$.  If $B\in\M_r$ is such that $A=ABA$ (see Lemma \ref{lem:regularity}), then it is easy to check that $X=X\star Y\star X$ where $Y=\tmat BCDE$ for any (appropriately sized) $C,D,E$, completing the proof that $P=\RegMmnJ$.  The fact that $P$ is a subsemigroup follows immediately from Proposition \ref{prop:regularSija} and Lemma \ref{lem:regularity} (or directly from the rule $[M,A,N]\star[K,B,L]=[M,AB,L]$).~\epf

\begin{rem}
Part (iv) of the previous proposition also follows from \cite[Theorem 2.1]{Chinram2009}, but the rest of Proposition~\ref{prop:P1P2} appears to be new.
\end{rem}




Now that we have described the sets $P_1$, $P_2$, $P_3=P=P_1\cap P_2$, we may characterise Green's relations on~$\MmnJ$.  As in Section \ref{sect:partial}, if $\gK$ is one of $\R$, $\L$, $\H$, $\D$, $\J$, we will write $\gKJ$ for the Green's $\gK$-relation on $\MmnJ$.  Since $\MmnJ$ is not a monoid in general, these relations are defined, for $X,Y\in\Mmn$, by
\bit
\item $X\gRJ Y \ \iff \ [X=Y]$ or [$X=Y\star U$ and $Y=X\star V$ for some $U,V\in\Mmn$],
\eit
and so on.  
Since $\M$ is stable, so too is $\MmnJ$, so we have $\gJJ=\gDJ$ (see Proposition \ref{prop:stabilityMija} and Lemmas \ref{lem:Rol=LoR_S} and \ref{lem:green<M}).
We will continue to write $\R$, $\L$, $\H$, $\D$, $\J$ for the relations on $\M$ defined in Section \ref{sect:preliminaries}.  As in Section \ref{sect:partial}, if $\gK$ is one of $\R$, $\L$, $\H$, $\D=\J$, and if $X\in\Mmn$, we will write
\[
K_X=\set{Y\in\Mmn}{X\gK Y} \AND K_X^J=\set{Y\in\Mmn}{X\gKJ Y}
\]
for the $\gK$-class and $\gKJ$-class of~$X$ in $\Mmn$, respectively.  As noted in Section \ref{sect:partial}, $\gKJ\sub\gK$ for each $\gK$, and so $K_X^J\sub K_X$ for each $X$.
The next result follows immediately from Theorem \ref{thm:green_Sij}.  It also follows from Theorem 2.3, Lemma 2.4, and Corollaries 2.5--2.8 of \cite{Chinram2009}, but we prefer the current succinct description.

\ms
\begin{thm}\label{green_thm}
If $X\in\Mmn$, then   
\ms
\bmc2
\itemit{i} $R_X^J = \begin{cases}
R_X\cap P_1 &\text{if $X\in P_1$}\\
\{X\} &\text{if $X\in \Mmn\sm P_1$,}
\end{cases}$
\itemit{ii} $L_X^J = \begin{cases}
L_X\cap P_2 &\hspace{0.7mm}\text{if $X\in P_2$}\\
\{X\} &\hspace{0.7mm}\text{if $X\in \Mmn\sm P_2$,}
\end{cases}
\phantom{
\begin{cases}
a\\b\\c\\d
\end{cases}
}
$

\itemit{iii} $H_X^J = \begin{cases}
H_X &\hspace{6.8mm}\text{if $X\in P$}\\
\{X\} &\hspace{6.8mm}\text{if $X\in \Mmn\sm P$,}
\end{cases}$
\itemit{iv} $D_X^J = \begin{cases}
D_X\cap P &\text{if $X\in P$}\\
L_X^J &\text{if $X\in P_2\sm P_1$}\\
R_X^J &\text{if $X\in P_1\sm P_2$}\\
\{X\} &\text{if $X\in \Mmn\sm (P_1\cup P_2)$.}
\end{cases}$
\emc
The sets $P_1,P_2$ are described in Proposition \ref{prop:P1P2}, and the sets $R_X,L_X,H_X,D_X$ in Proposition \ref{lem:greenMmn}.  In particular, $R_X^J=L_X^J=H_X^J=D_X^J=\{X\}$ if $\rank(X)>r$.  If $X\in\Mmn\sm P$, then $H_X^J=\{X\}$ is a non-group $\gHJ$-class of $\MmnJ$.   \epfres
\end{thm}

Eggbox diagrams of some linear sandwich semigroups are given in Figures \ref{fig:V3212_V3322} and \ref{fig:V2322_V2422}.  As usual, grey boxes indicate group $\gHJ$-classes; a label of {\tt k} on such a group $\gHJ$-class indicates isomorphism to $\G_k$.  Note that the bottom diagram from Figure \ref{fig:V3212_V3322} is of a \emph{variant} of $\M_3(\Z_2)=\M_{33}(\Z_2)$.  The diagrams in the pdf version of this article may be zoomed in a long way.  The authors may be contacted for more such pictures.

\begin{figure}[ht]
\begin{center}
\includegraphics[width=\textwidth]{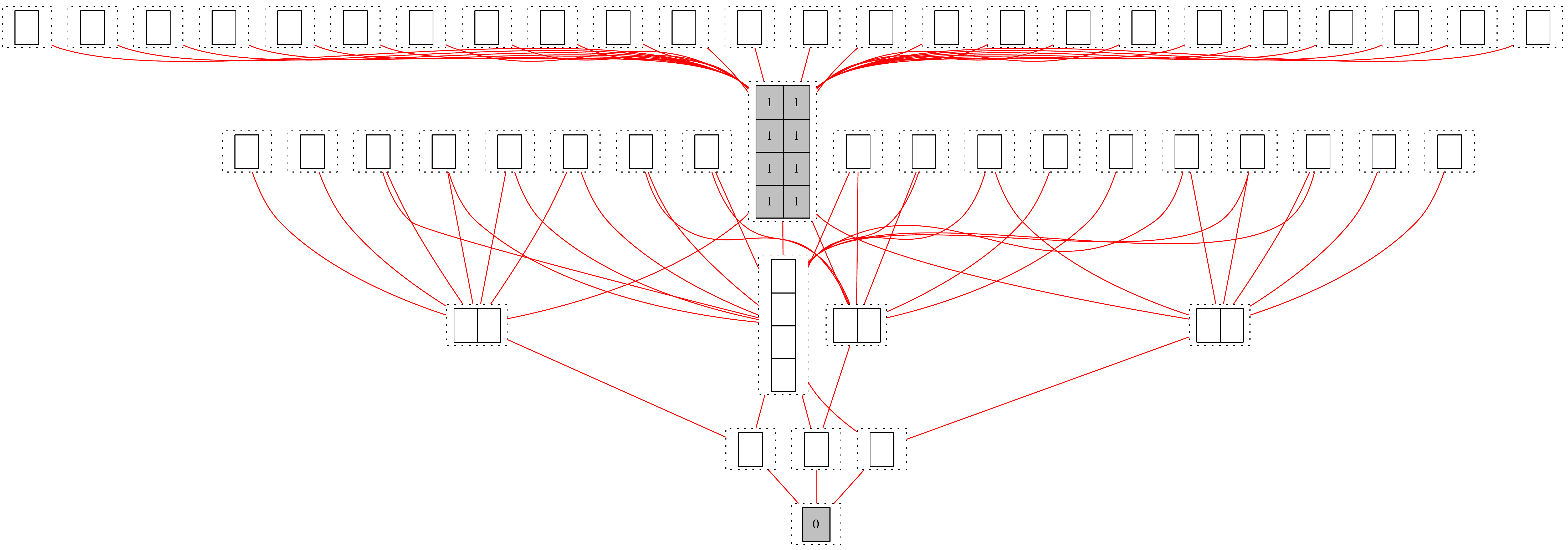}
\\~\\
\includegraphics[width=\textwidth]{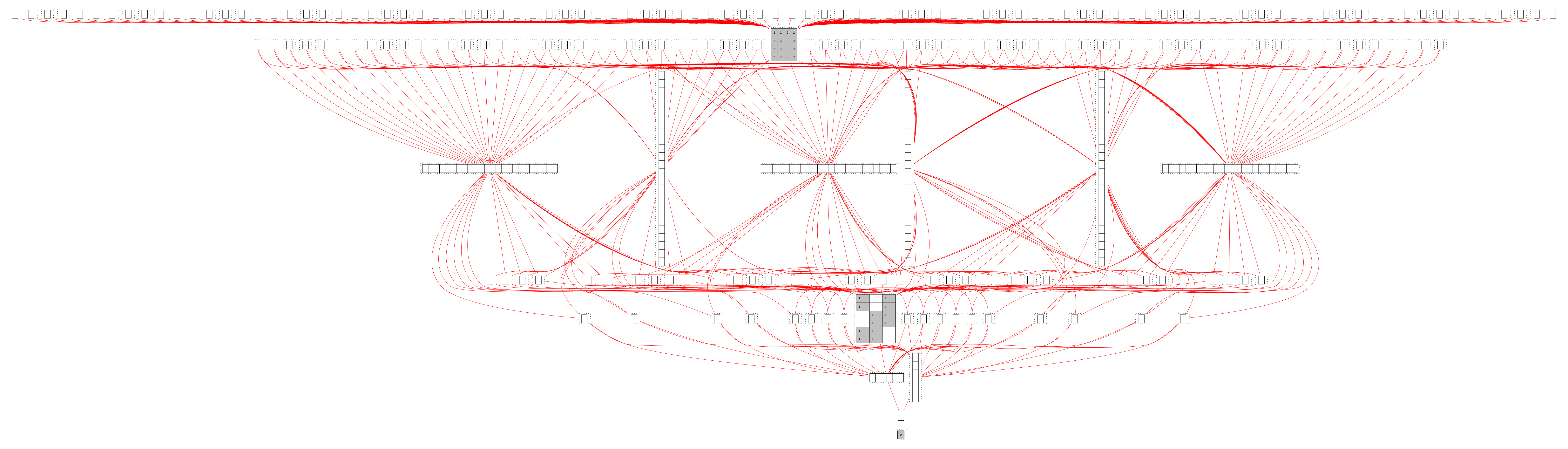}
    \caption{Egg box diagrams of the linear sandwich semigroups $\M_{32}^{J_{231}}(\Z_2)$ and $\M_{33}^{J_{332}}(\Z_2)$ (top and bottom, respectively).}
    \label{fig:V3212_V3322}
   \end{center}
 \end{figure}
 
 \begin{figure}[ht]
\begin{center}
\includegraphics[height=5cm]{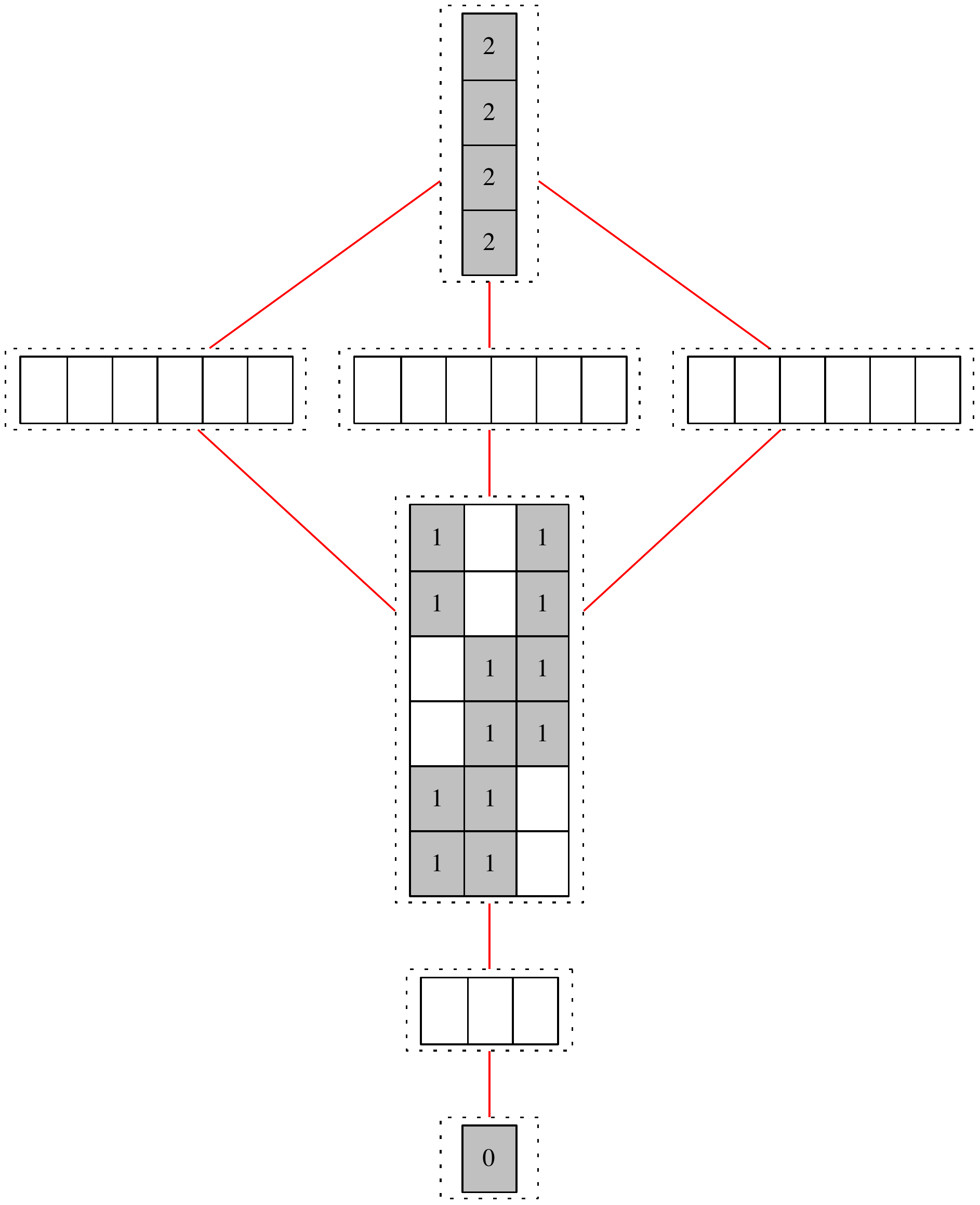} 
\qquad\qquad\qquad
\includegraphics[height=5cm]{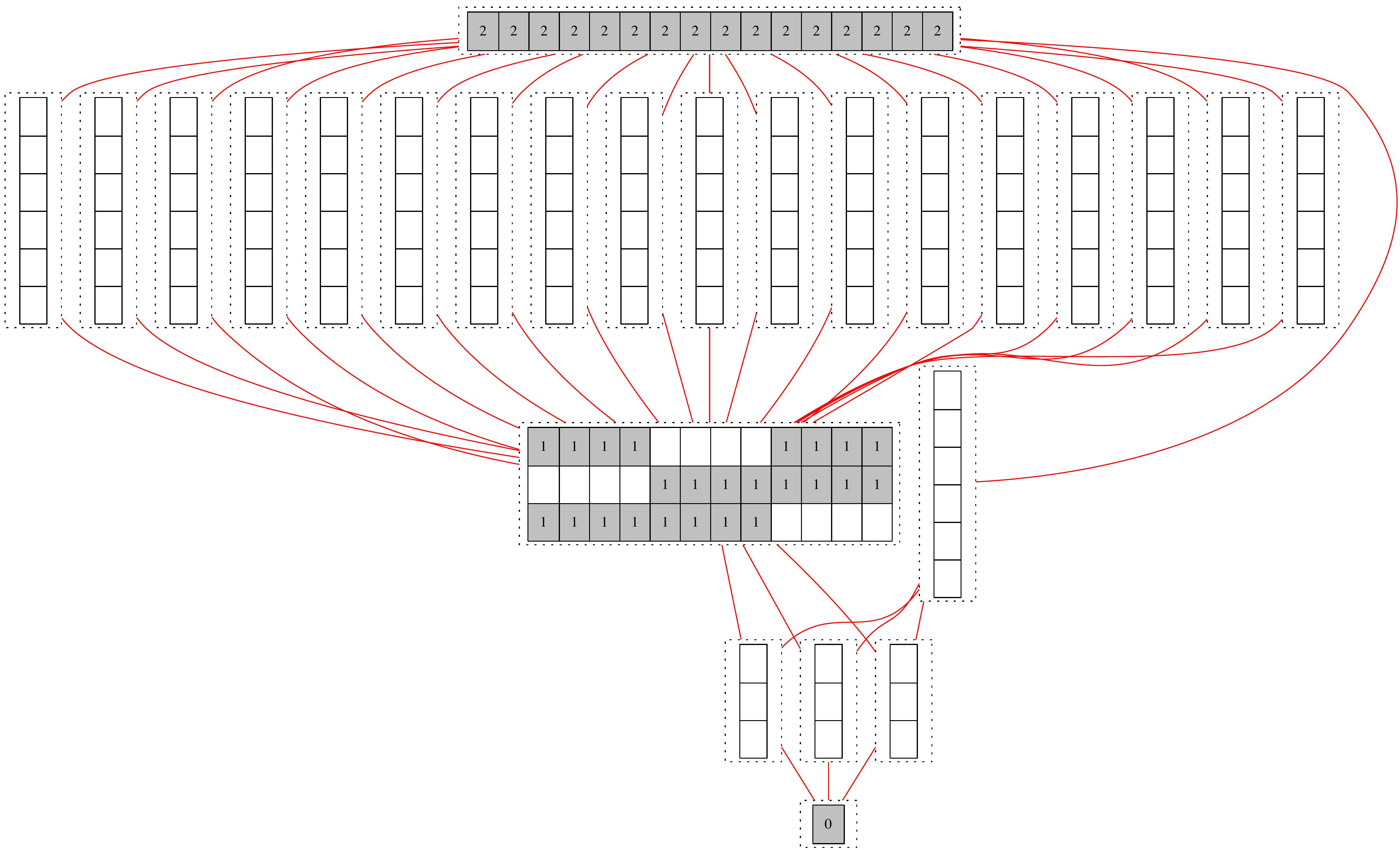} \\~\\
    \caption{Egg box diagrams of the linear sandwich semigroups $\M_{32}^{J_{232}}(\Z_2)$ and $\M_{24}^{J_{422}}(\Z_2)$ (left and right, respectively).}
    \label{fig:V2322_V2422}
   \end{center}
 \end{figure}


Theorem \ref{green_thm} yields an intuitive picture of the internal structure of $\MmnJ$.  Recall that the $\gD$-classes of $\Mmn$ are the sets $\DMmn s=\set{X\in\Mmn}{\rank(X)=s}$ for $0\leq s\leq l=\min(m,n)$.  If $r<l$, then each of the $\gD$-classes $\DMmn{r+1},\ldots,\DMmn l$ separates completely into singleton $\gDJ$-classes in $\MmnJ$.  (We will study these classes in more detail shortly.)  Next, note that $\DMmn0=\{O\}\sub P$ (as the zero matrix clearly belongs to both $P_1$ and $P_2$), so $\DMmn0$ remains a (regular) $\gDJ$-class of $\MmnJ$.
Now fix some $1\leq s\leq r$.  The $\gD$-class $\DMmn s$ is split into a single regular $\gDJ$-class, namely $\DMmn s\cap P$, and a number of non-regular $\gDJ$-classes.  Some of these non-regular $\gDJ$-classes are singletons, namely those of the form $D_X^J=\{X\}$ where $X\in \DMmn s$ belongs to neither $P_1$ nor $P_2$.  Some of the non-regular $\gDJ$-classes consist of one non-singleton $\gLJ$-class, namely those of the form $D_X^J=L_X^J=L_X\cap P_2$, where $X\in \DMmn s$ belongs to $P_2\sm P_1$; the $\gHJ$-classes contained in such a $\gDJ$-class are all singletons.  The remaining non-regular $\gDJ$-classes contained in $\DMmn s$ consist of one non-singleton $\gRJ$-class, namely those of the form $D_X^J=R_X^J=R_X\cap P_1$, where $X\in \DMmn s$ belongs to $P_1\sm P_2$; the $\gHJ$-classes contained in such a $\gDJ$-class are all singletons.  This is all pictured in Figure \ref{fig:green} for the $\D$-class $D_1(\M_{23})$ where $\F=\mathbb Z_3=\{0,1,2\}$ and $J=J_{321}=\left[\begin{smallmatrix}1&0\\0&0\\0&0\end{smallmatrix}\right]$; cf.~Figure \ref{fig:eggbox}. 

\nc{\fillboxbig}[2]{\draw[fill=gray!30](#1,#2)--(#1+1,#2)--(#1+1,#2+1)--(#1,#2+1)--(#1,#2);}

\begin{figure}[ht]
\begin{center}
\scalebox{.8}{
\begin{tikzpicture}[scale=0.9]
\draw[fill=lightgray!50,line width=.5mm] (0.75,2.5)--(14.25,2.5)--(14.25,8.5)--(0.75,8.5)--(0.75,2.5)-- (14.25,2.5); \bdmat{1.0}{7.9}100000
\dmat{1.0}{7.1}200000
\dmat{1.0}{5.9}100100
\dmat{1.0}{5.1}200200
\dmat{1.0}{3.9}100200
\dmat{1.0}{3.1}200100
\rldmat{1}{1.5-1}000100
\rrdmat{1}{1.5-1}000200
\dmat2{7.9}101000
\dmat2{7.1}202000
\dmat2{5.9}101101
\dmat2{5.1}202202
\dmat2{3.9}101202
\dmat2{3.1}202101
\rldmat{2}{1.5-1}000101
\rrdmat{2}{1.5-1}000202
\dmat3{7.9}102000
\dmat3{7.1}201000
\dmat3{5.9}102102
\dmat3{5.1}201201
\dmat3{3.9}102201
\dmat3{3.1}201102
\rldmat{3}{1.5-1}000102
\rrdmat{3}{1.5-1}000201
\dmat4{7.9}110000
\dmat4{7.1}220000
\dmat4{5.9}110110
\dmat4{5.1}220220
\dmat4{3.9}110220
\dmat4{3.1}220110
\rldmat{4}{1.5-1}000110
\rrdmat{4}{1.5-1}000220
\dmat5{7.9}111000
\dmat5{7.1}222000
\dmat5{5.9}111111
\dmat5{5.1}222222
\dmat5{3.9}111222
\dmat5{3.1}222111
\rldmat{5}{1.5-1}000111
\rrdmat{5}{1.5-1}000222
\dmat6{7.9}112000
\dmat6{7.1}221000
\dmat6{5.9}112112
\dmat6{5.1}221221
\dmat6{3.9}112221
\dmat6{3.1}221112
\rldmat{6}{1.5-1}000112
\rrdmat{6}{1.5-1}000221
\dmat7{7.9}120000
\dmat7{7.1}210000
\dmat7{5.9}120120
\dmat7{5.1}210210
\dmat7{3.9}120210
\dmat7{3.1}210120
\rldmat{7}{1.5-1}000120
\rrdmat{7}{1.5-1}000210
\dmat8{7.9}121000
\dmat8{7.1}212000
\dmat8{5.9}121121
\dmat8{5.1}212212
\dmat8{3.9}121212
\dmat8{3.1}212121
\rldmat{8}{1.5-1}000121
\rrdmat{8}{1.5-1}000212
\dmat{9}{7.9}122000
\dmat{9}{7.1}211000
\dmat{9}{5.9}122122
\dmat{9}{5.1}211211
\dmat{9}{3.9}122211
\dmat{9}{3.1}211122
\rldmat{9}{1.5-1}000122
\rrdmat{9}{1.5-1}000211
\dmat{10.666}{8.0}010000
\dmat{10.666}{7.0}020000
\dmat{10.666}{6.0}010010
\dmat{10.666}{5.0}020020
\dmat{10.666}{4.0}010020
\dmat{10.666}{3.0}020010
\dmat{10.666}{2.0-1}000010
\dmat{10.666}{1.0-1}000020
\dmat{11.666}{8.0}011000
\dmat{11.666}{7.0}022000
\dmat{11.666}{6.0}011011
\dmat{11.666}{5.0}022022
\dmat{11.666}{4.0}011022
\dmat{11.666}{3.0}022011
\dmat{11.666}{2.0-1}000011
\dmat{11.666}{1.0-1}000022
\dmat{12.666}{8.0}012000
\dmat{12.666}{7.0}021000
\dmat{12.666}{6.0}012012
\dmat{12.666}{5.0}021021
\dmat{12.666}{4.0}012021
\dmat{12.666}{3.0}021012
\dmat{12.666}{2.0-1}000012
\dmat{12.666}{1.0-1}000021
\dmat{13.666}{8.0}001000
\dmat{13.666}{7.0}002000
\dmat{13.666}{6.0}001001
\dmat{13.666}{5.0}002002
\dmat{13.666}{4.0}001002
\dmat{13.666}{3.0}002001
\dmat{13.666}{2.0-1}000001
\dmat{13.666}{1.0-1}000002
%
\foreach \x in {4,6} \draw (0.75,\x+.5)--(14.25,\x+.5); 
\foreach \x in {3,...,7} \foreach \y in {14.5,16,17.5,19} \draw (0.75+.15+\y,\x+.5)--(2.25-.15+\y,\x+.5);  
\foreach \x in {1,...,8} \draw (\x*1.5+0.75,2.5)--(\x*1.5+0.75,8.5);  
\foreach \x in {0.5,1,1.5,2,...,8,8.5} \draw (\x*1.5+0.75,2+.35-1)--(\x*1.5+0.75,1-.35-1); 
\foreach \x in {1,2.5,4,5.5} \foreach \y in {0,-1} {\draw[line width=.5mm] (15-.6+\x,1-.35+\y)--(15+.6+\x,1-.35+\y)-- (15+.6+\x,1+.35+\y)--(15-.6+\x,1+.35+\y)--(15-.6+\x,1-.35+\y)--(15+.6+\x,1-.35+\y);} 
\foreach \x in {1,2.5,4,5.5} {\draw[line width=.5mm] (15-.6+\x,1-.5+2)--(15+.6+\x,1-.5+2)-- (15+.6+\x,8.5)--(15-.6+\x,8.5)--(15-.6+\x,1-.5+2)--(15+.6+\x,1-.5+2);} 
\draw[line width=.5mm] (0.75,2.5-.15-1)--(14.25,2.5-.15-1)--(14.25,0.5+.15-1)--(0.75,0.5+.15-1) --(0.75,2.5-.15-1)--(14.25,2.5-.15-1); 
\draw[|-|] (0.75,9.0)--(14.25,9.0);
\draw[|-|] (20.25+1,9.0)--(14.25+1,9.0);
\draw[|-|] (.25,.5-1)--(.25,2.5-1);
\draw[|-|] (.25,8.5)--(.25,2.5);
\draw(7.5,9.3)node{$\sub P_1\phantom{\sub}$};
\draw(16.75+1.5,9.3)node{$\not\sub P_1\phantom{\sub}$};
\draw(.25,5.5)node[left]{$\sub P_2\phantom{}$};
\draw(.25,1.5-1.0)node[left]{$\not\sub P_2\phantom{}$};
\end{tikzpicture}
}
    \caption{A $\gD$-class $D_1(\M_{23}(\mathbb Z_3))$ breaks up into $\gDJ$-classes in $\M_{23}^{J}(\mathbb Z_3)$, where $J=J_{321}$.  Group $\gH^J$-classes are shaded grey; the idempotent of such a group is the upper of the two matrices.  (cf.~Figure \ref{fig:eggbox}.)} 
    \label{fig:green}
   \end{center}
 \end{figure}

It will be important to have a description of the partial order $\leq$ on the $\gDJ$-classes of $\MmnJ$.


\ms
\begin{prop}\label{prop:DorderMmnJ}
Let $X,Y\in\Mmn$.  Then $D_X^J\leq D_Y^J$ in $\MmnJ$ if and only if one of the following holds:
\ms
\bit
\begin{multicols}{2}
\itemit{i} $X=Y$,
\itemit{ii} $\rank(X)\leq\rank(JYJ)$,
\itemit{iii} $\Row(X)\sub\Row(JY)$,
\itemit{iv} $\Col(X)\sub\Col(YJ)$.
\end{multicols}
\end{itemize}
\end{prop}


\pf Note that $D_X^J\leq D_Y^J$ if and only if one of the following holds:
\ms
\bmc2
\item[(a)] $X=Y$,
\item[(b)] $X=UJYJV$ for some $U,V\in\Mmn$,
\item[(c)] $X=UJY$ for some $U\in \Mmn$,
\item[(d)] $X=YJV$ for some $V\in\Mmn$.
\emc
The equivalences (b) $\Leftrightarrow$ (ii), (c) $\Leftrightarrow$ (iii), and (d) $\Leftrightarrow$ (iv) all follow from Lemma \ref{lem:green<M}.  \epf

The description of the order on $\gDJ$-classes of $\MmnJ$ from Proposition \ref{prop:DorderMmnJ} may be simplified in the case that one of $X,Y$ is regular.

\ms
\begin{prop}\label{prop:DorderP}
Let $X,Y\in\Mmn$.
\bit
\itemit{i} If $X\in P$, then $D_X^J\leq D_Y^J\iff\rank(X)\leq\rank(JYJ)$.
\itemit{ii} If $Y\in P$, then $D_X^J\leq D_Y^J\iff\rank(X)\leq\rank(Y)$.
\eit
The regular $\gDJ$-classes of $\MmnJ$ form a chain: $D_0^J<\cdots<D_r^J$, where $$D_s^J=\DMmn s\cap P=\set{X\in P}{\rank(X)=s} \qquad\text{for each $0\leq s\leq r$.}$$
\end{prop}

\pf As in the proof of Proposition \ref{prop:DorderMmnJ}, $D_X^J\leq D_Y^J$ if and only if one of (a--d) holds.
Suppose first that $X\in P$, so $X=XJZJX$ for some $Z\in\Mmn$.  Then (a) implies $X=XJZ(JYJ)ZJX$, (c) implies $X=U(JYJ)ZJX$, and (d) implies $X=XJZ(JYJ)V$.  So, in each of cases (a--d), we deduce that $\rank(X)\leq\rank(JYJ)$.  So $D_X^J\leq D_Y^J$ implies $\rank(X)\leq\rank(YJY)$.  Proposition \ref{prop:DorderMmnJ} gives the reverse implication.  

Next, suppose $Y\in P$.  Now, each of (a--d) implies $\rank(X)\leq\rank(Y)$.  Conversely, if $\rank(X)\leq\rank(Y)$, then Proposition~\ref{prop:DorderMmnJ} gives $D_X^J\leq D_Y^J$, since $\rank(Y)=\rank(JYJ)$.
The statement about regular $\gDJ$-classes follows quickly from (ii).~\epf

The linear ordering on the regular $\gDJ$-classes may be seen by inspecting Figures \ref{fig:V3212_V3322} and \ref{fig:V2322_V2422}; see also Figure~\ref{fig:R}.
As an immediate consequence of Proposition \ref{prop:DorderP}, we may classify the isomorphism classes of sandwich semigroups on the set~$\Mmn$;  the $m=n$ case of the next result was proved in \cite{JCK2010}.

\ms
\begin{cor}\label{cor:MmnAcongMmnB}
Let $A,B\in\Mnm$.  Then $\MmnA\cong\MmnB$ if and only if $\rank(A)=\rank(B)$.
\end{cor}

\pf Put $r=\rank(A)$ and $s=\rank(B)$.  By Proposition \ref{prop:DorderP} and Lemma \ref{lem:MmnAMmnB}(ii), $\MmnA\cong\Mmn^{J_{nmr}}$ and $\MmnB\cong\Mmn^{J_{nms}}$ have $r+1$ and $s+1$ regular $\D^A$- and $\D^B$-classes, respectively.  So $\MmnA\cong\MmnB$ implies $r=s$.  The converse was proved in Lemma~\ref{lem:MmnAMmnB}(ii).~\epf


\begin{rem}
It is possible to have $\M_{mn}^A\cong\M_{kl}^B$ even if $(m,n)\not=(k,l)$, although we would of course still need $\rank(A)=\rank(B)$ by Proposition \ref{prop:DorderP}.  For example, if $O=O_{nm}$ is the $n\times m$ zero matrix, then $\M_{mn}^O$ is a \emph{zero semigroup} ($X\star Y=O_{mn}$ for all $X,Y\in\M_{mn}$).  Two such zero semigroups $\Mmn^O$ and $\M_{kl}^O$ are isomorphic if and only if they have the same cardinality; that is, if and only if $\F$ is infinite or $\F$ is finite and $mn=kl$.  We will return to the problem of distinguishing non-isomorphic $\M_{mn}^A$ and $\M_{kl}^B$ in Theorem \ref{thm:classification}.  See Figure \ref{fig:V2202}.
\end{rem}

\begin{figure}[ht]
\begin{center}
\includegraphics[width=12cm]{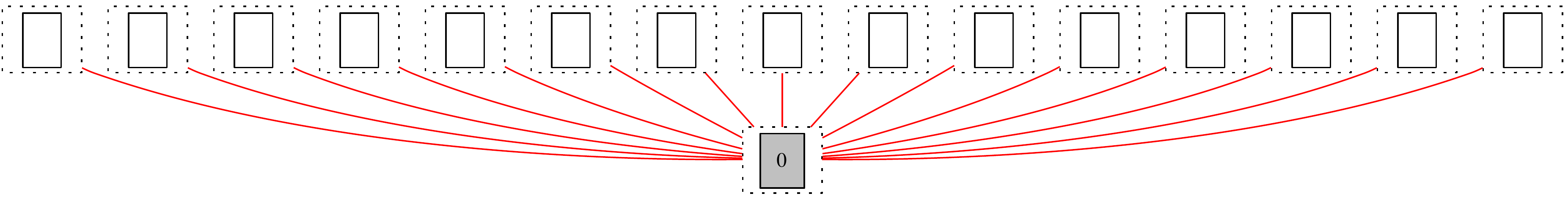}
    \caption{Egg box diagram of the linear sandwich semigroup $\M_{22}^{O_{22}}(\Z_2)$ or, equivalently, $\M_{21}^{O_{12}}(\F_4)$.  Both are zero semigroups of size $16$.}
    \label{fig:V2202}
   \end{center}
 \end{figure}

The next result describes the maximal $\gDJ$-classes of $\MmnJ$.  See also Figures \ref{fig:V3212_V3322} and \ref{fig:V2322_V2422}.

\ms\ms\ms
\begin{prop}\label{prop_maximalD}
\bit
\itemit{i} If $r=\min(m,n)$, then $D_r^J=D_r\cap P=\set{X\in P}{\rank(X)=r}$ is the unique maximal $\gDJ$-class of~$\MmnJ$, and is a subsemigroup of $\MmnJ$.
\itemit{ii} If $r<\min(m,n)$, then the maximal $\gDJ$-classes of $\MmnJ$ are those of the form $D_X^J=\{X\}$ with $\rank(X)>r$.
\eitres
\end{prop}

\pf Part (i) follows immediately from Proposition \ref{prop:DorderP}(ii), the rule $[M,A,N]\star[K,B,L]=[M,AB,L]$, and the fact that $\G_r=D_r(\M_r)$ is a subgroup of $\M_r$.  

For (ii), let $X\in\Mmn$.  Suppose first that $\rank(X)>r$ and that $D_X^J\leq D_Y^J$.  Then condition (ii) from Proposition \ref{prop:DorderMmnJ} does not hold, since $\rank(JYJ)\leq\rank(J)=r<\rank(X)$.  Similarly, $\rank(JY)<\rank(X)$ and $\rank(YJ)<\rank(X)$, so neither (iii) nor (iv) holds.  Having eliminated (ii--iv), we deduce that (i) must hold; that is, $X=Y$, so $D_X^J=\{X\}$ is indeed maximal.  
Conversely, suppose
$\rank(X)\leq r$, and let $Y=\tmat {I_r}OOD$, where $D\not=O$.  Then $\rank(Y)>r$, so $D_Y^J=\{Y\}$ is maximal by the previous paragraph.  But also $JYJ=J$, and it follows that $\rank(X)\leq r=\rank(J)=\rank(JYJ)$, so that $D_X^J< D_Y^J=\{Y\}$, whence $D_X^J$ is not maximal.  \epf

The description of the maximal $\gDJ$-classes from Proposition \ref{prop_maximalD} allows us to obtain information about generating sets for $\MmnJ$ and, in the case of finite $\F$, about $\rank(\MmnJ)$.  
In order to avoid confusion when discussing generation, if $\Om\sub\Mmn$, we will write $\la \Om\raJ$ for the subsemigroup of $\MmnJ$ generated by $\Om$, which consists of all products $X_1\star\cdots\star X_k$, with $k\geq1$ and $X_1,\ldots,X_k\in \Om$.  If $\Si\sub\M_k$ for some $k$, we will continue to write $\la\Si\ra$ for the subsemigroup of $\M_k$ generated by $\Si$.  For convenience, we will state two separate results, according to whether $r=\min(m,n)$ or $r<\min(m,n)$.  The next lemma will be useful as the inductive step in the proofs of both Theorems \ref{thm:rankMmnJ} and \ref{thm:rankMmnJ_r=m}.  Recall that $\{e_{m1},\ldots,e_{mm}\}$ is the standard basis of $V_m=\F^m$.

\ms
\begin{lemma}\label{lem:ind_step}
Suppose $X\in \DMmn s$, where $0\leq s\leq l-1$ and $l=\min(m,n)$.  Then $X=Y\star Z$ for some $Y\in \DMmn l$ and $Z\in\DMmn{s+1}$.
\end{lemma}

\pf 
Let $\B=\{v_1,\ldots,v_n\}$ be a basis of $V_n$ such that $\{v_{s+1},\ldots,v_n\}$ is a basis of $\ker(\lam_X)$.  
%
%
Consider the linear transformation $\be\in\Homnm$ defined by
\[
\be(v_i) = \begin{cases}
e_{mi} &\text{if $1\leq i\leq s$}\\
0 &\text{if $s<i<n$}\\
e_{mm} &\text{if $i=n$,}
\end{cases}
\]
noting that $\rank(\be)=s+1$.  The proof now breaks into two cases, depending on whether $r<m$ or $r=m$.

{\bf Case 1.}  Suppose first that $r<m$.  Let $\al\in\Homnm$ be any linear transformation of rank $l$ that extends the map $e_{ni}\mt \lam_X(v_i)$ ($1\leq i\leq s$).  One easily checks that $\al\circ\lam_J\circ\be=\lam_X$.

{\bf Case 2.}  Now suppose $r=m$.  Recall that we are assuming that $r=m=n$ does not hold,
so $r=m<n$.  This time, define we let $\al$ be any linear transformation of rank $m=l$ that extends the map $e_{ni}\mt\lam_X(v_i)$ ($1\leq i\leq s$), $e_{nr}=e_{nm}\mt0$.  Then, again, one easily checks that $\al\circ\lam_J\circ\be=\lam_X$.  \epf

\ms
\begin{thm}\label{thm:rankMmnJ}
Suppose $r<l=\min(m,n)$.  Then $\MmnJ=\la \Om\raa$, where $\Om=\set{X\in\Mmn}{\rank(X)>r}$.  Further, any generating set for $\MmnJ$ contains $\Om$. 
If $|\F|=q<\infty$, then
\[
\rank(\MmnJ)=|\Om|=\sum_{s=r+1}^{l} \qbin ms \qbin ns q^{{s\choose2}}(q-1)^s\qfact s.
\]
\end{thm}

\pf For convenience, we will assume that $l=m\leq n$.  The other case will follow by duality.  We will also denote $D_s(\Mmn)$ simply by $D_s$ for each $0\leq s\leq m$.  Consider the statement:
\nss
\begin{quote}
$H(s)$: ~ $\la \Om\raa$ contains $D_s\cup\cdots\cup D_m=\set{X\in\Mmn}{\rank(X)\geq s}$.
\end{quote}
\nss
Note that $\Om=D_{r+1}\cup\cdots\cup D_m$, so $H(s)$ is clearly true for $r+1\leq s\leq m$.  Lemma \ref{lem:ind_step} shows that $H(s+1)$ implies $H(s)$ for all $0\leq s\leq m-1$.  So we conclude that $H(s)$ is true for all $0\leq s\leq m$.
In particular, $H(0)$ says that $\MmnJ=\la \Om\raa$.

Since $\{X\}$ is a maximal $\gDJ$-class for any $X\in \Om$, it follows that any generating set of $\MmnJ$ must contain~$\Om$.  Thus, $\Om$ is the minimal generating set with respect to both size and containment, so $\rank(\MmnJ)=|\Om|$.  The formula for $|\Om|$ with $|\F|$ finite follows from Lemma \ref{lem:combinatorics_Mmn}. \epf

In order to consider the case in which $r=\min(m,n)$, we first prove an intermediate result.  There is a dual version of the following lemma (dealing with the case in which $r=n<m$), but we will not state it.

\ms
\begin{lemma}\label{lem:r=m}
If $r=m<n$, then 

~\ \ \emph{(i)} $P_2=\MmnJ$, \qquad
\emph{(ii)} $P=P_1$ is a left ideal of $\MmnJ$, \qquad
\emph{(iii)} $\gLJ=\L$ in $\MmnJ$.
%
\end{lemma}

\pf Let $X\in\MmnJ$.  As noted earlier, in the $2\times2$ block description, $X=\tmat ABCD$ (where $A\in\M_r$, and so on), the matrices $C$ and $D$ are empty (since $r=m$).  So we write $X=[A\ B]$.  Note that $J=\thmat IO$, so $JX=\thmat IO[A\ B]=\tmat ABOO$.  It follows that $\Row(JX)=\Row(X)$ and, since $X\in\MmnJ$ was arbitrary, this completes the proof of (i).

We immediately deduce $P=P_1$ from (i).  As in Proposition~\ref{prop:P1P2}, the regular elements of $\MmnJ$ are of the form $[A\ AN]$ where $A\in\M_r$ and $N\in\M_{r,n-r}$.  We denote such a regular element by $[A,N]$.  The proof of~(ii) concludes with the easily checked observation that $[A\ B]\star[C,N]=[AC,N]$.

Part (iii) follows quickly from (i) and Theorem \ref{green_thm}(ii). \epf


\ms
\begin{thm}\label{thm:rankMmnJ_r=m}
Suppose $r=\min(m,n)$ where $m\not=n$.  If $|\F|=q<\infty$, then 
\[
\rank(\MmnJ)
= \qbin Ll,
\]
where $L=\max(m,n)$ and $l=\min(m,n)$.
\end{thm}

\pf Again, it suffices to assume that $r=m<n$, so $l=m$ and $L=n$.  We keep the notation of the previous proof.

Let $\Om$ be an arbitrary generating set for $\MmnJ$.  Let $X\in D_m(\Mmn)$ be arbitrary.  We claim that $\Om$ must contain some element of $L_X^J=L_X$.  Indeed, consider an expression $X=Y_1\star\cdots\star Y_k$, where $Y_1,\ldots,Y_k\in\Om$.    If $k=1$, then $Y_1=X\in L_X$ and the claim is established, so suppose $k\geq2$.  
Since $D_m(\Mmn)$ is a maximal $\gDJ$-class, we must have $Y_k\in D_m(\Mmn)$.  So $Y_k \gDJ X = (Y_1\star\cdots\star Y_{k-1})\star Y_k$, whence $Y_k\gLJ (Y_1\star\cdots\star Y_{k-1})\star Y_k = X$, by stability.  By Lemma \ref{lem:r=m}(iii), this completes the proof of the claim.
In particular, $|\Om|$ is bounded below by the number of $\gL$-classes contained in $D_m(\Mmn)$, which is equal to $\tqbin nm$, by Lemma~\ref{lem:combinatorics_Mmn}.  Since $\Om$ was an arbitrary generating set, it follows that $\rank(\MmnJ)\geq\tqbin nm=\tqbin Ll$.

To complete the proof, it remains to check that there exists a generating set of the desired cardinality.  
For each $N\in\M_{m,n-m}$, choose some $A_N\in\G_r$ such that $\set{A_N}{N\in\M_{m,n-m}}$ generates $\G_m$, and put $X_N=[A_N,N]\in D_m^J$.  (This is possible since $|\M_{m,n-m}|=q^{m(n-m)}\geq2$, and $\rank(\G_m)\leq2$ by Theorem~\ref{thm:waterhouse}.)  It is easy to see that $\Om_1=\set{X_N}{N\in\M_{m,n-m}}$ is a cross-section of the $\L$-classes in $D_m^J$.  Also, choose some cross-section $\Om_2=\set{Y_i}{i\in I}$ of the $\L$-classes contained in $D_m(\Mmn)\sm D_m^J$.  Then
$
\Om = \Om_1\cup\Om_2 
$
is a cross-section of the $\L$-classes contained in $D_m(\Mmn)$.  Since, therefore, $|\Om|=\tqbin nm$, the proof will be complete if we can show that $\MmnJ=\la\Om\raJ$.  By Lemma \ref{lem:ind_step}, it suffices to show that $\la\Om\raJ$ contains $D_m(\Mmn)$.  So suppose $Z\in D_m(\Mmn)$.  Assume first that $Z\in D_m^J$, and write  $Z=[B,L]$, noting that $B\in\G_r$.  Choose $N_1,\ldots,N_k\in\M_{m,n-m}$ such that $BA_L^{-1}=A_{N_1}\cdots A_{N_k}$.  Then one easily checks that $Z=X_{N_1}\star\cdots\star X_{N_k}\star X_L$.  Now, suppose $Z$ is not regular.  Choose $i\in I$ such that $Z\L Y_i$.  By Lemma \ref{lem:greenMmn}, $Z=UY_i$ for some $U\in\G_m$.  But then $Z=[U\ V]\star Y_i$ for any $V\in\M_{m,n-m}$.  Since $\rank(U)=m$, we have $[U\ V]\in D_m^J\sub\la\Om\raJ$, whence $Z\in\la\Om\raJ$, completing the proof. \epf


\begin{rem}
By inspecting Figures \ref{fig:V3212_V3322} and \ref{fig:V2322_V2422}, the reader may use Theorems \ref{thm:rankMmnJ} and \ref{thm:rankMmnJ_r=m} to locate the elements from a minimal generating set for $\MmnJ$.
\end{rem}

\section{Connection to (non-sandwich) matrix semigroups}\label{sect:non-sandwich}

Recall that $J=J_{nmr}=\tmat {I_r}{O_{r,n-r}}{O_{m-r,r}}{O_{m-r,n-r}}\in\M_{nm}$.  Now let $K=J^T=J_{mnr}=\tmat {I_r}{O_{r,m-r}}{O_{n-r,r}}{O_{n-r,m-r}}\in\M_{mn}$.  So Lemma \ref{lem:MmnAMmnB} says that $\MnmK$ and $\MmnJ$ are anti-isomorphic.  Also, since $J=JKJ$ and $K=KJK$, Theorem \ref{thm:diamonds} says that we have the following commutative diagrams of semigroup homomorphisms where, for clarity, we write $\cdot$ for (non-sandwich) matrix multiplication:
\[
\includegraphics{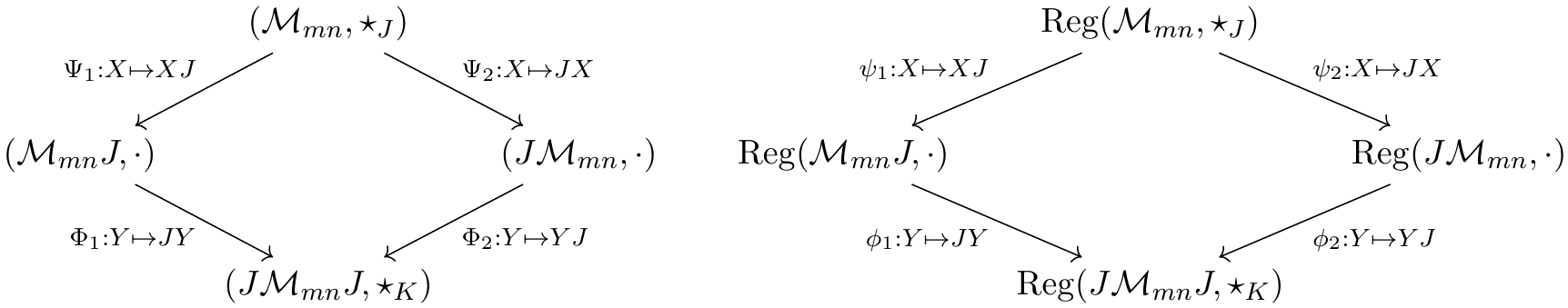}
\]
In this section, we show that the various semigroups appearing in the above diagrams are all (equal to or isomorphic to) certain well-known (non-sandwich) matrix semigroups, and explore the consequences for the structure of the sandwich semigroups $\MmnJ$.  First, we have a simple observation.

\ms
\begin{lemma}\label{lem:JMJ}
We have $(J\Mmn J,\star_K)=\Reg(J\Mmn J,\star_K)\cong(\M_r,\cdot)$.
\end{lemma}

\pf Let $X=\tmat ABCD\in\Mmn$.  We have already observed that, whether $X$ is regular or not, $JXJ=\tmat AOOO\in\Mnm$.  The result follows quickly from the fact that $\tmat AOOO\star_K \tmat EOOO=\tmat{AE}OOO$. \epf

For integers $k\geq1$ and $0\leq l\leq k$, we write 
\begin{align*}
\cC_k(l) &= \set{X\in\M_k}{\col_{l+1}(X)=\cdots=\col_k(X)=O}, \\
\cR_k(l) &= \set{X\in\M_k}{\row_{l+1}(X)=\cdots=\row_k(X)=O}.
\intertext{(As before, without causing confusion, we write $O$ for any zero matrix when the dimensions are clear from context.)  These matrix semigroups have been studied in a number of contexts (see for example \cite{NK2007,Sullivan2008}), along with their associated isomorphic semigroups of linear transformations}
\mathcal K_k(l) &= \bigset{\al\in\End(V_k)}{\ker(\al)\supseteq W_{kl}^\perp},\\
\mathcal I_k(l) &= \bigset{\al\in\End(V_k)}{\im(\al)\sub W_{kl}}.
\end{align*}
Here we have written $W_{kl}^\perp=\Span\{e_{k,l+1},\ldots,e_{kk}\}$.  Clearly, $\Ckl$ and $\Rkl$ are anti-isomorphic. 

\ms
\begin{lemma}
We have $\MMNJ=\Cmr$ and $\JMMN=\Rnr$.
\end{lemma}

\pf Let $X=\tmat ABCD\in\Mmn$.  We have already observed that $XJ=\tmat AOCO\in\M_m$ and $JX=\tmat ABOO$, and the result quickly follows. \epf

\begin{rem}\label{rem:r=m2}
A typical element $X\in\Rkl$ may be written as $X=\tmat ABOO$, where $A\in\M_l$, $B\in\M_{l,k-l}$ and so on.  One easily checks that multiplication of matrices in this form obeys the rule $\tmat ABOO \tmat EFOO = \tmat {AE}{AF}OO$.  Comparing this to the discussion in Remark \ref{rem:r=m}, we see that $\Rkl$ is isomorphic to the sandwich semigroup $\M_{lk}^J$ where $J=J_{kll}\in\M_{kl}$.  (A dual statement holds for the matrix semigroups $\Ckl$.)  Thus, every result we obtain for linear sandwich semigroups leads to analogous results for the semigroups $\Rkl$ and $\Ckl$.  For example, we deduce from Theorem \ref{thm:rankMmnJ_r=m} that 
$
\rank(\Ckl) = \rank(\Rkl) = \tqbin kl
$
if $|\F|=q<\infty$.
Note that the sandwich semigroups $\MmnJ$ pictured in Figure \ref{fig:V2322_V2422} satisfy $r=\min(m,n)$, so Figure \ref{fig:V2322_V2422} essentially pictures eggbox diagrams of $\mathcal C_3(2)$ and $\mathcal R_4(2)$.
\end{rem}


\ms
\begin{rem}
Similarly, one may think of an arbitrary linear sandwich semigroup $\MmnJ$ itself as a (non-sandwich) matrix semigroup, as noted by Thrall in \cite{Thrall1955} and slightly adapted as follows.  Consider the set $\mathscr M$ of all $(m+n-r)\times (m+n-r)$ matrices
that may be written in $3\times3$ block form
$\tmatt OOOBAODCO$,
where $A\in\M_r$, $D\in\M_{m-r,n-r}$ (and from which the dimensions of the other sub-matrices may be derived).  One easily checks that the matrices from $\mathscr M$ multiply according to the rule
$\tmatt OOOBAODCO \tmatt OOOFEOHGO = \tmatt OOO{AF}{AE}O{CF}{CE}O$,
so that 
$\tmat ABCD \mt \tmatt OOOBAODCO$
determines an isomorphism $(\Mmn,\star_J)\to(\mathscr M,\cdot)$.  Note also that
\begin{align*}
\mathscr M &= \cR_{m+n-r}^*(m) \cap \cC_{m+n-r}(n)
\end{align*}
where here we write
$\cR_k^*(l) = \set{X\in\M_k}{\row_1(X)=\cdots=\row_{k-l}(X)=O}$.
(It is easily seen that the map $\tmat ABOO\to\tmat OOBA$ determines an isomorphism $\Rkl\to \cR_k^*(l)$.)
Since using this isomorphic copy $\mathscr M$ of $\MmnJ$ does not appear to confer any obvious advantage, we will make no further reference to it.
\end{rem}

The regular elements of $\Ckl$ and $\Rkl$ 
(and also of $\mathcal I_k(l)$ and $\mathcal K_k(l)$) 
were classified in \cite{NK2007}.  The next result, which gives a much simpler description of these regular elements, may be deduced from \cite[Theorems~3.4 and 3.8]{NK2007} (and vice versa), but we include a simple proof for convenience.  

\ms
\begin{prop}
\label{cor:RegMNJ}
The regular elements of the semigroups $\Cmr=\MMNJ$ and $\Rnr=\JMMN$ are given by
\begin{align*}
\Reg(\Cmr) = \Reg(\MMNJ) = PJ &= \set{X\in\MMNJ}{\rank(JX)=\rank(X)} \\
\Reg(\Rnr) = \Reg(\JMMN) = JP &= \set{X\in\JMMN}{\rank(XJ)=\rank(X)}.
\end{align*}
\end{prop}

\pf We just prove the second statement as the other is dual.  Let $X=\tmat ABCD\in\Mmn$, and put $X'=\tmat ABOO\in\Mmn$.  Then $JX=JX'=\tmat ABOO\in\M_n$ (where the zero matrices in the last expression have $n-r$ rows).  Since $X'$ clearly belongs to $P_2$ (by Proposition \ref{prop:P1P2}), we have $\JMMN\sub JP_2$.  Next, note that $KJ=J_{mmr}$, so that $KJY=Y$ for all $Y\in\M_{mn}$ of the form $Y=\tmat ABOO$.  Now suppose $X\in\Mmn$ is such that $JX\in\RegJMMN$.  As above, we may assume that $X=\tmat ABOO$.  So $(JX)=(JX)(JY)(JX)$ for some $Y\in\Mmn$.  But then
$
X = K(JX) = K(JXJYJX)=XJYJX = X\star Y\star X,
$
so that, in fact, $X\in \RegMmnJ=P$.  This completes the proof that $\RegJMMN\sub JP$.  The reverse inclusion is easily checked.

Now suppose $X=JY$ where $Y=\tmat A{AN}OO\in P$.  Then $X=JY=\tmat A{AN}OO$ (with appropriately sized zero matrices), so $\rank(X)=\rank(JY)=\rank(Y)=\rank(JYJ)=\rank(XJ)$, where we have used Proposition~\ref{prop:P1P2}.  Conversely, suppose $X\in\JMMN$ is such that $\rank(XJ)=\rank(X)$.  As before, we may assume that $X=JY$ where $Y\in P_2$.  Then $\rank(Y)=\rank(JY)=\rank(X)=\rank(XJ)=\rank(JYJ)$, so that $Y\in P$.  This completes the proof. \epf

\begin{rem}
As always, the condition $\rank(JX)=\rank(X)$, for $X\in\Mmn$, is equivalent to saying that rows $\row_{r+1}(X),\ldots,\row_m(X)$ belong to $\Span\{\row_1(X),\ldots,\row_r(X)\}$, with a dual statement holding for the condition $\rank(XJ)=\rank(X)$.  The regular elements of the corresponding semigroups of linear transformations are given by 
\begin{align*}
\Reg(\K_m(r)) &= \set{\al\in\K_m(r)}{\im(\al)\cap W_{mr}^\perp=\{0\}}, \\
\Reg(\I_n(r)) &= \set{\al\in\I_n(r)}{\im(\al)=\al(W_{nr})}.
\end{align*}
\end{rem}

Putting together all the above, we have proved the following.  (In the following statement, we slightly abuse notation by still denoting the map $\mathcal C_m(r)=\Mmn J\to\M_r$ by $\Phi_1$ and so on.)

\ms
\begin{thm}\label{thm:diamondsMmnJ}
We have the following commutative diagrams of semigroup epimorphisms:
\[
\includegraphics{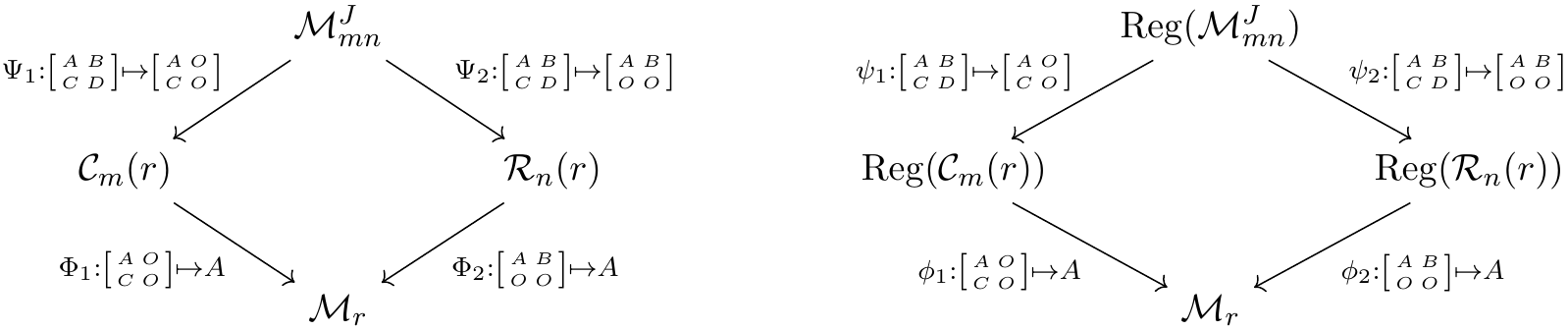}
\]
~ \epfres
\end{thm}

The remaining results of this section concern the regular subsemigroup $P=\Reg(\MmnJ)$.  From now on, we denote by $\phi=\phi_1\circ\psi_1=\phi_2\circ\psi_2$ the induced epimorphism $\phi:P\to\M_r$.  Also, for $X=\tmat ABCD\in P$, we write $\Xb=\phi(X)=A$.  The next result shows how the second commutative diagram from Theorem \ref{thm:diamondsMmnJ} may be used to identify $\RegMmnJ$ as a special kind of subdirect product of $\Reg(\Cmr)$ and $\Reg(\Rnr)$.

\ms
\begin{prop}\label{mono_prop}
There is an embedding
\[
\psi : \RegMmnJ\to \Reg(\Cmr)\times\Reg(\Rnr): X\mt (\psi_1(X),\psi_2(X))=(XJ,JX).
\]
As such, $P=\RegMmnJ$ is (isomorphic to) a pullback product of $PJ=\Reg(\Cmr)$ and $JP=\Reg(\Rnr)$.  Namely,
\[
P\cong\im(\psi)=\bigset{(\psi_1(X),\psi_2(X))}{X\in P} = \bigset{(Y,Z)\in PJ\times JP}{\phi_1(Y)=\phi_2(Z)}.
\]
\end{prop}

\pf Clearly, $\psi$ is a homomorphism.  Now let $X=[M,A,N]$ and $Y=[K,B,L]$ be elements of $P$ 
with $\psi(X)=\psi(Y)$.  Then
\[
\left(\tmat{A}{O}{MA}{O},\tmat A{AN}OO \right) = \psi(X) = \psi(Y) =  \left(\tmat{B}{O}{KB}{O},\tmat B{BL}OO \right).
\]
Comparing various coordinates, we deduce $A=B$, $MA=KB$ and $AN=BL$, giving 
$X=\tmat A{AN}{MA}{MAN}=\tmat B{BL}{KB}{KBL}=Y$, completing the proof that $\psi$ is injective.

To prove the statement about $\im(\psi)$, let $X\in P$ and put $Y=\psi_1(X)=XJ$ and $Z=\psi_2(X)=JX$.  Then $\phi_1(Y)=JY=JXJ=ZJ=\phi_2(Z)$.  Conversely, suppose $(Y,Z)\in PJ\times JP$ satisfies $JY=ZJ$.  Say $Y=UJ$ and $Z=JV$, where $U=[M,A,N]$ and $V=[K,B,L]$ belong to $P$.  Then $JY=JUJ=\tmat AOOO$ and $ZJ=JVJ=\tmat BOOO$, giving $A=B$.  But then $(Y,Z)=\psi(X)$, where $X=[M,A,L]\in P$.  \epf

\begin{rem}
We note that the previous result does not lift to a similar identification of $\MmnJ$ as a pullback product of $\Cmr$ and $\Rnr$ because the induced map
\[
\Psi:\MmnJ\to\Cmr\times\Rnr:X\mt(\Psi_1(X),\Psi_2(X))=(XJ,JX)
\]
is not injective.  Indeed, if $X=\tmat ABCD\in\Mmn$, then $\Psi(X)=\left(\tmat{A}{O}{C}{O},\tmat A{B}OO \right)$, with $\tmat ABCE$ mapping to the same pair for any other $E\in\M_{m-r,n-r}$.
\end{rem}

\ms
\begin{rem}
More generally, given a partial semigroup $(S,\cdot,I,\lam,\rho)$, the epimorphisms $\Psi_1$ and $\Psi_2$ from Theorem \ref{thm:diamonds}(v) allow for the definition of a map
\[
\Psi:(S_{ij},\star_a) \to (S_{ij}a,\cdot)\times (aS_{ij},\cdot):x\mt(xa,ax).
\]
To say that $\Psi$ is injective is to say that, for all $x,y\in S_{ij}$, $xa=ya$ and $ax=ay$ together imply $x=y$.  Compare this to the notion of a \emph{weakly reductive} semigroup $S$, in which, for every $x,y\in S$, the assumption that $xa=ya$ and $ax=ay$ for all $a\in S$ implies $x=y$.  See for example \cite[Definition 1.42]{Nagy2001}.
\end{rem}

We conclude this section with a simple but important observation that shows that $P=\RegMmnJ$ is a homomorphic image of the direct product of a rectangular band by the (non-sandwich) matrix semigroup~$\M_r$.  (Recall that a \emph{rectangular band} is a semigroup of the form $S\times T$ with product $(s_1,t_1)(s_2,t_2)=(s_1,t_2)$.)  Its proof is routine, relying on Proposition~\ref{prop:P1P2} and the rule $[M,A,N]\star[K,B,L] = [M,AB,L]$.  For the statement, recall that the \emph{kernel} of a semigroup homomorphism $\phi:S\to T$ (not to be confused with the kernel of a linear transformation) is the congruence $\ker(\phi)=\set{(x,y)\in S\times S}{\phi(x)=\phi(y)}$.  (A congruence on a semigroup $S$ is an equivalence relation $\sim$ for which $x_1\sim y_1$ and $x_2\sim y_2$ together imply $x_1x_2\sim y_1y_2$ for all $x_1,x_2,y_1,y_2\in S$; the quotient $S/{\sim}$ of all $\sim$-classes is a semigroup under the induced operation.  The first homomorphism theorem for semigroups states that any semigroup homomorphism $\phi:S\to T$ induces an isomorphism $S/\ker(\phi)\cong \im(\phi)$.)

\ms
\begin{prop}\label{prop:MAN}
Consider the semigroup $U=\M_{m-r,r}\times\M_r\times\M_{r,n-r}$ under the operation $\diamond$ defined by
\[
(M,A,N)\diamond(K,B,L) = (M,AB,L).
\]
Define an equivalence $\sim$ on $U$ by
\[
(M,A,N) \sim (K,B,L) \ \ \iff \ \ \text{$A=B$, $MA=KB$ and $AN=BL$.}
\]
Then $\sim$ is a congruence on $U$, and the map
\[
\xi:U\to P=\RegMmnJ:(M,A,N)\mt[M,A,N]=\mat A{AN}{MA}{MAN}
\]
is an epimorphism with $\ker(\xi)={\sim}$.  In particular, $P\cong U/{\sim}$. \epfres
\end{prop}

\section{The regular subsemigroup}\label{sect:RegMmnJ}

In this section, we continue to study the subsemigroup $P=\RegMmnJ$
consisting of all regular elements of~$\MmnJ$.  Eggbox diagrams of $P=\Reg(\M_{43}^J(\Z_2))$ are given in Figure \ref{fig:R} for values of $0\leq\rank(J)\leq3$; more examples can be seen by inspecting the regular $\gDJ$-classes in Figures~\ref{fig:V3212_V3322} and \ref{fig:V2322_V2422}.  
\begin{figure}[ht]
\begin{center}
\rotatebox[origin=c]{270}{
\includegraphics[width=0.55cm]{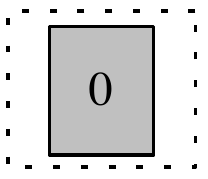}
\qquad
\includegraphics[width=1.3cm]{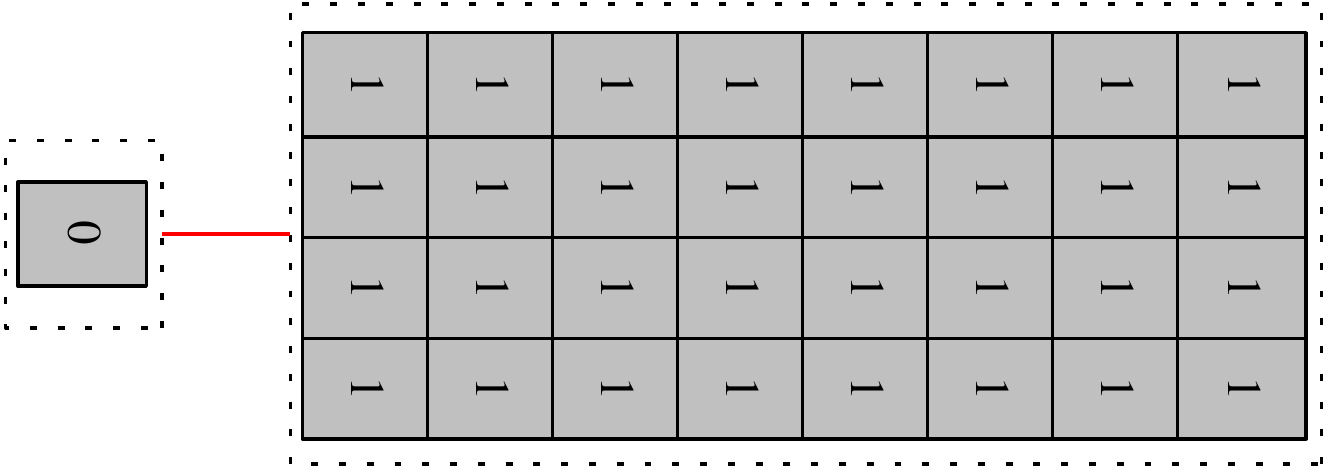}
\qquad
\includegraphics[width=1.85cm]{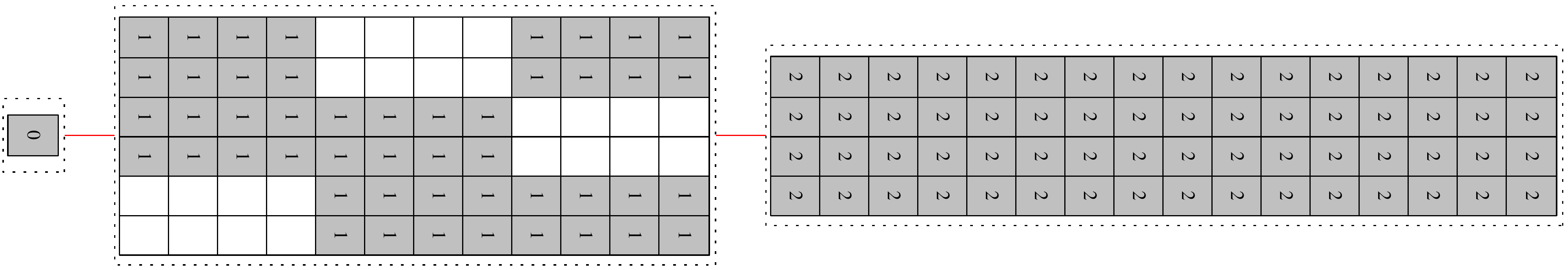}
\qquad
\includegraphics[width=2cm]{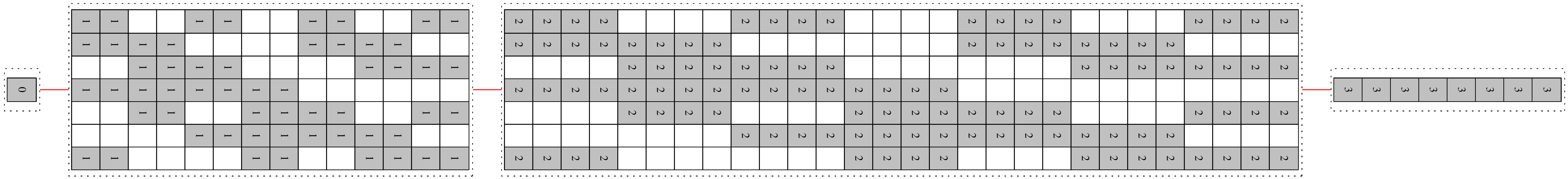}
}
\vspace{-4.5cm}
    \caption{Egg box diagrams (drawn sideways) of the regular linear sandwich semigroups $P=\Reg(\M_{43}^{J}(\Z_2))$, where $\rank(J)=0,1,2,3$ (top to bottom).}
    \label{fig:R}
   \end{center}
 \end{figure}
Comparing Figure \ref{fig:R} with Figure \ref{fig:M0...M3}, which pictures the full linear monoids $\M_r(\Z_2)$ for $0\leq r\leq 3$, an interesting pattern seems to emerge: namely, that $P=\Reg(\M_{43}^J(\Z_2))$ appears to be a kind of ``inflation'' of $\M_r$, where $r=\rank(J)$.  One of the goals of this section is to explain this phenomenon, and we do so by further exploring the map $$\phi:P\to\M_r:X= [M,A,N]\mt \Xb=A$$ defined after Theorem \ref{thm:diamondsMmnJ}.
We also calculate $|P|$, $\rank(P)$, and the number and sizes of various Green's classes.  As before, we assume that $J=J_{nmr}=\tmat {I_r}OOO\in\M_{nm}$.  Since $\MmnJ$ is just a zero semigroup if $r=0$, we generally assume that $r\geq1$.

Now, Theorem \ref{green_thm} enables us to immediately describe Green's relations on $P=\RegMmnJ$.  Since~$P$ is a regular subsemigroup of $\MmnJ$, the $\gR$, $\gL$, $\gH$ relations on $P$ are just the restrictions of the corresponding relations on $\MmnJ$ (see for example \cite{Hig,Howie}), and it is easy to check that this is also true for the $\gD=\gJ$ relation in this case.  So if $X\in P$ and $\gK$ is one of $\gR$, $\gL$, $\gH$, $\gD$, we will continue to write $\gKJ$ for the $\gK$ relation on~$P$, and $K_X^J$ for the $\gKJ$-class of $X$ in $P$.  
Parts (i--iv) of the next result also appear in \cite[Theorem~2.3]{Chinram2009}.


\ms
\begin{cor}\label{cor_green}
If $X\in P$, then 
\bit
\itemit{i} $R_X^J =R_X\cap P = \set{Y\in P}{\Col(X)=\Col(Y)}$, 
\itemit{ii} $L_X^J  =L_X\cap P = \set{Y\in P}{\Row(X)=\Row(Y)}$, 
\itemit{iii} $H_X^J   = H_X\cap P=H_X = \set{Y\in P}{\Col(X)=\Col(Y)\text{\emph{ and }}\Row(X)=\Row(Y)}$, 
\itemit{iv} $D_X^J  =D_X\cap P = \set{Y\in P}{\rank(X)=\rank(Y)}$.
\eit
The $\gDJ$-classes of $P$ form a chain: $D_0^J<\cdots<D_r^J$, where $D_s^J=\set{X\in P}{\rank(X)=s}$ for each $0\leq s\leq r$.~\epfres
\end{cor}

Also, the regularity of $P$ means that $P$ inherits the stability property from $\MmnJ$.  
The next result gives some combinatorial information about the size of $P$, and of various Green's classes in~$P$, in the case that $\F$ is finite.  Recall that $\{e_{k1},\ldots,e_{kk}\}$ is the standard basis of $V_k=\F^k$ and that $W_{ks}=\Span\{e_{k1},\ldots,e_{ks}\}$ for each $0\leq s\leq k$.

\ms
\begin{prop}\label{prop:DXJ_combinatorics}
Suppose $|\F|=q<\infty$.  Let $X\in P$ with $\rank(X)=s$.  Then
\bit
\itemit{i} $|R_X^J|=q^{s(n-r)}q^{{s\choose2}}(q-1)^s\qfact s\tqbin rs$,
\itemit{ii} $|L_X^J|=q^{s(m-r)}q^{{s\choose2}}(q-1)^s\qfact s\tqbin rs$,
\itemit{iii} $|H_X^J|=|\G_s|= q^{{s\choose2}}(q-1)^s\qfact s$,
\itemit{iv}$D_X^J=D_s^J$ is the union of:
\begin{itemize}\begin{multicols}{3}
\itemit{a} $q^{s(m-r)}\tqbin rs$ $\gRJ$-classes, 
\itemit{b} $q^{s(n-r)}\tqbin rs$ $\gLJ$-classes, 
\itemit{c} $q^{s(m+n-2r)}\tqbin rs^2$ $\gHJ$-classes,
\end{multicols}\end{itemize}
\itemit{v} $|D_X^J|=|D_s^J|=q^{s(m+n-2r)}q^{{s\choose2}}(q-1)^s\qfact s\tqbin rs^2$.
\eit
Consequently,
$
\ds{|P|=|\RegMmnJ|=\sum_{s=0}^r q^{s(m+n-2r)}q^{{s\choose2}}(q-1)^s\qfact s\tqbin rs^2}.
$
\end{prop}

\pf We start with (i).  Since $|R_X^J|=|R_Y^J|$ for all $Y\in D_X^J=D_s^J$, we may assume $X=J_{mns}$.  Now, $\Col(X)=W_{ms}$.  By Proposition~\ref{prop:P1P2} and Corollary \ref{cor_green}, we have
\begin{align*}
R_X^J &= \set{Y\in P}{Y\gRJ X} \\
&= \set{Y\in\Mmn}{\Col(X)=\Col(Y),\ \Col(Y)=\Col(YJ),\ \Row(Y)=\Row(JY)} \\
&= \set{Y\in\Mmn}{\Col(Y)=\Col(YJ)=W_{ms}},
\end{align*}
since if $\Col(Y)=W_{ms}$, then $Y$ is of the form $Y=\tmat AB{O_{m-s,r}}{O_{m-s,n-r}}$ for some $A\in\M_{sr}$ and $B\in\M_{s,n-r}$, in which case $JY=\tmat AB{O_{n-s,r}}{O_{n-s,n-r}}$ automatically has the same row space as $Y$.  

Now consider some $Y\in R_X^J$.  As noted above, we must have $Y=\tmat ABOO$ for some $A\in\M_{sr}$ and $B\in\M_{s,n-r}$.  Since $YJ=\tmat AOOO$, the condition $\Col(YJ)=W_{ms}$ is equivalent to $\Col(A)=V_s$.  In particular, there is no restriction on the entries of $B$, so $B$ may be chosen (arbitrarily, and independently of $A$) in $q^{s(n-r)}$ ways.  Also, $\dim(\Row(A))=\dim(\Col(A))=s$.  So $A$ may be specified by listing its rows (in order), which are $s$ linearly independent row vectors from $\F^r$.  The number of possible choices for $A$ is therefore $(q^r-1)(q^r-q)\cdots(q^r-q^{s-1})=q^{{s\choose2}}(q-1)^s\qfact s\tqbin rs$.  Multiplying these two values gives (i).

Part (ii) is dual to (i).  
%
%
%
%
%
Part (iii) follows directly from Corollary \ref{cor_green}(iii) and Lemma \ref{lem:combinatorics_Mmn}(iii).  Parts (a) and (b) of (iv) follow by dividing $|L_X^J|$ and $|R_X^J|$ by $|H_X^J|$, respectively.  Part (c) follows from (a) and (b).  Part (v) follows from (iii) and part (c) of (iv).
The formula for $|P|$ is obtained by adding the sizes of the $\gDJ$-classes.~\epf

Recall that, for $X=[M,A,N]\in P$, we write $\Xb=\phi(X)=A\in\M_r$.  We extend this notation to subsets of~$P$, so if $\Om\sub P$, we write $\Omb=\set{\Xb}{X\in\Om}$.
We now show how the epimorphism $\phi:P\to\M_r$ may be used to relate Green's relations on the semigroups $P$ and $\M_r$.
If $X,Y\in P$ and $\gK$ is one of $\gR$, $\gL$, $\gH$, $\gD$, we say $X\gKh Y$ if $\Xb\gK\Yb$ (in~$\M_r$).  Denote by $\Kh_X=\phi^{-1}(K_{\Xb})=\set{Y\in P}{X\gKh Y}$ the $\gKh$-class of $X$ in~$P$.  We first need a technical result.  

\ms
\begin{lemma}\label{lem:|RhX|}
Let $X,Y\in P$.  If $X\gDh Y$, then $|\Rh_X|=|\Rh_Y|$ and $|\Lh_X|=|\Lh_Y|$.
\end{lemma}

\pf By duality, it suffices to prove the statement about $\gRh$-classes.  Now, $X\gDh Y$ means that $X\gLh W\gRh Y$ for some $W\in P$.  Since $\Rh_Y=\Rh_W$, we may assume without loss of generality that $X\gLh Y$.
Write $X=[M,A,N]$ and $Y=[K,B,L]$.  By definition $X\gLh Y$, means that $A\L B$ in $\M_r$, so $A=UB$ for some $U\in\G_r$ by Lemma \ref{lem:greenMmn}.  Now let $Z=[M',A',N']\in\Rh_X$, and define $\al(Z)=[M'U,U^{-1}A',N']$.  It is easy to check that for any other representation $Z=[M'',A'',N'']$, we have $[M''U,U^{-1}A'',N'']=[M'U,U^{-1}A',N']$, so that $\al(Z)$ is well-defined.  Also,
\begin{align*}
Z\gRh X\ \implies\ A'\R A 
\ \implies\ U^{-1}A'\R U^{-1}A =B 
\ \implies\ \al(Z) \gRh Y.
\end{align*}
Thus $\al$ is a map $\al:\Rh_X\to\Rh_Y$.  It is easy to check that $\Rh_Y\to\Rh_X:[M',A',N']\mt[M'U^{-1},UA',N']$ is the inverse mapping of $\al$.  We conclude that $|\Rh_X|=|\Rh_Y|$. \epf

For the proof of the next result, we note that stability of $\M$ implies that $A^2\D A \iff A^2\H A$ for all $A\in\M_r$.  We also use the fact that an $\H$-class $H$ of a semigroup $S$ is a group if and only if $x^2\in H$ for some (and hence for all) $x\in H$ \cite[Theorem 2.2.5]{Howie}.
Recall that a $k\times l$ \emph{rectangular band} is a semigroup of the form $S\times T$ with product $(s_1,t_1)(s_2,t_2)=(s_1,t_2)$, where $|S|=k$ and $|T|=l$.  A $k\times l$ \emph{rectangular group} with respect to a group $G$ is a direct product of a $k\times l$ rectangular band with $G$.

For the proof of the next result (and elsewhere), it will be convenient to define a number of equivalence relations.  For $A\in\M_r$, we define equivalences $\sim_A$ and $\approx_A$ on $\M_{m-r,r}$ and $\M_{r,n-r}$ (respectively) by
\[
M_1\sim_AM_2 \ \iff \ M_1A=M_2A \AND N_1\approx_AN_2 \ \iff \ AN_1=AN_2.
\]



\ms
\begin{thm}\label{inflation_thm}
Suppose $|\F|=q<\infty$.  Let $X\in P=\RegMmnJ$ and put $s=\rank(X)$. 
\bit
\itemit{i} $\Rh_X$ is the union of $q^{s(m-r)}$ $\gRJ$-classes of $P$.
\itemit{ii} $\Lh_X$ is the union of $q^{s(n-r)}$ $\gLJ$-classes of $P$.
\itemit{iii} $\Hh_X$ is the union of $q^{s(m+n-2r)}$ $\gHJ$-classes of $P$, each of which has size $|\G_s|=q^{{s\choose2}} (q-1)^s\qfact s$.  The map $\phi:P\to\M_r$ is injective when restricted to any $\gHJ$-class of $P$.
\itemit{iv} If $H_{\Xb}$ is a non-group $\gH$-class of $\M_r$, then each $\gHJ$-class of $P$ contained in $\Hb_X$ is a non-group.
\itemit{v} If $H_{\Xb}$ is a group $\gH$-class of $\M_r$, then each $\gHJ$-class of $P$ contained in $\Hb_X$ is a group isomorphic to $\G_s$; further, $\Hb_X$ is a $q^{s(m-r)}\times q^{s(n-r)}$ rectangular group with respect to $\G_s$.
\itemit{vi} $\gDh=\gDJ$ and $\Dh_X=D_X^J=D_s^J=\set{Y\in P}{\rank(Y)=s}$ is the union of:
\begin{itemize}
\itemit{a} $\tqbin rs$\ $\gRh$-classes (and the same number of $\gLh$-classes) of $P$,
\itemit{b} $q^{s(m-r)}\tqbin rs$\ $\gRJ$-classes of $P$,
\itemit{b} $q^{s(n-r)}\tqbin rs$\ $\gLJ$-classes of $P$,
\itemit{d} $\tqbin rs^2$\ $\gHh$-classes of $P$,
\itemit{e} $q^{s(m+n-2r)}\tqbin rs^2$\ $\gHJ$-classes of $P$.
\eitres
\eitres
\end{thm}

\pf First observe that if $\rho:S\to T$ is an epimorphism of semigroups, and if $K$ is a $\gK$-class of $T$ where $\gK$ is one of $\gR$, $\gL$, $\gH$, then $\rho^{-1}(K)$ is a union of $\gK$-classes of $S$.  Throughout the proof, we write
\[
X=[M,A,N]=\mat A{AN}{MA}{MAN},
\]
so $A\in\M_r$ satisfies $\rank(A)=\rank(JXJ)=\rank(X)=s$.  We note that $\gDh=\gDJ$ immediately follows.
\bit
\item[(i)] By the first observation, it suffices to count the number of $\gRJ$-classes contained in $\Rh_X$.  
Since $|\Rh_X|=|\Rh_Y|$ for all $Y\in D_s^J$ by Lemma \ref{lem:|RhX|}, it follows that each $\gRh$-class of $D_s^J$ contains the same number of $\gRJ$-classes.  By Lemma \ref{lem:combinatorics_Mmn}, $D_s(\M_r)$ is the union of $\tqbin rs$ $\R$-classes (and the same number of $\L$-classes), so it follows that $D_s^J$ is the union of $\tqbin rs$ $\gRh$-classes (and the same number of $\gLh$-classes).  By Proposition \ref{prop:DXJ_combinatorics}, $D_s^J$ is the union of $q^{s(m-r)}\tqbin rs$ $\gRJ$-classes.  Dividing these, it follows that each $\gRh$-class of $D_s^J$ is the union of $q^{s(m-r)}$ $\gRJ$-classes.

\item[(ii)] This is dual to (i).

\item[(iii)] The statement concerning the number of $\gHJ$-classes contained in $\Hb_X$ follows immediately from (i) and~(ii), and the size of these $\gHJ$-classes was given in Proposition \ref{prop:DXJ_combinatorics}.  Next, for any $B\in\M_r$ with $B\H A$, it is easy to check that $[M,B,N] \gHJ X$.  So the set $\Om=\set{[M,B,N]}{B\in H_A}$ is contained in~$\gHJ_X$.  Since $|\Om|=|H_A|=|\G_s|=|H_X^J|$, we see that $H_X^J=\Om$.  For any $Z=[M,B,N]\in\Om=H_X^J$, we have $\phi(Z)=B$, so it follows that $\phi|_{H_X^J}$ is injective.

\item[(iv)] Suppose $H_A=H_{\Xb}$ is a non-group $\gH$-class of $\M_r$, and let $Y=[K,B,L]\in\Hh_X$ be arbitrary.  Since $Y\gHh X$, it follows that $B=\Yb\gH \Xb=A$.  Since $H_B=H_A$ is not a group, we have $B^2\not\in H_B$, whence $B^2\not\in D_B$ and $\rank(B^2)<\rank(B)=\rank(A)=s$.  But then $Y^2=[K,B^2,L]\not\in D_s^J=D_Y^J$, so that $Y^2\not\in H_Y^J$, and we conclude that $H_Y^J$ is not a group.

\item[(v)] Suppose $H_{\Xb}$ is a group.  Then $\Yb^2\in H_{\Xb}$ for any $Y\in \Hh_X$, so $\rank(Y\star Y)=\rank(\Yb^2)=\rank(\Yb)=\rank(Y)$, giving $Y\star Y\gDJ Y$, so that $Y\star Y\in H_Y^J$ and $H_Y^J$ is a group.  By (iii), the restriction of $\phi$ to $H_Y^J$ yields an isomorphism onto $H_{\Yb}\cong\G_s$.

Let $E\in\M_r$ be the identity element of the group $H_A$.  
Let $\sM_E\sub\M_{m-r,r}$ (resp., $\sN_E\sub\M_{r,n-r}$) be a cross-section of the $\sim_E$-classes (resp., $\approx_E$-classes) in $\M_{m-r,r}$ (resp., $\M_{r,n-r}$).  
%
%
It is easy to check that every $Y\in\Hh_X$ may be uniquely represented as $Y=[K,B,L]$ for some $K\in\sM_E$, $B\in H_A$ and $L\in\sN_E$.  
It follows that the map
\[
\sM_E\times H_A\times\sN_E\to\Hh_A:(K,B,L)\mt[K,B,L]
\]
is a well-defined isomorphism, where the (rectangular group) product on $\sM_E\times H_A\times\sN_E$ is defined by $(K_1,B_1,L_1)\cdot(K_2,B_2,L_2)=(K_1,B_1B_2,L_2)$.  We have already observed that $H_A\cong\G_s$, and the dimensions of the rectangular band $\sM_E\times\sN_E$ follow from parts (i--iii) together with the observation that $\set{[K,B,L]}{B\in H_A}$ is an $\gHJ$-class contained in $\Hh_X$ for each $K\in\sM_E$ and $L\in\sN_E$.

\item[(vi)] We have already noted that $\gDh=\gDJ$.  We proved (a) while proving (i), above.  Parts (b), (c) and (e) were proved in Proposition \ref{prop:DXJ_combinatorics}.  Part (d) follows from (a).
\epf
\eit


The previous result explains the ``inflation'' phenomenon discussed at the beginning of this section; see also Figure \ref{fig:R}.  As an immediate corollary of Theorem \ref{inflation_thm}, we may now completely classify the isomorphism classes of finite linear sandwich semigroups.



\ms
\begin{thm}\label{thm:classification}
Let $\F_1$ and $\F_2$ be two finite fields with $|\F_1|=q_1$ and $|\F_2|=q_2$, let $m,n,k,l\geq1$, and let $A\in D_r(\M_{nm})$ and $B\in D_s(\M_{lk})$. 
The following are equivalent:
\bit
\itemit{i} $\MMN\cong\MKL$,
\itemit{ii} one of the following holds:
\begin{itemize}
\itemit{a} $r=s=0$ and $q_1^{mn}=q_2^{kl}$, or
\itemit{b} $r=s\geq1$, $(m,n)=(k,l)$, and $q_1=q_2$.
\end{itemize}
\eit
Further, if $r\geq1$, then $\MMN\cong\MKL$ if and only if $\RegMMN\cong\RegMKL$.
\end{thm}

\pf Again, if $r\not=s$, then counting the regular $\D^A$- and $\D^B$-classes shows that $\MMN\not\cong\MKL$.  For the remainder of the proof, we assume $r=s$.
Suppose first that $r=s=0$.  Then $\MMN$ and $\MKL$ are both zero semigroups and so are isomorphic if and only if their sizes, $q_1^{mn}$ and $q_2^{kl}$, are equal.  
For the remainder of the proof, we assume $r=s\geq1$, and write $D_t^A$ and $D_t^B$ for the relevant regular $\D^A$- and $\D^B$-classes in $\MMN$ and $\MKL$ for each $0\leq t\leq r=s$.  

By Theorem \ref{inflation_thm}(v), any group $\H^A$-class contained in $D_1^A$ is isomorphic to $\G_1(\F_1)\cong\F_1^\times$, the multiplicative group of $\F_1$.  Since $|\F_1^\times|=q_1-1$ and $|\F_2^\times|=q_2-1$, it follows that if $q_1\not=q_2$, then $\RegMMN\not\cong\RegMKL$ and, hence, $\MMN\not\cong\MKL$.  
Now suppose $q_1=q_2$ (so $\F_1\cong\F_2$), and write $q=q_1$.  By Theorem~\ref{inflation_thm}(vi), $D_1^A$ (resp., $D_1^B$) contains $q^{m-r}\tqbin r1$ $\R^A$-classes (resp., $q^{k-r}\tqbin r1$ $\R^B$-classes).  
It follows that if $m\not=k$ (or, dually, if $n\not=l$), then $\RegMMN\not\cong\RegMKL$ and, hence, $\MMN\not\cong\MKL$.  
Conversely, if (b) holds, then $\MMN\cong\MKL$ by Lemma \ref{lem:MmnAMmnB}(ii).

For the final statement, first note that $\MMN\cong\MKL$ clearly implies $\RegMMN\cong\RegMKL$.  In the previous paragraph, we showed that the negation of (b) implies $\RegMMN\not\cong\RegMKL$.  This completes the proof. \epf

\begin{rem}
Of course, if $\rank(A)=\rank(B)=0$, then $\RegMMN=\{O_{mn}\}\cong\RegMKL=\{O_{kl}\}$, regardless of $m,n,k,l,q_1,q_2$.  So the final clause of Theorem \ref{thm:classification} does not hold for $r=0$.
\end{rem}

\ms
\begin{rem}\label{rem:infinite_classification}
The infinite case is not as straight-forward, since $|\M_{mn}(\F)|=|\F|$ for all $m,n\geq1$, and since it is possible for two non-isomorphic fields $\F_1,\F_2$ to have isomorphic multiplicative groups $\F_1^\times,\F_2^\times$ (for example, $\mathbb Q$ and $\mathbb Z_3(x)$ both have multiplicative groups isomorphic to $\mathbb Z_2\oplus F$, where $F$ is a free abelian group of countably infinite rank).  So we have the following isomorphisms:
\bit
\item[(i)] $\MMN\cong\MKL$ if $m,n,k,l\geq1$, $|\F_1|=|\F_2|$, and $\rank(A)=\rank(B)=0$ --- indeed, both sandwich semigroups are zero semigroups of size $|\F_1|=|\F_2|$;
\item[(ii)] $\MMN\cong\M_{mn}^B(\F_2)$ if $\F_1^\times\cong\F_2^\times$ and $\rank(A)=\rank(B)=1$ --- indeed, 
%
when $J=J_{nm1}=\tmat{I_1}OOO$, sandwich products $X\star_J Y$ involve only field multiplication and no addition:
\[
\left[
\begin{matrix}
a_{11} & a_{12} & \cdots & a_{1n} \\
a_{21} & a_{22} & \cdots & a_{2n} \\
\vdots & \vdots &\ddots &\vdots \\
a_{m1} & a_{m2} & \cdots & a_{mn} \\
\end{matrix}
\right]
\star_J
\left[
\begin{matrix}
b_{11} & b_{12} & \cdots & b_{1n} \\
b_{21} & b_{22} & \cdots & b_{2n} \\
\vdots & \vdots &\ddots &\vdots \\
b_{m1} & b_{m2} & \cdots & b_{mn} \\
\end{matrix}
\right]
=
\left[
\begin{matrix}
a_{11}b_{11} & a_{11}b_{12} & \cdots & a_{11}b_{1n} \\
a_{21}b_{11} & a_{21}b_{12} & \cdots & a_{21}b_{1n} \\
\vdots & \vdots &\ddots &\vdots \\
a_{m1}b_{11} & a_{m1}b_{12} & \cdots & a_{m1}b_{1n} \\
\end{matrix}
\right].
\]
\eitres
We leave it as an open problem to completely classify the isomorphism classes of linear sandwich semigroups over infinite fields.  But we make two simple observations:
\bit
\item[(iii)] as in the proof of Theorem \ref{thm:classification}, if $\MMN\cong\M_{kl}^B(\F_2)$, then we must have $\rank(A)=\rank(B)$;
\item[(iv)] if $\MMN\cong\M_{kl}^B(\F_2)$ with $\rank(A)=\rank(B)=r\geq2$, we must have $\F_1\cong\F_2$ (since the maximal subgroups of $\MMN$ are isomorphic to $\G_s(\F_1)$ for $0\leq s\leq r$, and since $\G_s(\F_1)\cong\G_s(\F_2)$ implies $\F_1\cong\F_2$ for $s\geq2$ \cite{Dieudonne1971}).
\eit
\end{rem}

In what follows, the top $\gDJ$-class of $P=\RegMmnJ$ plays a special role.  We write $D$ for this $\gDJ$-class, so
\[
D=D_r^J = \phi^{-1}(\G_r) = \set{X\in P}{\rank(X)=r}.
\]
As a special case of Theorem \ref{inflation_thm}(v), $D$ is a $q^{r(m-r)}\times q^{r(n-r)}$ rectangular group with respect to $\G_r$.  Since~$D$ is the pre-image of $\G_r$ under the map $\phi:P\to\M_r$, we may think of $D$ as a kind of ``inflation'' of $\G_r$, the group of units of $\M_r$.  In fact, more can be said along these lines.  Recall again that the \emph{variant} of a semigroup~$S$ with respect to an element $a\in S$ is the semigroup $S^a$ with underlying set $S$ and operation $\star_a$ defined by $x\star_ay=xay$ for all $x,y\in S$.  
Recall also that an element $a\in S$ of a (necessarily regular) semigroup $S$ is \emph{regularity preserving} if the variant $S^a$ is regular.  The set $\RP(S)$ of all regularity preserving elements of $S$ was studied in \cite{KL2001,Hickey1983}; we will not go into the details here, but it was explained in \cite{KL2001} that $\RP(S)$ is a useful alternative to the group of units in the case that $S$ is not a monoid (as with $P$ when $r=m=n$ does not hold).  Because of this, it is significant that $D$ is equal to $\RP(P)$, the set of all regularity preserving elements of $P=\RegMmnJ$, as we will soon see.  We now state a result from \cite{KL2001} concerning regularity preserving elements.  Recall that an element $u$ of a semigroup $S$ is a \emph{mididentity} if $xuy=xy$ for all $x,y\in S$ \cite{Yamada1955}; of course for such an element, $\star_u$ is just the original semigroup operation.  

\ms
\begin{prop}[Khan and Lawson \cite{KL2001}]\label{RPS_prop}
Let $S$ be a regular semigroup.  
\bit
\itemit{i} An element $a\in S$ is regularity preserving if and only if $a\gH e$ for some regularity preserving idempotent $e\in E(S)$. (In particular, $\RP(S)$ is a union of groups.)
\itemit{ii} An idempotent $e\in E(S)$ is regularity preserving if and only if $fe \gR f \gL ef$ for all idempotents $f\in E(S)$.
\itemit{iii} Any mididentity is regularity preserving. \epfres
\eitres
\end{prop}


In order to avoid confusion when discussing idempotents, if $\Om\sub\Mmn$, we will write
\[
\EJ(\Om)=\set{X\in\Om}{X=X\star X}
\]
for the set of idempotents from $\Om$ with respect to the $\star$ operation on $\MmnJ$.  If $\Si\sub \M_k$ for some $k$, we will continue to write $E(\Si)=\set{A\in\Si}{A=A^2}$ for the set of idempotents from $\Si$ with respect to the usual matrix multiplication.  
%
%
%

\ms\ms
\begin{lemma}\label{lemma_aea}
\bit
\itemit{i} $\EJ(\MmnJ)=\EJ(P)=\set{[M,A,N]}{A\in E(\M_r),\ M\in\M_{m-r,r},\ N\in\M_{r,n-r}}$.
\itemit{ii} $\EJ(D)=\set{[M,I_r,N]}{M\in\M_{m-r,r},\ N\in\M_{r,n-r}}$ is a $q^{r(m-r)}\times q^{r(n-r)}$ rectangular band.
\itemit{iii} Each element from $\EJ(D)$ is a mididentity for both $\MmnJ$ and $P$.
\itemit{iv} $D=\RP(P)$ is the set of all regularity-preserving elements of $P$.
\eitres
\end{lemma}

\pf Note that all idempotents are regular.  If $X=[M,A,N]\in P$, then $X\star X=[M,A^2,N]$, so $X=X\star X$ if and only if $A=A^2$, giving~(i).  Part (ii) follows from (i), 
since $I_r$ is the only idempotent from the group~$\G_r=D_r(\M_r)$.  Using (ii), it is easy to check by direct computation that $X\star Y\star Z=X\star Z$ for all $X,Z\in \Mmn$ and $Y\in\EJ(D)$, giving (iii).  Finally, to prove (iv), note that by Proposition \ref{RPS_prop}(i), it suffices to show that $\EJ(\RP(P))=\EJ(D)$.  By (iii) and Proposition \ref{RPS_prop}(iii), we have $\EJ(D)\sub\EJ(\RP(P))$.  Conversely, suppose $X\in\EJ(\RP(P))$.  Let $Y\in\EJ(D)$.  By Proposition \ref{RPS_prop}(ii), and the fact that $\gLJ\sub\gDJ$, $X\star Y\gDJ Y$.  It follows that $r=\rank(Y)=\rank(XJY)\leq\rank(X)\leq r$, giving $\rank(X)=r$, and $X\in\EJ(D)$.
This shows that $\EJ(\RP(P))\sub \EJ(D)$, and completes the proof.~\epf

We may now calculate the rank of $P=\RegMmnJ$ in the case of finite $\F$.  For the following proof, recall from \cite{HHR} that the \emph{relative rank} $\rank(S:U)$ of a semigroup $S$ with respect to a subset $U\sub S$ is defined to be the minimum cardinality of a subset $V\sub S$ such that $S=\la U\cup V\ra$.  

\ms
\begin{thm}\label{thm:rankP}
Suppose $|\F|=q<\infty$.  If $1\leq r\leq\min(m,n)$ and we do not have $r=m=n$, then $$\rank(P)=\rank(\RegMmnJ)=q^{r(L-r)}+1,$$ where $L=\max(m,n)$.
\end{thm}

\pf Since $D$ is a subsemigroup of $P$ and $P\sm D$ is an ideal, it quickly follows that $\rank(P)=\rank(D)+\rank(P:D)$.  It is well-known \cite{Ruskuc1994} that a rectangular group $R=(S\times T)\times G$ satisfies $\rank(R)=\max\big\{|S|,|T|,\rank(G)\big\}$.  Since $D$ is a $q^{r(m-r)}\times q^{r(n-r)}$ rectangular group with respect to $\G_r$, and since $\rank(\G_r)\leq2$ by Theorem~\ref{thm:waterhouse}, it immediately follows that $\rank(D)=q^{r(L-r)}$.  Since $\la D\raJ=D\not=P$ (as $r\geq1$), we have $\rank(P:D)\geq1$, so the proof will be complete if we can show that $P=\la D\cup\{X\}\raJ$ for some $X\in P$.  
With this in mind, let $X\in D_{r-1}^J$ be arbitrary.  Note that $\Db=\set{\Yb}{Y\in D}=\G_r$, and $\Xb\in D_{r-1}(\M_r)$.  It follows from Theorem~\ref{thm:waterhouse} that $\M_r=\la\Db\cup\{\Xb\}\ra$.  Now let $Y=[M,A,N]\in P$ be arbitrary.  Choose $Z_1,\ldots,Z_k\in D\cup\{X\}$ such that $A=\Zb_1\cdots\Zb_k$.  Then $Y=[M,I_r,N]\star Z_1\star\cdots\star Z_k\star[M,I_r,N]$, with $[M,I_r,N]\in D$. \epf 

\begin{rem}\label{rem:r=m3}
If $r=0$, then $P=\{O\}$, while if $r=m=n$, then $P=\M_n$.  So $\rank(P)$ is trivial in the former case, and well-known in the latter (see Theorem~\ref{thm:waterhouse}).  
%
As in Remark \ref{rem:r=m2}, we deduce that 
$
\rank(\Reg(\Ckl)) = \rank(\Reg(\Rkl)) = q^{l(k-l)} + 1
$
for $|\F|=q<\infty$.
\end{rem}

\section{The idempotent generated subsemigroup}\label{sect:EMmnJ}


In this section, we investigate the idempotent generated subsemigroup $\la\Ea(\MmnJ)\raa$ of $\MmnJ$; we write $\EmnJ$ for this idempotent generated subsemigroup.  Our main results include a proof that $\EmnJ=(P\sm D)\cup \Ea(D)$ and a calculation of $\rank(\EmnJ)$ and $\idrank(\EmnJ)$; in particular, we show that these two values are equal.
Since the solution to every problem we consider is trivial when $r=0$, and well-known when $r=m=n$, we will continue to assume that $r\geq1$ and that $r=m=n$ does not hold.  To simplify notation, we will write $E=\Ea(\MmnJ)=\Ea(P)$, so $\EmnJ=\la E\raa$.  
We begin by calculating $|E|$ in the case of finite $\F$, for which we need the following formulae for $|E(D_s(\M_r))|$.  Although the next result might already be known, we are unaware of a reference and include a simple proof for convenience.

\ms
\begin{lemma}\label{lem:EDmMr}
Suppose $|\F|=q<\infty$.  If $0\leq s\leq r$, then $|E(D_s(\M_r))|=q^{s(r-s)}\tqbin rs $.  Consequently,
\[
|E(\M_r)| = \sum_{s=0}^r q^{s(r-s)}\qbin rs.
\]
\end{lemma}

\pf To specify an idempotent endomorphism $\al\in\Endr$ of rank $s$, we first choose $W=\im(\al)$, which is a subspace of dimension $s$ and may be chosen in $\tqbin rs$ ways, and we note that $\al$ must map $W$ identically.  If $\{v_1,\ldots,v_r\}$ is an arbitrary basis for $V_r$, such that $\{v_1,\ldots,v_s\}$ is a basis of $W$, then $\al$ may map each of $v_{s+1},\ldots,v_r$ arbitrarily into $W$, and there are $(q^s)^{r-s}$ ways to choose these images. \epf

\begin{prop}\label{prop:enumeration_E}
Suppose $|\F|=q<\infty$.  If $0\leq s\leq r$, then $|\EJ(D_s^J)|=q^{s(m+n-r-s)}\tqbin rs$.  Consequently,
\[
|\EJ(\MmnJ)| = \sum_{s=0}^r q^{s(m+n-r-s)}\qbin rs.
\]
\end{prop}

\pf 
Parts (iv) and (v) of Theorem \ref{inflation_thm} say that an $\gHJ$-class $H_X^J\sub D_s^J$ is a group (so contains an idempotent) if and only if $H_{\Xb}$ is a group $\H$-class of $\M_r$, and that there are $q^{s(m-r)}\times q^{s(n-r)}$ idempotents of $D_s^J$ corresponding to each rank $s$ idempotent of $\M_r$, of which there are $q^{s(r-s)}\tqbin rs$ by Lemma \ref{lem:EDmMr}.  The result quickly follows. \epf

We now describe the idempotent generated subsemigroup of $\MmnJ$.

\ms
\begin{thm}\label{thm:EmnJ}
We have $\EmnJ=\la\EJ(\MmnJ)\raJ=(P\sm D)\cup\EJ(D)$.
\end{thm}

\pf Suppose $X_1,\ldots,X_k\in E=\EJ(\MmnJ)$, and write $X_i=[M_i,A_i,N_i]$ for each $i$.  So $A_i\in E(\M_r)$ for each $i$.  Then $X_1\star\cdots\star X_k = [M_1,A_1\cdots A_k,N_k]$.  If any of $A_1,\ldots,A_k$ belongs to $\MrGr$, then so too does $A_1\cdots A_k$, so that $X_1\star\cdots\star X_k\in P\sm D$.  If all of $A_1,\ldots,A_k$ belong to $\G_r$, then $A_1=\cdots=A_k=I_r$, so $X_1\star\cdots\star X_k=[M_1,I_r,N_k]\in\EJ(D)$.  This shows that $\EmnJ\sub(P\sm D)\cup\EJ(D)$.  Conversely, it suffices to show that $P\sm D\sub\EmnJ$, so suppose $X\in P\sm D$, and write $X=[M,A,N]$.  Since $X\not\in D$, we must have $\rank(A)=\rank(X)<r$.  But then $A\in\MrGr$, so that $A=B_1\cdots B_l$ for some $B_1,\ldots,B_l\in E(\M_r)$ by Theorem \ref{thm_MnGn}.  It follows that $X=[M,B_1,N]\star\cdots\star[M,B_l,N]$, with all $[M,B_i,N]\in E$. \epf 

\begin{rem}
Recall (see Theorem \ref{thm_MnGn}) that ${\E_n=\la E(\M_n)\ra = (\MnGn)\cup\{I_n\}}$.  Theorem \ref{thm:EmnJ} is a pleasing analogue of that result, since $\{I_n\}=E(\G_n)$, where $\G_n$ is the top $\gD$-class of $\M_n$.  Also, $\G_n=G(\M_n)=\RP(\M_n)$ and, while $P$ has no group of units as it is not a monoid, it is still the case that $D=\RP(P)$.
\end{rem}



Now that we have described the elements of the semigroup $\EmnJ$, the next natural task is to calculate its rank and idempotent rank.  



\ms
\begin{thm}\label{thm:rankEmnJ}
Suppose $|\F|=q<\infty$.  Then $$\rank(\EmnJ)=\idrank(\EmnJ)=q^{r(L-r)}+(q^r-1)/(q-1),$$ where $L=\max(m,n)$.
\end{thm}

\pf As in the proof of Theorem \ref{thm:rankP}, we have $\rank(\EmnJ)=\rank(\Ea(D))+\rank(\EmnJ:\Ea(D))$.  Since $\EJ(D)$ is a $q^{r(m-r)}\times q^{r(n-r)}$ rectangular band (see Lemma \ref{lemma_aea}(ii)), we again deduce from \cite{Ruskuc1994} that $\rank(\Ea(D))=\idrank(\Ea(D)) = q^{r(L-r)}$.  So it remains to show that:
\bit
\item[(i)] there exists a set $\Om\sub E$ of size $(q^r-1)/(q-1)$ such that $\EmnJ=\la\Ea(D)\cup \Om\raa$, and
\item[(ii)] if $\Si\sub P$ satisfies $\EmnJ=\la\Ea(D)\cup \Si\raa$, then $|\Si|\geq(q^r-1)/(q-1)$.
\eit
By Theorem \ref{thm_MnGn}, we may choose some set $\Ga\sub E(\M_r)$ with $\la\Ga\ra=\MrGr$ and $|\Ga|=(q^r-1)/(q-1)$.  For each $A\in\Ga$, choose any $M_A\in\M_{m-r,r}$ and $N_A\in\M_{r,n-r}$, and put $\Om=\set{[M_A,A,N_A]}{A\in\Ga}$.  Since $\EmnJ=(P\sm D)\cup\EJ(D)$, the proof of (i) will be complete if we can show that $P\sm D\sub\la\Ea(D)\cup \Om\raa$.  So let $X=[K,B,L]\in P\sm D$, and write $B=A_1\cdots A_k$ where $A_1,\ldots,A_k\in\Ga$.  Then
\[
X=[K,I_r,L]\star[M_{A_1},A_1,N_{A_1}]\star\cdots\star[M_{A_k},A_k,N_{A_k}]\star[K,I_r,L]\in\la\Ea(D)\cup \Om\raa,
\]
as required.  Next, suppose $\EmnJ=\la\Ea(D)\cup \Si\raa$, where $\Si\sub\EmnJ\sm \Ea(D)=P\sm D$.  We will show that $\Sib$ generates $\MrGr$.  Indeed, let $A\in\MrGr$ be arbitrary, and choose any $X\in P$ such that $\Xb=A$.  Since $\rank(X)=\rank(A)<r$, it follows that $X\in P\sm D\sub\EmnJ$.  Consider an expression $X=Y_1\star\cdots\star Y_k$, where $Y_1,\ldots,Y_k\in \Ea(D)\cup \Si$.  Now, $A=\Xb=\Yb_1\cdots\Yb_k$.  If any of the $Y_i$ belongs to $\Ea(D)$, then $\Yb_i=I_r$, so the factor $\Yb_i$ is not needed in the product $A=\Yb_1\cdots\Yb_k$.  After cancelling all such factors, we see that $A$ is a product of elements from $\Sib$.  Since $A\in\MrGr$ was arbitrary, we conclude that $\MrGr=\la\Sib\ra$.  In particular, $|\Si|\geq|\Sib|\geq\rank(\MrGr)=(q^r-1)/(q-1)$, giving (ii). \epf

\begin{rem}\label{rem:r=m4}
As in Remarks \ref{rem:r=m2} and \ref{rem:r=m3}, we deduce from the results of this section that for $|\F|=q<\infty$,
\bit
\item $\Ckl$ (and $\Rkl$) has $\sum_{s=0}^l q^{s(k-s)}\tqbin ls$ idempotents,
\item the semigroup generated by $E(\Ckl)$ (and the semigroup generated by $E(\Rkl)$) has rank and idempotent rank equal to $q^{l(k-l)} + (q^l-1)/(q-1)$.
\eit
\end{rem}

\section{Ideals}\label{sect:ideals}

In this final section, we consider the ideals of $P=\RegMmnJ$.  In particular, we show that each of the proper ideals is idempotent generated, and we calculate the rank and idempotent rank, showing that these are equal.  Although the next result is trivial if $r=0$ and well-known if $r=m=n$ (see Theorem \ref{thm_ideals_Mn}), the statement is valid for those parameters.

\ms
\begin{thm}\label{thm:ideals}
The ideals of $P=\RegMmnJ$ are precisely the sets 
\[
I_s^J=D_0^J\cup\cdots\cup D_s^J=\set{X\in P}{\rank(X)\leq s} \qquad\text{for $0\leq s\leq r$,}
\]
and they form a chain: $I_0^J\sub\cdots\sub I_r^J$.  If $0\leq s<r$, then $I_s^J = \la \Ea(D_s^J) \raa$ is generated by the idempotents in its top $\gDJ$-class, and if $|\F|=q<\infty$, then
\[
\rank(I_s^J)=\idrank(I_s^J)=q^{s(L-r)}\qbin rs, \qquad\text{where $L=\max(m,n)$.}
\]
\end{thm}

\pf 
For convenience, we will assume that $m\leq n$ throughout the proof, so that $L=n$.  (The other case will follow by duality.)

More generally, it may easily be checked that if the $\gJ$-classes of a semigroup $S$ form a chain, $J_0<\cdots<J_k$, then the ideals of $S$ are precisely the sets $I_h=J_0\cup\cdots\cup J_h$ for $0\leq h\leq k$ (and these obviously form a chain).  Now suppose $0\leq s<r$, let $\Ga\sub E(D_s(\M_r))$ be any idempotent generating set of $I_s(\M_r)$ (see Theorem~\ref{thm_ideals_Mn}), and put $\Om_\Ga=\set{[M,A,N]}{M\in\M_{m-r,r},\ A\in\Ga,\ N\in\M_{r,n-r}}$.  If $X=[M,A,N]\in I_s^J$ is arbitrary, then $A=B_1\cdots B_k$ for some $B_1,\ldots,B_k\in\Ga$, and it follows that $X=[M,B_1,N]\star\cdots\star[M,B_k,N]\in\la\Om_\Ga\raJ$.  Since $\Om_\Ga\sub\EJ(D_s^J)$, it follows that $I_s^J = \la \Ea(D_s^J) \raa$.


We now prove the statement about rank and idempotent rank.  Suppose $\Om$ is an arbitrary generating set for~$I_s^J$ where $0\leq s<r$.  Let $X\in D_s^J$ and consider an expression $X=Y_1\star\cdots\star Y_k$ with $Y_1,\ldots,Y_k\in\Om$.  Since $X=X\star Z\star X$ for some $Z\in D_s^J$, we may assume that $k\geq2$.  Since $I_{s-1}^J$ is an ideal of $I_s^J$ (we interpret $I_{s-1}^J=\emptyset$ if $s=0$), each of $Y_1,\ldots,Y_k$ must belong to $D_s^J=D_X^J$.  In particular, $Y_k\gDJ X=(Y_1\star\cdots\star Y_{k-1})\star Y_k$.  By stability, it then follows that $Y_k\gLJ (Y_1\star\cdots\star Y_{k-1})\star Y_k=X$.  Since $X\in D_s^J$ was arbitrary, it follows that $\Om$ contains at least one element from each $\gLJ$-class contained in $D_s^J$, and there are $q^{s(n-r)}\tqbin rs$ such $\gLJ$-classes, by Theorem~\ref{inflation_thm}(vi).  It follows that $\rank(I_s^J)\geq q^{s(n-r)}\tqbin rs 
=q^{s(L-r)}\tqbin rs$.

Since $\idrank(S)\geq\rank(S)$ for any idempotent generated semigroup $S$, the proof will be complete if we can find an idempotent generating set of $I_s^J$ of the specified size.  First, let $\Ga\sub E(D_s(\M_r))$ be such that ${\la\Ga\ra=I_s(\M_r)}$ and $|\Ga|=\tqbin rs$.  Fix some $A\in\Ga$, and let $\sim_A$ and $\approx_A$ be the equivalence relations on $\M_{m-r,r}$ and $\M_{r,n-r}$ defined before Theorem \ref{inflation_thm}, and let $\sM_A$ and $\sN_A$ be cross-sections of the equivalence classes of $\sim_A$ and $\approx_A$.  Let $\sM_A=\{M_1,\ldots,M_{q^{s(m-r)}}\}$ and $\sN_A=\{N_1,\ldots,N_{q^{s(n-r)}}\}$.  (We know $\sM_A$ and $\sN_A$ have the specified sizes by Theorem \ref{inflation_thm}.)  Put $Q=q^{s(n-r)}=q^{s(L-r)}$.  (Recall that we are assuming $m\leq n$.)  Extend $\sM_A$ arbitrarily to $\sM_A'=\{M_1,\ldots,M_Q\}$.  Now put $\Om_A=\set{[M_i,A,N_i]}{1\leq i\leq Q}$.  If $M\in\M_{m-r,r}$ and $N\in\M_{r,n-r}$ are arbitrary, then $M\sim M_i$ and $N\sim N_j$ for some $i,j$, and we have $[M,A,N]=[M_i,A,N_j]=[M_i,A,N_i]\star[M_j,A,N_j]\in\la\Om_A\raJ$.  Now put $\Om=\bigcup_{A\in\Ga}\Om_A$.  By the previous discussion, we see that $\la\Om\raJ$ contains $\Om_\Ga$, which is a generating set for $I_s^J$ (by the first paragraph of this proof), so $I_s^J=\la\Om\raJ$. Since $|\Om|=Q|\Ga|=q^{s(L-r)}\tqbin rs$, the proof is complete.  \epf


\begin{rem}
Again, we may deduce a corresponding statement for the ideals of the matrix semigroups $\Reg(\Ckl)$ and $\Reg(\Rkl)$; the reader may supply the details if they wish.
\end{rem}

\subsection*{Acknowledgements}

The first named author gratefully acknowledges the support of Grant No.~174019 of the Ministry of Education, Science, and Technological Development of the Republic of Serbia, and Grant No.~1136/2014 of the Secretariat of Science and Technological Development of the Autonomous Province of Vojvodina.  
The authors wish to thank Dr Attila Egri-Nagy for constructing the GAP \cite{GAP} code that enabled us to produce the eggbox diagrams from Figures \ref{fig:V3212_V3322}, \ref{fig:V2322_V2422}, \ref{fig:V2202} and \ref{fig:R}.


%
%
%
%

\footnotesize
\def\bibspacing{-1.1pt}
\bibliography{Sandwiches_bib}
\bibliographystyle{plain}
\end{document}